\documentclass[aap,MSNbibl,seceqn,citesort,dvips]{arximspdf}

\doi{10.1214/11-AAP781}
\volume{22}
\issue{2}
\pubyear{2012}
\firstpage{608}
\lastpage{669}

\makeatletter
\renewcommand{\epsilon}{\varepsilon}
\newcommand{\esssup}{\operatorname{ess\sup}}
\newcommand{\varliminf}{\operatorname{\underline{lim}}\limits}
\newproclaim{df}{Definition}[section]
\newtheorem{lem}{Lemma}[section]
\newtheorem{theo}{Theorem}[section]
\newtheorem{cor}{Corollary}[section]
\newtheorem{prop}{Proposition}[section]
\newproclaim{rem}{Remark}[section]
\newproclaim{example}{Example}
\newcommand{\fracd}[2]{({#1}/{#2})}
\newcommand{\fraca}[2]{{#1}/{#2}}
\makeatother

\begin{document}
\begin{frontmatter}

\title{Downside risk minimization via a large deviations~approach\thanksref{TT1}}
\runtitle{Downside risk minimization}

\begin{aug}
\author[A]{\fnms{Hideo} \snm{Nagai}\corref{}\ead[label=e1]{nagai@sigmath.es.osaka-u.ac.jp}}
\runauthor{H. Nagai}
\affiliation{Osaka University}
\address[A]{Division of Mathematical Science for\\ \quad  Social Systems\\
Graduate School of Engineering Science\\
Osaka University\\
Toyonaka, 560-8531\\ Japan\\
\printead{e1}} 
\end{aug}
\thankstext{TT1}{Supported in part by a Grant-in-aid for Scientific
Research 20340019 JSPS.}

\received{\smonth{9} \syear{2009}}
\revised{\smonth{4} \syear{2011}}

%
\begin{abstract}
We consider minimizing the probability of falling below a target growth
rate of the wealth process up to a time horizon $T$ in
an incomplete market model, and then study the asymptotic behavior of
minimizing probability as $T\to\infty$.
This problem can be closely related to an ergodic risk-sensitive
stochastic control problem in the risk-averse case. Indeed, in our main
theorem, we relate the former problem concerning the asymptotics for
risk minimization to the latter as its dual. As a result, we obtain an
expression of the limit value of the probability as the Legendre
transform of the value of the control problem, which is characterized
as the solution to an H-J-B equation of ergodic type, in the case of a
Markovian incomplete market model.
\end{abstract}

%
\end{frontmatter}

\section{Introduction}\label{sec1}
Risk management is a main topic in the study of finance. In the present
paper, we consider the problem of minimizing the downside risk
associated with an investor's total wealth in a certain incomplete
market model. More precisely, let $S_t^0$ be the price of a riskless
asset with the dynamics $dS^0_t=r_tS^0_t\,dt$, $(S_t^1,\ldots, S_t^m)$
the prices of the risky assets, and $N_t^i$, $ i=0,\ldots, m$, the
number of shares of $i$th security. Then the total wealth that the
investor possesses is defined as
\[
V_t=\sum_{i=0}^m N_t^iS_t^i,
\]
and we assume a self-financing condition,
\[
dV_t=\sum_{i=0}^m N_t^i\,dS_t^i.
\]
When setting the proportion of the portfolio invested in the $i$th
security as $h_t^i=\frac{N_t^i S_t^i}{V_t}$, we have
\[
\frac{dV_t}{V_t}=\sum_{i=0}^m h_t^i\,\frac{dS_t^i}{S_t^i},
\]
and the total wealth is denoted by $V_t=V_t(h)$, which is the solution
to this stochastic differential equation for a given strategy $h_t$.
Let us consider minimizing the probability
%
\begin{equation}\label{eq1.1}
P\biggl(\frac{1}{T}\log\frac{V_T(h)}{S_T^0}\leq\kappa\biggr)
\end{equation}
for a given target growth rate $\kappa$ by selecting portfolio choice
$h$. Let us make clear the meaning of the probability. If we choose
strategy $(h_t^0,h_t^1,\ldots, h_t^m)=(1,0,\ldots,0)$, then we have
\[
d\log V_t=\frac{dV_t}{V_t}=\frac{dS^0_t}{S_t^0}=d \log S_t^0.
\]
Thus, the probability is always $1$ for large time $T$ and $\kappa>0$.
Accordingly, in considering the above minimization, we investigate the
extent for which we can improve the probability by selecting a
strategy, as compared with the trivial strategy of investing the total
wealth in a riskless asset. The latter strategy is considered the
benchmark in terms of finance.

We shall consider the asymptotic behavior of the probability
%
\begin{equation}\label{eq1.2}
\varliminf_{T\to\infty}\frac{1}{T}\inf_{h_{\cdot}}\log P\biggl(\frac
{1}{T}\log\frac{V_T(h)}{S_T^0}\leq\kappa\biggr).
\end{equation}
According to the theory of large deviation, it is natural to relate
(\ref{eq1.2}) to
%
\begin{equation}\label{eq1.3}
{\hat\chi}(\gamma):=\varliminf_{T\to\infty}\frac{1}{T}\inf
_{h_{\cdot}}\log E\bigl[e^{\gamma\log(\fraca{V_T(h)}{S_T^0})}\bigr]
\end{equation}
for $\gamma<0$. Namely, as $T\to\infty$,
\[
\frac{1}{T}\inf_{h_{\cdot}}\log P \biggl (\frac{1}{T}\log\frac
{V_T(h)}{S_T^0}\in(-\infty,\kappa]  \biggr)  \to  -\inf_{k\in
(-\infty,\kappa]}\sup_{\gamma<0}\{\gamma k-{\hat\chi}(\gamma)\}
\]
is expected to hold since the Legendre transform $I(k)$ of ${\hat\chi
}(\gamma)$,
\[
I(k)=\sup_{\gamma<0}\{\gamma k-{\hat\chi}(\gamma)\},
\]
is regarded as the rate function of the asymptotics, if ${\hat\chi
}(\gamma)$ is a convex function; cf. \cite{DZ}. Note that we can see
from H\"older's inequality that\break $\log E[(\frac{V_T(h)}{S_T^0})^{\gamma
}]=\log E[e^{\gamma\log(\fraca{V_T(h)}{S_T^0})}]$ is a convex\vadjust{\goodbreak}
function of $\gamma$, but this does not always imply the convexity of
its infimum
%
\begin{equation}\label{eq1.4}
\inf_{h_{\cdot}}\log E\biggl[\biggl(\frac{V_T(h)}{S_T^0}\biggr)^{\gamma}\biggr] =\inf_{h_{\cdot}}\log
E\bigl[e^{\gamma\log(\fraca{ V_T(h)}{S_T^0})}\bigr].
\end{equation}
Therefore, the convexity of ${\hat\chi}(\gamma)$ cannot be
determined immediately and the above idea does not directly apply. In
the present paper, we will find the convexity of ${\hat\chi}(\gamma
)$ by identifying the solution of the H-J-B equation of ergodic type
with the limit value (\ref{eq1.3}); cf. Proposition~\ref{prop4.2} and Corollary~\ref{cor4.1}. Then
we shall see that the duality relation between (\ref{eq1.2}) and (\ref{eq1.3}) holds
under suitable conditions, as expected; cf. Theorem~\ref{theo2.4}.

Minimization (\ref{eq1.4}), which is equivalent to power utility maximization,
could be regarded as a risk-sensitive control problem. The infinite
time horizon counterpart of (\ref{eq1.4}) without a benchmark,
%
\begin{equation}\label{eq1.5}
\inf_{h_{\cdot}}\varliminf_{T\to\infty}\frac{1}{T}\log E\bigl[e^{\gamma\log
V_T(h)}\bigr],
\end{equation}
has been extensively studied as risk-sensitive control (e.g., \cite
{BP1,BPS2,FS1,FS2,FS3,KN,N1,NP,HS1,HS2}), and a benchmarked case has
recently been reported in \cite{DL}. From the viewpoint of stochastic
control theory, it may appear more natural, compared with the above
relationship between (\ref{eq1.2}) and (\ref{eq1.3}), to relate
%
\begin{equation}\label{eq1.6}
\inf_{h_{\cdot}}\varliminf_{T\to\infty}\frac{1}{T}\log E\bigl[e^{\gamma\log
(\fraca{V_T(h)}{S_T^0})}\bigr]
\end{equation}
to
%
\begin{equation}\label{eq1.7}
\inf_{h_{\cdot}} \varliminf_{T\to\infty}\frac{1}{T}\log P\biggl(\frac
{1}{T}\log\frac{V_T(h)}{S_T^0}\leq\kappa\biggr),
\end{equation}
which is considered in the present paper as well; cf. Theorem~\ref{theo2.5}.

We note that the problem relating (\ref{eq1.2}) to (\ref{eq1.3}) is thought to be
equivalent to considering
\[
\varliminf_{T\to\infty}\frac{1}{T}\inf_{h_{\cdot}}\log P\biggl(\frac
{1}{T}\log V_T(h)\leq\kappa\biggr)
\]
and
\[
{\check\chi}(\gamma):=\varliminf_{T\to\infty}\frac{1}{T}\inf
_{h_{\cdot}}\log E\bigl[e^{\gamma\log V_T(h)}\bigr]
\]
without a benchmark. However, the arguments in this article may be
simpler than in the case without a benchmark (cf. Remark~\ref{rem2.2}).

In previous papers \cite{HNS,N2}, we studied similar
asymptotic behavior without benchmarks for linear Gaussian models in
relation to the asymptotics of risk-sensitive portfolio optimization.
Indeed, we established a duality relation between these problems, and
as a result, an explicit expression of the limit value of the
probability minimizing downside risk for each case of full and partial
information. To obtain those results, the key analysis involved Poisson
equations derived by taking derivatives with respect to $\gamma$ of
the H-J-B equations of ergodic type corresponding to risk-sensitive
control over an infinite time horizon. Since the solutions of the H-J-B
equations can be explicitly expressed as quadratic functions by using
the solutions of the Riccati equations for linear Gaussian models,
analysis of the differentiability of the solutions of the Riccati
equations with respect to $\gamma$ was essential in those works.

In the present paper, we shall consider general diffusion market models
and discuss the above-mentioned duality relation between the
asymptotics of the minimization of downside risk and the risk-sensitive
stochastic control for large time. Since the solutions of H-J-B
equations of ergodic type do not always have explicit expressions, we
need to consider, in general, the differentiablity with respect to
$\gamma$ of H-J-B equations of ergodic type. The analysis is presented
in Sections~\ref{sec5} and~\ref{sec6} based on the results concerning H-J-B equations of
ergodic type and related stochastic control problems given in Sections~\ref{sec3} and~\ref{sec4}. Here, we mention the ongoing work of Hata and Sheu \cite
{HSh}, which is closely related to the present paper and examines
similar problems under the assumptions that $\alpha(x)$ in (\ref{eq2.2}) in
Section~\ref{sec2} is bounded and that $\beta(x)^*x\leq-c|x|^2$ for $ |x|\geq
R$ in place of (\ref{eq2.20}). We shall explain more precisely the
relationships between these papers in Remark~\ref{rem2.3}.

We note that maximization of an upside chance probability for the long
term was studied by Pham \cite{P1,P2} for continuous time
models, and then by Stettner \cite{Stet} for discrete time models, in
relation to risk-sensitive portfolio optimization in the risk-seeking
case. By regarding the maximization problem as large deviation control,
Pham established a duality relation between these two types of
problems. Explicit calculation of the limit value is given in the case
of 1-dimensional Gaussian models. The problem was later extended to a
nonlinear case by Hata and Sekine \cite{HS1,HS2} and also to
the partial information case by Hata and Iida \cite{HI} for
1-dimensional Gaussian models. See also related works \cite{Br1,Br2,St}. However, asymptotic estimates of downside risk
probabilities and upside chance probabilities cannot be obtained in
parallel. Indeed, obtaining the estimates of downside risk is rather
difficult than those of upside chance and further analysis of H-J-B
equations is required to show the estimates as was shown in \cite
{HNS}. Further note that large deviations control (\ref{eq1.7}) is an
unconventional optimization problem, and thus we need to employ a new
approach to study it.

\section{Setting up and main results}\label{sec2}
Consider a market model with $m+1$ securities and $n$ factors, where
the bond price is governed by the ordinary differential equation
%
\begin{equation}\label{eq2.1}
dS^0(t)=r(X_t)S^0(t)\,dt,\qquad S^0(0)=s^0.
\end{equation}
The other security prices and factors are assumed to satisfy the
stochastic differential equations
%
\begin{eqnarray}\label{eq2.2}
 dS^i(t)&=&S^i(t)\Biggl\{\alpha^i(X_t)\,dt+\sum_{k=1}^{n+m}\sigma^i_k(X_t)\,dW^k_t\Biggr\}, \nonumber
\\[-8pt]
\\[-8pt]
  S^i(0)&=&s^i, \qquad   i=1,\ldots,m,
\nonumber
\end{eqnarray}
and
%
\begin{eqnarray}\label{eq2.3}
 dX_t&=&\beta(X_t)\,dt +\lambda(X_t)\,dW_t, \nonumber
 \\[-8pt]
 \\[-8pt]
  X(0)&=&x,
\nonumber
\end{eqnarray}
where $W_t=(W_t^k)_{k=1,\ldots,(n+m)}$ is a standard $m+n$-dimensional
Brownian motion process on a probability space $(\Omega,{\mathcal
F},P)$. Let $N_t^i$ be the number of the shares of the $i$th security.
Then, the total wealth that the investor possesses is defined as
\[
V_t=\sum_{i=0}^m N_t^iS_t^i
\]
and the proportion of the portfolio invested in the $i$th security is
\[
h_t^i=\frac{N_t^iS_t^i}{V_t},\qquad i=0,1,2,\ldots,m.
\]
$N_t=(N_t^0,N_t^1,N_t^2,\ldots,N_t^m)$ [or, $h_t=(h_t^1,\ldots
,h_t^m)$] is called self-financing if
\[
dV_t=\sum_{i=0}^mN_t^i\,dS_t^i=\sum_{i=0}^m\frac{V_th_t^i}{S_t^i}\,dS_t^i.
\]
Thus, under the self-financing condition, we have
\begin{eqnarray*}
 \frac{dV_t}{V_t}&=&h_t^0r(X_t)\,dt+\sum_{i=1}^mh_t^i\Biggl\{\alpha
^i(X_t)\,dt+\sum_{j=1}^{n+m}\sigma_j^i(X_t)\,dW_t^j\Biggr\} \\
  &=&r(X_t)\,dt+\sum_{i=1}^mh_t^i\Biggl\{\bigl(\alpha^i(X_t)-r(X_t)\bigr)\,dt+\sum
_{j=1}^{n+m}\sigma_j^i(X_t)\,dW_t^j\Biggr\}.
\end{eqnarray*}
Here we note that $h_t$ is defined as an $m$-vector consisting of
$h_t^1,\ldots,h_t^m$ since $h_t^0=1-\sum_{i=1}^mh_t^i$ holds by definition.

The filtration that must be satisfied by admissible investment strategies
\[
{\mathcal G}_t=\sigma\bigl(S(u),X(u),  u\leq t\bigr)
\]
is relevant in the present problem, and we introduce the following
definition. 
%
\begin{df}$h(t)_{0\leq t\leq T}$ is said to be an investment
strategy if
$h(t)$ is an $R^m$ valued ${\mathcal G}_t$-progressively measurable
stochastic process such that
\[
P\biggl(  \int_0^T|h(s)|^2\,ds<\infty\biggr)=1.
\]
\end{df}

The set of all investment strategies is denoted by
${\mathcal H}(T)$.
For a given $h\in{\mathcal H}(T)$, the process $V_t=V_t(h)$
representing the total wealth of the investor at time $t$ is determined
by the stochastic differential equation shown above.
%
\begin{eqnarray}\label{eq2.5}
\frac{dV_t}{V_t}&=&
r(X_t)\,dt+h(t)^*\bigl(\alpha(X_t)-r(X_t){\mathbf1}\bigr)\,dt+h(t)^*\sigma(X_t)\,dW_t,\nonumber
\\[-8pt]
\\[-8pt]
 V_0&=&v_0,
\nonumber
\end{eqnarray}
where ${\mathbf1}=(1,1,\ldots,1)^*$.

We are interested in the asymptotics of minimization of a downside risk
for a given constant $\kappa$ in comparison with investing the whole
portfolio in a~riskless security as the benchmark
%
\begin{equation}\label{eq2.6}
J(\kappa):=\varliminf_{T\to\infty}\frac{1}{T}\inf_{h\in{\mathcal
H}(T)}\log P\biggl(\frac{1}{T}\log\frac{ V_T(h)}{ S_T^0}\leq\kappa
\biggr).
\end{equation}
We also examine downside risk minimization with the benchmark $S^0$
over an infinite time horizon,
%
\begin{equation}\label{eq2.7}
J_{\infty}(\kappa):=\inf_{h\in{\mathcal H}}\varliminf_{T\to\infty
}\frac{1}{T}\log P\biggl(\frac{1}{T}\log\frac{ V_T(h)}{ S_T^0}\leq\kappa
\biggr),
\end{equation}
where
\[
{\mathcal H}=\{h; h\in{\mathcal H}(T), \forall T\}.
\]
$J(\kappa)$ will be shown to be related to the following
risk-sensitive asset allocation problem with benchmark $S^0$. Namely,
for a given constant $\gamma< 0$, let us consider the asymptotics
%
\begin{equation}\label{eq2.8}
{\hat\chi}(\gamma)=\varliminf_{T\to\infty}\frac{1}{T}\inf_{h
\in{\mathcal A}(T)} J(v,x;h;T),
\end{equation}
where
%
\begin{equation}\label{eq2.9}
J(v,x;h;T)=\log E\biggl[\biggl(\frac{V_T(h)}{S_T^0}\biggr)^{\gamma}\biggr]=\log E\bigl[e^{\gamma
\log(\fraca{ V_T(h)}{S_T^0})}\bigr] ,
\end{equation}
and $h$ ranges over the set ${\mathcal A}(T)$ of all admissible
investment strategies defined by
\[
{\mathcal A}(T)=\bigl\{h\in{\mathcal H}(T); E\bigl[e^{\gamma\int
_0^Th_s^*\sigma(X_s)\,dW_s-\fracd{\gamma^2}{2}\int_0^Th_s^*\sigma
\sigma^*(X_s)h_sds}\bigr]=1\bigr\}.
\]
Then we shall see that (\ref{eq2.6}) could be considered the dual problem to~(\ref{eq2.8}),
while (\ref{eq2.7}) is the dual problem to risk-sensitive asset
allocation over an infinite time horizon,
%
\begin{equation}\label{eq2.10}
\chi_{\infty}(\gamma)=\inf_{h \in{\mathcal A}}\varliminf_{T\to
\infty}\frac{1}{T} J(v,x;h;T),
\end{equation}
where
\[
{\mathcal A}=\{h; h\in{\mathcal A}(T);\forall T\}.
\]
We shall consider these problems under the assumptions that
%
\begin{equation}\label{eq2.11}
\alpha \mbox{ and } \beta \mbox{ are globally Lipschitz, }  \lambda
\in C^2_b, \sigma,r\in C^1_b, \alpha, \beta\in C^1
\end{equation}
and
%
\begin{equation}\label{eq2.12}
\cases{\displaystyle
 c_1|\xi|^2\leq\xi^*\lambda\lambda^*(x)\xi\leq c_2|\xi|^2, &\quad
$c_1, c_2>0,   \xi\in R^n$, \vspace*{2pt}\cr\displaystyle
  c_1|\zeta|^2\leq\zeta^*\sigma\sigma^*(x)\zeta\leq c_2|\zeta
|^2,  & \quad    $\zeta\in R^m$,
}
\end{equation}
hold. In considering these problems, we first introduce the value function
%
\begin{equation}\label{eq2.13}
v(t,x)=\inf_{h_{\cdot}\in{\mathcal A}(T-t)}\log E\bigl[e^{\gamma\log(\fraca
{V_{T-t}(h)}{S_{T-t}^0})}\bigr].
\end{equation}
Note that
\[
e^{\gamma\log V_T}=v_0^{\gamma}e^{\gamma\int_0^T\{r(X_s)+h_s^*{\hat
\alpha}(X_s)-\fracd{1}{2}h_s^*\sigma\sigma^*(X_s)h_s\}\,ds+\gamma\int
_0^Th_s^*\sigma(X_s)\,dW_s},
\]
where ${\hat\alpha}(x)=\alpha(x)-r(x){\mathbf1}$. Therefore
\[
e^{\gamma( \log V_T-\log S_T^0)}=v_0^{\gamma}e^{\gamma\int_0^T\eta
(X_s,h_s)\,ds+\gamma\int_0^Th_s^*\sigma(X_s)\,dW_s-\fracd{\gamma
^2}{2}\int_0^Th_s^*\sigma\sigma^*(X_s)h_sds},
\]
where
\[
\eta(x,h)=h^*{\hat\alpha}(x)-\frac{1-\gamma}{2}h^*\sigma\sigma^*(x)h.
\]
By introducing a probability measure
\[
P^h(A)=E\bigl[e^{\gamma\int_0^Th_s^*\sigma(X_s)\,dW_s-\fracd{\gamma
^2}{2}\int_0^Th_s^*\sigma\sigma^*(X_s)h_sds}\dvtx A\bigr],
\]
the dynamics of the factor process can be written as
\[
dX_t=\{\beta(X_t)+\gamma\lambda\sigma^*(X_t)h_t\}\,dt +\lambda
(X_t)\,dW_t^h,\qquad X_0=x,
\]
with the new Brownian motion process $W_t^h$ defined by
\[
W_t^h:=W_t-\gamma\int_0^t\sigma^*(X_s)h_s\,ds,
\]
and the value function written as
%
\begin{equation}\label{eq2.14}
v(t,x)=\gamma\log v_0+\inf_{h_{\cdot}\in{\mathcal A}(T)}\log E^h\bigl[e^{\gamma
\int_0^{T-t}\eta(X_s,h_s)\,ds}\bigr].
\end{equation}
The H-J-B equation for the value function $v(t,x)$ is
\[
\cases{\displaystyle
\frac{\partial v}{\partial t}+\frac{1}{2}\operatorname{tr}[\lambda\lambda
^*D^2v]+\frac{1}{2}(Dv)^*\lambda\lambda^*Dv\vspace*{1pt}\cr\displaystyle
 \qquad {}+\inf_h\{[\beta+\gamma
\lambda\sigma^*h]^*Dv+\gamma\eta(x,h)\}=0,  \vspace*{2pt}\cr\displaystyle
v(T,x)=\gamma\log v_0,
}
\]
which is also written as
%
\begin{equation}\label{eq2.15}
\cases{\displaystyle
\frac{\partial v}{\partial t}+\frac{1}{2}\operatorname{tr}[\lambda\lambda
^*D^2v]+\beta_{\gamma}^*Dv+\frac{1}{2}(Dv)^*\lambda N_{\gamma
}^{-1}\lambda^*Dv-U_{\gamma}=0,  \vspace*{2pt}\cr\displaystyle
v(T,x)=\gamma\log v_0,
}
\end{equation}
where
\begin{eqnarray*}
\beta_{\gamma}&=&\beta+\frac{\gamma}{1-\gamma}\lambda\sigma
^*(\sigma\sigma^*)^{-1}{\hat\alpha}, \\
N_{\gamma}^{-1}&=&I+\frac{\gamma}{1-\gamma}\sigma^*(\sigma\sigma
^*)^{-1}\sigma
\end{eqnarray*}
and
\[
U_{\gamma}=-\frac{\gamma}{2(1-\gamma)}{\hat\alpha}^*(\sigma
\sigma^*)^{-1}{\hat\alpha}.
\]

\begin{rem}
\begin{eqnarray*}
&&\inf_{h\in R^m}\{[\gamma\lambda\sigma^*h]^*Dv+\gamma\eta(x,h)\}
 \\
&& \qquad =\inf_{h\in R^m}\biggl\{[\gamma\lambda\sigma^*h]^*Dv-\frac{\gamma
(1-\gamma)}{2}h^*\sigma\sigma^*h+\gamma h^*{\hat\alpha}\biggr\} \\
&& \qquad =\inf_{h\in R^m}\biggl\{-\frac{\gamma(1-\gamma)}{2}\biggl[h-\frac
{1}{1-\gamma}(\sigma\sigma^*)^{-1}({\hat\alpha}+\sigma\lambda
^*Dv)\biggr]^*\\
&& \qquad\hphantom{=\inf_{h\in R^m}\biggl\{} {}\times\sigma\sigma^*\biggl[h-\frac{1}{1-\gamma}(\sigma\sigma
^*)^{-1}({\hat\alpha}+\sigma\lambda^*Dv)\biggr] \\
&& \qquad\hphantom{=\inf_{h\in R^m}\biggl\{} {}+\frac{\gamma}{2(1-\gamma)}({\hat\alpha}+\sigma\lambda
^*Dv)^*(\sigma\sigma^*)^{-1}({\hat\alpha}+\sigma\lambda^*Dv)\biggr\}.
\end{eqnarray*}
Therefore, the function
\[
{\hat h}(t,x):=\frac{1}{1-\gamma}(\sigma\sigma^*)^{-1}({\hat\alpha
}+\sigma\lambda^*Dv)
\]
defines the generator of the optimal diffusion ${\hat L}$ for $
\inf_{h \in{\mathcal A}(T)} J(v,x;h;T)$:
\[
{\hat L}\psi:=\frac{1}{2}\operatorname{tr}[\lambda\lambda^*D^2\psi
]+\biggl[\beta+\frac{\gamma}{1-\gamma}\lambda\sigma^*(\sigma\sigma
^*)^{-1}({\hat\alpha}+\sigma\lambda^*Dv)\biggr]^*D\psi,
\]
which is seen in Proposition~\ref{prop2.1}.
\end{rem}

Set ${\bar v}=-v.$ Then,
%
\begin{equation}\label{eq2.16}\qquad
\cases{\displaystyle
\frac{\partial{\bar v}}{\partial t}+\frac{1}{2}\operatorname{tr}[\lambda
\lambda^*D^2{\bar v}]+\beta_{\gamma}^*D{\bar v}-\frac{1}{2}(D{\bar
v})^*\lambda N_{\gamma}^{-1}\lambda^*D{\bar v}+U_{\gamma}=0, \vspace*{2pt}\cr\displaystyle
{\bar v}(T,x)=-\gamma\log v_0.
}
\end{equation}
Since $I-\sigma^*(\sigma\sigma^*)^{-1}\sigma\geq0$, which is
easily seen by taking $\xi=\sigma^*\zeta+\mu, $ with $\mu$
orthogonal to the range of $\sigma^*$, and seeing that $\xi
^*(I-\sigma^*(\sigma\sigma^*)^{-1}\sigma)\xi=\mu^*\mu$, we have
%
\begin{equation}\label{eq2.17}
\frac{1}{1-\gamma}I\leq N^{-1}\leq I.
\end{equation}
As for existence of a solution to (\ref{eq2.16}) satisfying sufficient
regularities, we have the following results; cf. \cite{BFN,N1}.\vspace*{-3pt}
%
\begin{theo}[(\cite{BFN,N1})]\label{theo2.1} Assume (\ref{eq2.11}) and (\ref{eq2.12}). Then
H-J-B equation~(\ref{eq2.16}) has a solution such that
\begin{eqnarray*}
{\bar v}(t,x)+\gamma\log v_0&\geq&0,
\\[-2pt]
{\bar v}, \frac{\partial{\bar v}}{\partial t}, \frac{\partial
{\bar v}}{\partial x_k}, \frac{\partial^2 {\bar v}}{\partial
x_k\,\partial x_j}&\in& L^p(0,T;L^p_{\mathrm{loc}}(R^n)), \qquad  1<\forall p<\infty,
\\[-2pt]
\frac{\partial^2{\bar v}}{\partial t^2}, \frac{\partial^2{\bar
v}}{\partial x_k\,\partial t}, \frac{\partial^3 {\bar v}}{\partial
x_k\,\partial x_j\,\partial x_l},
\frac{\partial^3 {\bar v}}{\partial x_k\,\partial x_j\,\partial t}&\in&
L^p(0,T;L^p_{\mathrm{loc}}(R^n)),
\\[-2pt]
\frac{\partial{\bar v}}{\partial t}&\leq&0,
\end{eqnarray*}
and
\begin{eqnarray*}
 |\nabla{\bar v}|^2-k_0\frac{\partial{\bar v}}{\partial t}&\leq& C\bigl(
|\nabla Q_{\gamma}|_{2\rho}^2+|Q_{\gamma}|_{2\rho}^2+|\nabla
(\lambda\lambda^*)|^2_{2\rho} \\[-2pt]
&&\hphantom{C\bigl(}  {} +|\nabla\beta_{\gamma}|_{2\rho}+|\beta_{\gamma
}|_{2\rho}^2+|U_{\gamma}|_{2\rho}+|\nabla U_{\gamma}|_{2\rho}+1\bigr)
\end{eqnarray*}
for $x\in B_{\rho}$ and $t\in[0,T)$, where $Q_{\gamma}=\lambda
N_{\gamma}^{-1}\lambda^*$, $k_0=\frac{4(1+c)(1-\gamma)}{-\gamma},
 c>0$, $|f|_{2\rho}=\sup_{\{x;x\in B_{2\rho}\}}|f(x)|$, $C$ is a
universal constant and $B_{\rho}=\{x\in R^n; |x|<\rho\}$.\vspace*{-3pt}
\end{theo}

For ${\hat h}(t,x)$, we consider the stochastic differential equation
\[
dX_t=\{\beta(X_t)+\gamma\lambda\sigma^*(X_t){\hat h}(t,X_t)\}\,dt
+\lambda(X_t)\,dW_t^{\hat h},\qquad X_0=x,
\]
and define ${\hat h}_t:={\hat h}(t,X_t)$ for the solution $X_t$ of the
stochastic differential equation. Note that the solution of this
stochastic differential equation is obtained by the change of measure
from the solution of (\ref{eq2.3}). Indeed, we can see that $\nabla v$ has at
most linear growth under assumptions (\ref{eq2.11}) and~(\ref{eq2.12}) from the above
gradient estimates, and therefore,
\[
E\bigl[e^{\gamma\int_0^T{\hat h}(s,X_s)^*\sigma(X_s)\,dW_s-\fracd{\gamma
^2}{2}\int_0^T{\hat h}(s,X_s)^*\sigma\sigma^*(X_s){\hat h}(s,X_s)\,ds}\bigr]=1
\]
holds. Thus
\[
P^{\hat h}(A):=E\bigl[e^{\gamma\int_0^T{\hat h}(s,X_s)^*\sigma
(X_s)\,dW_s-\fracd{\gamma^2}{2}\int_0^T{\hat h}(s,X_s)^*\sigma\sigma
^*(X_s){\hat h}(s,X_s)\,ds};A\bigr]\vadjust{\goodbreak}
\]
defines a probability measure. Under this measure, $X_t$ turns out to
be a~solution of the above stochastic differential equation.

The following is a verification theorem, the proof of which is almost
the same as the proof of Proposition 2.1, \cite{N1}, and thus is
omitted here.
%
\begin{prop}[(\cite{N1})]\label{prop2.1}
Assume (\ref{eq2.11}) and (\ref{eq2.12}). Then ${\hat h}_t^{(\gamma,T)}\equiv{\hat
h}_t:={\hat h}(t, X_t)\in{\mathcal A}(T)$ and it is optimal
%
\begin{equation}\label{eq2.18}
 \qquad v(0,x)=\inf_{h_{\cdot}}\log E\bigl[e^{\gamma(\log V_T(h)-\log S_T^0)}\bigr]=\log
E\bigl[e^{\gamma(\log V_T({\hat h})-\log S_T^0)}\bigr].
\end{equation}
\end{prop}

Let us consider an H-J-B equation of ergodic type that is thought to be
the limit equation of (\ref{eq2.15}). Namely,
%
\begin{equation}\label{eq2.19}
\chi=\tfrac{1}{2}\operatorname{tr}[\lambda\lambda^*D^2 w]+\beta_{\gamma
}^*D w+\tfrac{1}{2}(D w)^*\lambda N_{\gamma}^{-1}\lambda^*D
w-U_{\gamma}.
\end{equation}
Set
\[
G(x):=\beta(x)-\lambda\sigma^*(\sigma\sigma^*)^{-1}{\hat\alpha}(x),
\]
and assume that
%
\begin{equation}\label{eq2.20}
G(x)^*x\leq-c_G|x|^2+c_G',\qquad c_G, c_G'>0,
\end{equation}
and
%
\begin{equation}\label{eq2.21}
{\hat\alpha}^*(\sigma\sigma^*)^{-1}{\hat\alpha}\to \infty
 \qquad \mbox{as }   |x| \to\infty.
\end{equation}
Under these assumptions, we have a solution to the H-J-B equation of
ergodic type, and the proof is given in Proposition~\ref{prop3.1} in Section~\ref{sec3}.
%
\begin{prop} Assume (\ref{eq2.11}), (\ref{eq2.12}), (\ref{eq2.20}) and (\ref{eq2.21}). Then (\ref{eq2.19}) has a
solution $(\chi, w^{(\gamma)})$ such that $ w\in C^2(R^n)$,
\[
w(x) \to-\infty   \qquad \mbox{as }   |x| \to\infty
\]
and the solution satisfying this condition is unique up to additive
constants with respect to $w$.
\end{prop}

We further assume that
%
\begin{equation}\label{eq2.22}
{\hat\alpha}^*(\sigma\sigma^*)^{-1}{\hat\alpha}\geq
c_0|x|^2-c_0',\qquad c_0,  c_0'>0.
\end{equation}
Then we have the following theorem, and the proof is given after
Proposition~\ref{prop4.2} in Section~\ref{sec4}.
%
\begin{theo}\label{theo2.2} Under assumptions (\ref{eq2.11}), (\ref{eq2.12}), (\ref{eq2.20}) and (\ref{eq2.22}), we have
\[
{\hat\chi}(\gamma)=\varliminf_{T\to\infty}\frac
{1}{T}v(0,x;T)=\chi(\gamma).
\]
\end{theo}

The following results are important to prove our main results, and the
proofs are given in Lemma~\ref{lem6.3}, Lemma~\ref{lem7.1} and
Corollary~\ref{cor4.1}.\vadjust{\goodbreak}
%
\begin{theo}\label{theo2.3} Let $(\chi,w^{(\gamma)})$ be a solution to (\ref{eq2.19}). Then
under the assumptions of Theorem~\ref{theo2.2}, $\chi(\gamma)$ and $w^{(\gamma
)}$ are differentiable with respect to $\gamma$ and $\chi(\gamma)$
is convex. Their derivatives satisfy
%
\begin{eqnarray}\label{eq2.23}\qquad
\chi'(\gamma)&=&\frac{1}{2}\operatorname{tr}[\lambda\lambda^*D^2
w_{\gamma}]+\bigl(\beta_{\gamma}^*+\bigl(D w^{(\gamma)}\bigr)^*\lambda N_{\gamma
}^{-1}\lambda^*\bigr)D w_{\gamma} \nonumber
\\[-8pt]
\\[-8pt]
&&{} +\frac{1}{2(1-\gamma)^2}\bigl\{{\hat\alpha}+\sigma
\lambda^*Dw^{(\gamma)}\bigr\}^*(\sigma\sigma^*)^{-1}\bigl\{{\hat\alpha
}+\sigma\lambda^*Dw^{(\gamma)}\bigr\},
\nonumber
\end{eqnarray}
where $w_{\gamma}=\frac{\partial w^{(\gamma)}}{\partial\gamma}$.
Furthermore,
\[
\lim_{\gamma\to-\infty}\chi'(\gamma)=0.
\]
\end{theo}

%
\begin{rem}\label{rem2.2} It is important to know the limit value $\lim_{\gamma
\to-\infty}\chi'(\gamma)$ since it determines the left endpoint of
the interval of the target growth rate~$\kappa$, which makes $J(\kappa
)$ finite. Here we compare the results above with those to be expected
for the case without a benchmark, considering asymptotics
\[
{\check\chi}(\gamma):=\varliminf_{T\to\infty}\frac{1}{T}\inf
_{h_{\cdot}}\log E\bigl[e^{\gamma\log V_T(h)}\bigr].
\]
The H-J-B equation of ergodic type of this problem becomes
\[
{\check\chi}=\tfrac{1}{2}\operatorname{tr}[\lambda\lambda^*D^2 {\check
w}]+\beta_{\gamma}^*D{\check w}+\tfrac{1}{2}(D {\check w})^*\lambda
N_{\gamma}^{-1}\lambda^*D {\check w}-U_{\gamma}+\gamma r(x),
\]
and we can obtain its derivative
\begin{eqnarray*}
{\check\chi}'(\gamma)&=&\frac{1}{2}\operatorname{tr}[\lambda\lambda^*D^2
{\check w}_{\gamma}]+\bigl(\beta_{\gamma}^*+(D {\check w})^*\lambda
N_{\gamma}^{-1}\lambda^*\bigr)D {\check w}_{\gamma} \\
&& {}+\frac{1}{2(1-\gamma)^2}\{{\hat\alpha}+\sigma
\lambda^*D{\check w}\}^*(\sigma\sigma^*)^{-1}\{{\hat\alpha}+\sigma
\lambda^*D{\check w}\}+r
\end{eqnarray*}
through almost the same arguments as the current ones provided to
obtain the results in the present article. The difference appears in
considering the asymptotics of $\chi'(\gamma)$ as $\gamma\to-\infty
$. Indeed,
\[
\lim_{\gamma\to-\infty}{\check\chi}'(\gamma)=\lim_{\gamma\to
-\infty}\int r(x){\check m}_{\gamma}(dx)<\infty
\]
could be seen as in \cite{HSh}, where ${\check m}_{\gamma}(dx)$ is
the invariant measure of ${\check L}$-diffusion process and ${\check
L}$ is defined by
\[
{\check L}\psi=\tfrac{1}{2}\operatorname{tr}[\lambda\lambda^*D^2\psi
]+(\beta_{\gamma}+\lambda N_{\gamma}^{-1}\lambda^*D{\check
w})^*D\psi.
\]
Note that ${\check L}$ corresponds to ${\bar L}$ defined by (\ref{eq4.16}) in
the present paper and
can be shown to be ergodic under suitable conditions in a manner
similar to the proof of Proposition~\ref{prop4.3}.
\end{rem}

Now we can state our main theorems. The proofs are given in
Sections~\ref{sec7}
and~\ref{sec8}.\vadjust{\goodbreak}
%
\begin{theo}\label{theo2.4}
Under the assumptions of Theorem~\ref{theo2.2}, for $0<\kappa<{\hat\chi}'(0-)$,
%
\begin{equation}\label{eq2.24}
J(\kappa)=-\inf_{k\in(-\infty,\kappa]}\sup_{\gamma<0}\{\gamma
k-{\hat\chi}(\gamma)\}=-\sup_{\gamma<0}\{\gamma\kappa-{\hat\chi
}(\gamma)\}.
\end{equation}
Moreover, for $\gamma(\kappa)$ such that ${\hat\chi}'(\gamma
(\kappa))=\kappa\in(0,{\hat\chi}'(0-))$, take a strategy ${\hat
h}_t^{(\gamma(\kappa),T)}$ defined in Proposition~\ref{prop2.1}. Then,
\[
J(\kappa)=\lim_{T\to\infty}\frac{1}{T}\log P\biggl(\frac{1}{T}\log
\frac{V_T({\hat h}^{(\gamma(\kappa),T)})}{S_T^0}\leq\kappa\biggr).
\]
For $\kappa<0$,
\[
J(\kappa)=-\sup_{\gamma<0}\{\gamma\kappa-{\hat\chi}(\gamma)\}
=-\infty.
\]
\end{theo}

For the solution $w=w^{(\gamma)}$ to H-J-B equation ergodic type
(\ref{eq2.19}), let us set
\[
{\bar h}(x)=\frac{1}{1-\gamma}(\sigma\sigma^*)^{-1}({\hat\alpha
}+\sigma\lambda^*Dw)(x).
\]
Further consider the stochastic differential equation
%
\begin{equation}\label{eq2.25}
 \qquad dX_t=\{\beta(X_t)+\gamma\lambda\sigma^*(X_t){\bar h}(X_t)\}\,dt
+\lambda(X_t)\,dW_t^{\bar h},\qquad X_0=x,
\end{equation}
and define ${\bar h}_t^{(\gamma(\kappa))}:={\bar h}(X_t)$ for the
solution $X_t$ of the stochastic differential equation. Then we have
the following theorem.
%
\begin{theo}\label{theo2.5}
Under the assumptions of Theorem~\ref{theo2.2}, let $0<\kappa<{\hat\chi}'(0-)$
and $\gamma(\kappa)$ be the same as above. We also assume that
%
\begin{equation}\label{eq2.26}\qquad
\bigl(Dw^{(\gamma)}\bigr)^*\lambda\sigma^*(\sigma\sigma^*)^{-1}\sigma
\lambda^*Dw^{(\gamma)}<{\hat\alpha}^*(\sigma\sigma^*)^{-1}{\hat
\alpha}, \qquad    \gamma=\gamma(\kappa).
\end{equation}
Then
\[
J_{\infty}(\kappa)=J(\kappa)=-\inf_{k\in(-\infty,\kappa]}\sup
_{\gamma<0}\{\gamma k-{\hat\chi}(\gamma)\}=-\sup_{\gamma<0}\{
\gamma\kappa-{\hat\chi}(\gamma)\}
\]
and
\[
J(\kappa)=\lim_{T\to\infty}\frac{1}{T}\log P\biggl(\log\frac
{V_T(h^{(\gamma(\kappa))})}{S_T^0}\leq\kappa T\biggr).
\]
\end{theo}

%
\begin{rem}\label{rem2.3}
In our previous paper \cite{HNS}, we studied similar problems without
benchmarks in the case of linear Gaussian models. Specifically, we
discussed the case where $\alpha(x)=Ax+a$, $\beta(x)=Bx+b$, $\sigma
(x)\equiv\sigma$, $\lambda(x)\equiv\lambda$ and $r(x)\equiv r$, in
which $A$, $B$, $\sigma$ and $\lambda$ (resp., $a$ and $b$) are all
constant matrices (resp., vectors), and $r$ is a constant. Under the
main assumption that
\[
G:=B-\lambda\sigma^*(\sigma\sigma^*)^{-1}A     \mbox{ is stable},\vadjust{\goodbreak}
\]
which corresponds to (\ref{eq2.20}) above, we obtained results similar to
Theorems~\ref{theo2.4} and~\ref{theo2.5}. The present paper is a natural extension to
general diffusion incomplete market models. On the other hand, Hata and
Sheu \cite{HSh} treat the case where $\alpha(x)$ is bounded, and
$\beta(x)^*x\leq-c|x|^2$ for $|x|\geq R$, in which linear Gaussian
models are excluded. In that case, $U_{\gamma}$ becomes bounded and
they employ quite different methods from ours to analyze H-J-B equation~(\ref{eq2.19}), while assumption (\ref{eq2.22}) is crucial in our settings. For that
reason our theorems do not include the case where $\alpha(x)$ is bounded.
\end{rem}

%
\begin{rem}
The generator of the optimal diffusion process governed by (\ref{eq2.25}) for
risk-sensitive control problem (\ref{eq2.10}) is defined by
\[
L_{\infty}\psi:=\frac{1}{2}\operatorname{tr}[\lambda\lambda^*D^2 \psi
]+\biggl[\beta_{\gamma}^*+\frac{\gamma}{1-\gamma}(Dw)^*\lambda\sigma
^*(\sigma\sigma^*)^{-1}\sigma\lambda^*\biggr]D\psi.
\]
On the other hand, in proving Theorem~\ref{theo2.2} we introduce another type of
stochastic control problem (\ref{eq4.9}) with (\ref{eq4.7}). The generator of the
optimal diffusion process for this problem is defined by (\ref{eq4.16}).
\[
{\bar L}\psi=\tfrac{1}{2}\operatorname{tr}[\lambda\lambda^*D^2 \psi
]+[\beta_{\gamma}^*+(Dw)^*\lambda N_{\gamma}^{-1}\lambda^*]D\psi,
\]
where $w$ is a solution to H-J-B equation (\ref{eq2.19}) of ergodic type. Then
we note that ${\bar L}$ is related to $L_{\infty}$ through the gauge
transform,
\[
[e^{-w}L_{\infty}e^w]\varphi=\bigl[{\bar L}-\bigl(\gamma\eta-\chi(\gamma
)\bigr)\bigr]\varphi.
\]
Further, we see that $\psi_{\infty}:=e^w$ is an eigenfunction of
$L_{\infty}+\gamma\eta$
\[
(L_{\infty}+\gamma\eta)\psi_{\infty}=\chi(\gamma)\psi_{\infty}
\]
for the principal eigenvalue $\chi(\gamma)$; cf. \cite{Donvar1}.
Note that ${\bar L}$ is ergodic as is seen in Proposition~\ref{prop4.3}, while
$L_{\infty}$ is not always ergodic.
\end{rem}

\begin{example*} We assume (\ref{eq2.11}) and (\ref{eq2.12}) and that $\beta
(x)=B(x)x+b(x)$, $\alpha(x)=A(x)x+a(x)$ with an $m\times n$ (resp.,
$n\times n$) matrix-valued bounded function $A$ (resp., $B$), and an
$m$ (resp., $n$)-vector-valued bounded function~$a$ (resp., $b$) such that:
\begin{longlist}[(iii)]
\item[(i)]     $A^*A(x)\geq CI_n,   \exists C>0;$
\item[(ii)]     the real parts of all eigenvalues of
$(B^*-A^*A)(x)$   is less than
 $     -C_B$,   $C_B>0;$
\item[(iii)]    $ \operatorname{\mathit{Range}}
( \lambda^*-\sigma^*A ) \subset\operatorname{\mathit{Kernel}}( \sigma)$.
\end{longlist}
In this case
\begin{eqnarray*}
G(x)^*x&\equiv& \beta(x)^*x-{\hat\alpha}^*(\sigma\sigma
^*)^{-1}\lambda^*(x)x \\
&=& \bigl(B(x)x+b(x)\bigr)^*-\bigl(A(x)x+a(x)-r(x){\mathbf1}\bigr)^*(\sigma\sigma
^*)^{-1}\sigma\lambda^*(x)x \\
&=&\bigl(B(x)x+b(x)\bigr)^*x-x^*A^*(x)\bigl((\sigma\sigma^*)^{-1}\sigma\lambda
^*-A\bigr)(x)x \\
&&{} +x^*A^*A(x)x-\bigl(a(x)-r(x){\mathbf1}\bigr)^*(\sigma\sigma
^*)^{-1}\sigma\lambda^*(x)x \\
&=&x^*(B^*-A^*A)(x)x+b(x)^*x-\bigl(a(x)-r(x){\mathbf1}\bigr)^*(\sigma\sigma
^*)^{-1}\sigma\lambda^*(x)x,
\end{eqnarray*}
and we see that (\ref{eq2.20}) holds. Furthermore, (\ref{eq2.22}) holds because of (i).
\end{example*}

\section{H-J-B equations of ergodic type}\label{sec3}
Instead of (\ref{eq2.19}), we shall study an H-J-B equation of ergodic type for
${\bar w}=-w^{(\gamma)}$.
%
\begin{equation}\label{eq3.1}
-\chi=\tfrac{1}{2}\operatorname{tr}[\lambda\lambda^*D^2{\bar w}]+\beta
_{\gamma}^*D{\bar w}-\tfrac{1}{2}(D{\bar w})^*\lambda N_{\gamma
}^{-1}\lambda^*D{\bar w}+U_{\gamma}.
\end{equation}
%
\begin{prop}\label{prop3.1} Assume (\ref{eq2.11}), (\ref{eq2.12}), (\ref{eq2.20}) and (\ref{eq2.21}). Then (\ref{eq3.1}) has a
solution $(\chi,{\bar w})$ such that ${\bar w}\in C^2(R^n)$,
\[
{\bar w}(x) \to\infty   \qquad \mbox{as }   |x| \to\infty,
\]
and the solution satisfying this condition is unique up to additive
constants with respect to ${\bar w}$.
\end{prop}

%
\begin{rem}
The following notation is useful for the task at hand. Set $\Sigma
:=(\sigma\sigma^*)^{-1}\sigma$. Then
\[
\Sigma^*=\sigma^*(\sigma\sigma^*)^{-1}, \qquad    \Sigma\Sigma
^*=(\sigma\sigma^*)^{-1}, \qquad  \Sigma^*(\Sigma\Sigma^*)^{-1}\Sigma
=\sigma^*(\sigma\sigma^*)^{-1}\sigma.
\]
Moreover, we see that
\[
\Sigma N_{\gamma}^{-1}=\frac{1}{1-\gamma}\Sigma,  \qquad   N_{\gamma
}=I-\gamma\Sigma^*(\Sigma\Sigma^*)^{-1}\Sigma=I-\gamma\sigma
^*(\sigma\sigma^*)^{-1}\sigma.
\]
\end{rem}

To prove Proposition~\ref{prop3.1}, we first consider the H-J-B equation of
discounted type,
%
\begin{equation}\label{eq3.2}
\epsilon v_{\epsilon}=\tfrac{1}{2}\operatorname{tr}[\lambda\lambda^*D^2
v_{\epsilon}]+\beta_{\gamma}^*D v_{\epsilon}-\tfrac{1}{2}(D
v_{\epsilon})^*\lambda N_{\gamma}^{-1}\lambda^*D v_{\epsilon
}+U_{\gamma}.
\end{equation}
Note that (\ref{eq3.2}) can be written as
%
\begin{eqnarray}\label{eq3.3}
  \qquad   \epsilon v_{\epsilon}&=&\tfrac{1}{2}\operatorname{tr}[\lambda\lambda^*D^2
v_{\epsilon}]+G^*D v_{\epsilon}-\tfrac{1}{2}(\lambda^* D v_{\epsilon
}-\Sigma^*{\hat\alpha})^* N_{\gamma}^{-1}(\lambda^*D v_{\epsilon
}-\Sigma^*{\hat\alpha})\nonumber\hspace*{-15pt}
\\[-9pt]
\\[-9pt]
&&  {}+\tfrac{1}{2}{\hat\alpha}\Sigma\Sigma
^*{\hat\alpha}.
\nonumber
\end{eqnarray}
Then, we consider the linear equation
%
\begin{equation}\label{eq3.4}
\epsilon\varphi_{\epsilon}=L\varphi_{\epsilon}+\tfrac{1}{2}{\hat
\alpha}\Sigma\Sigma^*{\hat\alpha},
\end{equation}
where
\[
L\varphi=\tfrac{1}{2}\operatorname{tr}[\lambda\lambda^*D^2\varphi
]+G^*D\varphi.
\]
Under assumptions (\ref{eq2.11}), (\ref{eq2.12}) and (\ref{eq2.20}), (\ref{eq3.4}) has a solution
$\varphi_{\epsilon}\in C^2(R^n)$. Indeed, set $\psi_1(x)=c|x|^2$,
$c>0$. Then, by taking $c$ to be sufficiently large, we can see that
there exists $R_0$ such that for $R>R_0$,
\[
L\psi_1+\tfrac{1}{2}{\hat\alpha}\Sigma\Sigma^*{\hat\alpha
}<0 \qquad \mbox{in }  B_R^c.
\]
Therefore, when setting
\[
\Phi_{\epsilon}(x)=\frac{M}{\epsilon}+\psi_1(x), \qquad     M=\sup
_{x\in B_R}\biggl|L\psi_1(x)+\frac{1}{2}{\hat\alpha}\Sigma\Sigma^*{\hat
\alpha}(x)\biggr|,
\]
$\Phi_{\epsilon}(x)$ turns out to be a supersolution to (\ref{eq3.4}), and we
can see that there exists a solution $\varphi_{\epsilon}\in C^2(R^n)$
to (\ref{eq3.4}) such that $0\leq\varphi_{\epsilon}\leq \Phi_{\epsilon
}(x)$ since $v\equiv0$ is a subsolution.

We note that $\varphi_{\epsilon}(x)$ is a supersolution to (\ref{eq3.2})
which is the same equation as (\ref{eq3.3}).\vspace*{-2pt}
%
\begin{lem}\label{lem3.1}
Under the assumptions of Proposition~\ref{prop3.1}, (\ref{eq3.2}) has a solution such
that $v_{\epsilon}\in C^2(R^n)$ and $0\leq v_{\epsilon}\leq\varphi
_{\epsilon}$.\vspace*{-2pt}
\end{lem}

\begin{pf}
In proving the existence of the solution, we introduce a Dirichlet
problem on $B_R$, $R>0$:
%
\begin{eqnarray}\label{eq3.5}
 \qquad  \epsilon v_{\epsilon}&=&\tfrac{1}{2}\operatorname{tr}[\lambda\lambda^*D^2
v_{\epsilon}]+\beta_{\gamma}^*D v_{\epsilon}-\tfrac{1}{2}(D
v_{\epsilon})^*\lambda N_{\gamma}^{-1}\lambda^*D v_{\epsilon
}+U_{\gamma} \quad    \mbox{in }  B_R, \nonumber\hspace*{-15pt}
\\[-9pt]
\\[-9pt]
  \qquad v_{\epsilon}(x)&=&\varphi_{\epsilon},  \qquad    x\in  \partial B_R.
\nonumber
\end{eqnarray}
Owing to Theorem 8.3 (\cite{Lad}, Chapter 4), Dirichlet problem (\ref{eq3.4})
has a solution $v_{\epsilon}$. We extend $v_{\epsilon}$ to the whole
Euclidean space as
\[
v_{\epsilon, R}=\cases{\displaystyle
v_{\epsilon}(x),& \quad $x\in B_R$, \cr\displaystyle
\varphi_{\epsilon},& \quad $x\in B_R^c$.
}
\]
Then we can see that $v_{\epsilon,R}$ is nonincreasing with respect to
$R$. Indeed, for $R<R'$, $v_{\epsilon,R}$ is a supersolution to (\ref{eq3.3})
in $B_{R'}$, and
we have
\begin{eqnarray*}
 &&\epsilon(v_{\epsilon,R}-v_{\epsilon, R'})\\
 && \qquad \geq\tfrac{1}{2}\operatorname
{tr}[\lambda\lambda^*D^2(v_{\epsilon,R}-v_{\epsilon,R'})]+\beta
_{\gamma}^*D(v_{\epsilon,R}-v_{\epsilon,R'}) \\
&& \qquad  \quad
{}-\tfrac{1}{2}(D v_{\epsilon,R})^*\lambda N_{\gamma}^{-1}\lambda^*D
v_{\epsilon,R}+
\tfrac{1}{2}(D v_{\epsilon,R'})^*\lambda N_{\gamma}^{-1}\lambda^*D
v_{\epsilon,R'}    \qquad \mbox{in }  B_R' \\
&& \qquad =
\tfrac{1}{2}\operatorname{tr}[\lambda\lambda^*D^2( v_{\epsilon
,R}-v_{\epsilon,R'})]+\beta_{\gamma}^*D( v_{\epsilon,R}-v_{\epsilon
,R'}) \\
&& \qquad  \quad {}
-\tfrac{1}{2}(D v_{\epsilon,R}+D v_{\epsilon,R'})^*\lambda N_{\gamma
}^{-1}\lambda^*D( v_{\epsilon,R}-v_{\epsilon,R'}).
\end{eqnarray*}
Therefore, from the maximum principle (cf. Theorem 3.1 in \cite{GT})
we see that
%
\begin{equation}\label{eq3.6}
v_{\epsilon,R}-v_{\epsilon,R'}\geq0
\end{equation}
since $ v_{\epsilon,R}(x)=v_{\epsilon,R'}(x),   x\in\partial
B_{R'}$. We further note that
%
\begin{equation}\label{eq3.7}
v_{\epsilon,R}\geq0
\end{equation}
for each $R$ because $\psi_0(x)\equiv0$ is a subsolution to (\ref{eq3.2}),
and the maximum principle again applies.

Similar to the proof of Lemma~2.6 in \cite{KaiN}, we have the
following gradient estimate: for each $R$ and $r<\frac{R}{2}$,
%
\begin{equation}\label{eq3.8}
\|\nabla v_{\epsilon,R}\|_{L^{\infty}(B_r)}\leq M_r,
\end{equation}
where $M_r$ is a constant independent of $R$, $\epsilon$. Thus, when
taking a sequence $R_n$ such that $R_n\uparrow\infty$, $v_{\epsilon,
R_n}$ forms a family of uniformly bounded and equicontinuous functions.
Thus we can choose a subsequence $v_{\epsilon,R_{n_k}}$ converging to
a continuous function $v_{\epsilon}$. Furthermore, since
%
\begin{equation}\label{eq3.9}
\| v_{\epsilon, R_n}\|_{H^1(B_r)}\leq M_r'
\end{equation}
for a positive constant $M_r'$ independent of $R_n$ and $\epsilon$, it
converges weakly in $H^1_{\mathrm{loc}}(R^n)$ to $v_{\epsilon}$ by taking a
subsequence if necessary. By similar arguments to Lemma~{6.8} in \cite
{KaiN}, the convergence can be strengthened as $\nabla v_{\epsilon
,R_{n_k}}$ converges strongly in $L^2_{\mathrm{loc}}(R^n)$ to $\nabla
v_{\epsilon}$. As a result we can see from the regularity theorems
that we have a solution $v_{\epsilon}\in C^2(R^n)$ to (\ref{eq3.2}). Since
$v_{\epsilon,R}\leq\varphi_{\epsilon}$, for each $R>0$ from the
maximum principle as well as (\ref{eq3.7}), we see that $0\leq v_{\epsilon
}\leq\varphi_{\epsilon}$.
\end{pf}

Set
\[
\psi_{\delta}(x):=e^{\delta|x|^2}, \qquad\delta>0.
\]
Then, by taking $\delta$ to be sufficiently small, we can see that
there exists $R_1$ such that for $R>R_1$,
\[
L\psi_{\delta}(x) <-1\qquad\mbox{in }   B_R^c.
\]
Therefore, we see that $L$ and $\psi_{\delta}$ satisfy assumption
(\ref{eq9.3}) in the last section.
Set $K(x;\psi_{\delta})=-L\psi_{\delta}$,
\[
F_{\psi}:=\biggl\{u(x)\in W^{2,p}_{\mathrm{loc}}(R^n); \esssup\limits_{x\in B_R^c}\frac
{|u(x)|}{\psi_{\delta}(x)}<\infty\biggr\}
\]
and
\[
F_K:=\biggl\{f(x)\in L^{\infty}_{\mathrm{loc}}(R^n); \esssup\limits_{x\in B_R^c}\frac
{|f(x)|}{K(x;\psi_{\delta})}<\infty\biggr\}.
\]
Then for $f\in F_K$ there exists a solution $\varphi\in F_{\psi}$ to
\[
0=L\varphi+f
\]
if and only if
\[
\int f(x)m(dx)=0,
\]
where $m(dx)$ is an invariant measure for $L$; cf. Proposition~\ref{prop9.4} in
in the \hyperref[appm]{Appendix}.
Therefore, setting
%
\begin{equation}\label{eq3.10}
\chi_0=\int\frac{1}{2}{\hat\alpha}\Sigma\Sigma^*{\hat\alpha
}(x)m(dx),
\end{equation}
there exists a solution $\varphi_0\in F_{\psi}$ to
\[
\chi_0=L\varphi_0+ \tfrac{1}{2}{\hat\alpha}\Sigma\Sigma^*{\hat
\alpha}(x),
\]
and it is known that $\epsilon\varphi_{\epsilon}$ converges to $\chi
_0$ as $\epsilon\to0$ uniformly on each compact set.

Now we can prove Proposition~\ref{prop3.1}.
\begin{pf*}{Proof of Proposition~\ref{prop3.1}}
We first note that
\[
0\leq v_{\epsilon}\leq\varphi_{\epsilon}
\]
because of Lemma~\ref{lem3.1}. Therefore, we have
\[
\|\epsilon v_{\epsilon}\|_{L^{\infty}(B_r)}\leq K_r,
\]
where $K_r$ is a constant independent of $\epsilon$. Moreover,
\[
\|\nabla v_{\epsilon}\|_{L^{\infty}(B_r)}\leq K_r'
\]
for a positive constant $K_r'$ independent of $\epsilon$ in view of (\ref{eq3.8}).
Thus, similarly to the proof of Theorem 3.1 in \cite{KaiN}, we can
prove the existence of the solution $(-\chi,{\bar w})$ to (\ref{eq3.1}) such
that ${\bar w}\in W^{2,p}_{\mathrm{loc}}$. From regularity theorems we see that
${\bar w}\in C^2(R^n)$. The proof of uniqueness is similar to the proof
of Lemma 3.2 in \cite{N0}.
\end{pf*}

Now we have the following proposition.
%
\begin{prop} Under the assumptions of Proposition~\ref{prop3.1}, the solution
${\bar w}$ to (\ref{eq3.1}) satisfies
%
\begin{equation}\label{eq3.11}
|\nabla{\bar w}(x)|^2\leq c(|x|^2+1),
\end{equation}
where $c$ is a positive constant.
If we further assume (\ref{eq2.22}), then, for each $\gamma_0<0$, there exists
a positive constant $c(\gamma_0)$ such that the nonnegative solution
${\bar w}(x)={\bar w}(x;\gamma)$, $\gamma\leq\gamma_0$, satisfies
%
\begin{equation}\label{eq3.12}
{\bar w}(x)\geq c(\gamma_0)|x|^2, \qquad|x|\geq\exists R'.
\end{equation}
\end{prop}

\begin{pf} Set $Q_{\gamma}:=\lambda N_{\gamma}^{-1}\lambda^*$.
Then we shall prove for each $x_0\in R^n$ that
%
\begin{eqnarray}\label{eq3.13}
 |\nabla{\bar w}|^2(x_0)&\leq& K\biggl(|\nabla Q_{\gamma}|^2_r+\frac
{1}{r^2}|Q_{\gamma}|^2_r+|\beta_{\gamma}|^2_r\nonumber
\\[-8pt]
\\[-8pt]  &&\hphantom{K\biggl(}{}+|U_{\gamma
}|_r+|\nabla U_{\gamma}|_r+\frac{|\beta_{\gamma}|_r}{r}+|\nabla
\beta_{\gamma}|_r+c\biggr)
\nonumber
\end{eqnarray}
for positive constants $K$ and $c$, where $|f|_r=|f|_{L^{\infty
}(B_r(x_0))}$. Note that (\ref{eq3.13}) implies (\ref{eq3.11}) because of our
assumptions on the coefficients $\sigma, \lambda, \beta, \alpha$
and~$r$.

We have $\chi(\gamma)\leq0$ since $\epsilon v_{\epsilon}\geq0$ and
$ \epsilon v_{\epsilon}\to-\chi(\gamma)\leq\chi_0$ as $\epsilon
\to0$. In the following $\beta_{\gamma}$, $Q_{\gamma}$ and
$U_{\gamma}$ are abbreviated to $\beta$, $Q$ and $U$, respectively.
$|\cdot|_r$ is abbreviated to $|\cdot|$.

By differentiating (\ref{eq3.1}) with respect to $x_k$, we have
%
\begin{eqnarray}\label{eq3.14}
 0&=&\tfrac{1}{2}(\lambda\lambda^*)^{ij}D_{ijk}w+\tfrac{1}{2}(\lambda
\lambda^*)^{ij}_kD_{ij}w+\beta^iD_{ik}w+\beta^i_kD_iw\nonumber
\\[-8pt]
\\[-8pt]
&&{} -D_iwQ^{ij}D_{jk}w-\tfrac{1}{2}D_iwQ^{ij}_kD_jw+U_k.
\nonumber
\end{eqnarray}
Set
\[
F=|\nabla w|^2=\sum_{k=1}^n|D_kw|^2.
\]
Then we have
\begin{eqnarray*}
&&-\frac{1}{2}(\lambda\lambda^*)^{ij}D_{ij}F-\beta
^iD_iF+Q^{ij}D_iwD_jF \\
&& \qquad =-(\lambda\lambda^*)^{ij}D_{jk}wD_{ik}w\\
&& \qquad \quad  {}-D_kw\{(\lambda\lambda
^*)^{ij}D_{ijk}w+2\beta^iD_{ik}w-2Q^{ij}D_jwD_{ik}w\} \\
&& \qquad
=-(\lambda\lambda^*)^{ij}D_{jk}wD_{ik}w\\
&& \qquad \quad  {}+D_kw\{(\lambda\lambda
^*)^{ij}_kD_{ij}w+2\beta_k^iD_iw-D_iwQ_k^{ij}D_jw+2U_k\} \\
&& \qquad \leq-\frac{1}{2nc_2}\{(\lambda\lambda^*)^{ij}D_{ij}w\}^2-\frac
{1}{2}(\lambda\lambda^*)^{ij}D_{jk}wD_{ik}w+\frac{c}{2\delta
}|\nabla w|^2+\frac{c\delta}{2}|D^2w|^2 \\
&& \qquad \quad  {}+2|\nabla\beta||\nabla w|^2+|\nabla Q|^2|\nabla w|^3+2|\nabla
U||\nabla w| \\
&& \qquad \leq-\frac{1}{2nc_2}(-2\chi-2\beta
^iD_iw+D_iwQ^{ij}D_jw-2U)^2+\frac{2c_2}{\delta}|\nabla w|^2 \\
&& \qquad \quad {}+2|\nabla\beta||\nabla w|^2+|\nabla Q|^2|\nabla w|^3+2|\nabla
U||\nabla w|.
\end{eqnarray*}
Here we have used (\ref{eq3.14}) and the matrix inequality
\[
(\operatorname{tr}[AB])^2\leq n C\operatorname{tr}[AB^2]
\]\vfill\eject
for symmetric matrix $B$ and nonnegative definite symmetric matrix $A$,
where $C$ is the maximum eigenvalue of $A$. Set
\[
\tau(x):=\cases{\displaystyle
\biggl(\frac{|x-x_0|^2}{r^2}-1\biggr)^2,& \quad  |$x-x_0|\leq r$,\vspace*{2pt}\cr\displaystyle
0,& \quad $|x-x_0|>r$.
}
\]
Then $\operatorname{tr}[\lambda\lambda^*D^2\tau]\geq-\frac{4n}{r^2}c_2,
$ $(D\tau)^*\lambda\lambda^*D\tau\leq\frac{16c_2}{r^2}\tau $
and $|D\tau|^2\leq\frac{16 c_2}{c_1 r^2}\tau$. Let $x$ be the
maximum point of $\tau F$ in $B_r(x_0)$. Then $D(\tau F)(x)=0$ and
$\operatorname{tr}[\lambda\lambda^*\times D^2(\tau F)](x)\leq0$. Therefore, from
the maximum principle we have
\begin{eqnarray*}
0&\leq&-\frac{1}{2}(\lambda\lambda^*)^{ij}D_{ij}(\tau F)-\beta
^iD_i(\tau F)+Q^{ij}D_jwD_i(\tau F) \\
&=&\tau\biggl\{-\frac{1}{2}(\lambda\lambda^*)^{ij}D_{ij}F-\beta
^iD_iF+Q^{ij}D_jwD_iF\biggr\} \\
&&{}-\frac{1}{2}(\lambda\lambda^*)^{ij}D_{ij}\tau F-(\lambda
\lambda^*)^{ij}D_i\tau D_jF-(\beta^iD_i\tau)F+(Q^{ij}D_jwD_i\tau
)F \\
&\leq&\tau\biggl[-\frac{1}{2nc_2}(-2\chi-2\beta
^iD_iw+D_iwQ^{ij}D_jw-2U)^2+\frac{2c_2}{\delta}|\nabla w|^2 \\
&&\hspace*{83pt}{}+2|\nabla\beta||\nabla w|^2+|\nabla Q|^2|\nabla w|^3+2|\nabla
U||\nabla w|\biggr] \\
&&{}-F\biggl\{\frac{1}{2}(\lambda\lambda^*)^{ij}D_{ij}\tau-\frac{(\lambda
\lambda^*)^{ij}D_i\tau D_j\tau}{\tau} -\beta^iD_i\tau
+Q^{ij}D_jwD_i\tau\biggr\}.
\end{eqnarray*}
Since $\frac{1}{1-\gamma}\lambda\lambda^*\leq Q\leq\lambda\lambda
^*$, by taking $\delta$ to be sufficiently small,
\[
c(\gamma)|Dw|^2\leq-2\beta^*Dw+(Dw)^*QDw+\frac{1}{\delta}|\beta
|^2\leq(c_2+1)|Dw|^2+\biggl(1+\frac{1}{\delta}\biggr)|\beta|^2
\]
for a positive constant
%
\begin{equation}\label{eq3.15}
c(\gamma)=\frac{c_1}{1-\gamma}-\delta>0.
\end{equation}
Therefore, it follows that
\begin{eqnarray*}
0&\leq&-\tau\biggl(-2\beta^*Dw+(Dw)^*QDw+\frac{1}{\delta}|\beta
|^2\biggr)^2 \\
&&{}+2\tau\biggl(-2\beta^*Dw+(Dw)^*QDw+\frac{1}{\delta}|\beta|^2\biggr)\biggl(\frac
{1}{\delta}|\beta|^2+2U+2\chi\biggr) \\
&&{}-\tau\biggl(\frac{1}{\delta}|\beta|^2+2U+2\chi\biggr)^2+\tau(2|\nabla
\beta||\nabla w|^2+|\nabla Q||\nabla w|^3+2|\nabla U||\nabla
w|) \\
&&{}+\frac{2nc_2}{r^2}F+\frac{16c_2}{r^2}F+|\beta|\frac{4\sqrt
{c_2}}{\sqrt{c_1}r}\tau^{1/2}F+\frac{4\sqrt{c_2}}{\sqrt{c_1}r}\tau
|Q|F^{3/2}
 \\
&\leq& -\tau c(\gamma)^2|\nabla w|^4+2\tau\biggl\{(c_2+1)\nabla w|^2+\biggl(1+\frac
{1}{\delta}\biggr)|\beta|^2\biggr\}\biggl(\frac{1}{\delta}|\beta|^2+2U\biggr) \\
&&{}+\biggl(\frac{2nc_2}{r^2}+\frac{16 c_2}{r^2}\biggr)F+|\beta|\frac{4\sqrt
{c_2}}{\sqrt{c_1}r}\tau^{1/2}F+\frac{4\sqrt{c_2}}{\sqrt{c_1}r}\tau
|Q|F^{3/2} \\
&&{}+\tau(2|\nabla\beta|F+|\nabla Q|F^{3/2}+2|\nabla U|F^{1/2}).
\end{eqnarray*}
We can assume $F\geq|\beta|^2$ and $F\geq|\nabla U|$; thus,
\begin{eqnarray*}
 0&\leq&-c(\gamma)^2 \tau F^2+2\biggl(c_2+2+\frac{1}{\delta}\biggr)\tau F\biggl( \frac
{1}{\delta}|\beta|^2+2U\biggr) \\
&&{}
+\biggl(\frac{2nc_2}{r^2}+\frac{16 c_2}{r^2}\biggr)F+|\beta|\frac{4\sqrt
{c_2}}{\sqrt{c_1}r}
\tau^{1/2}F+\frac{4\sqrt{c_2}}{\sqrt{c_1}r}\tau|Q|F^{3/2} \\
&&{}+\tau(2|\nabla\beta|F+|\nabla Q|F^{3/2}+2F^{3/2}).
\end{eqnarray*}
Accordingly, we have
\begin{eqnarray*}
 0&\leq&-c(\gamma)^2\tau F+\biggl(|\nabla Q|+\frac{4\sqrt{c_2}}{\sqrt
{c_1}r}|Q|+2\biggr)(\tau F)^{1/2} \\
&&{}+2\biggl(c_2+2+\frac{1}{\delta}\biggr)\biggl(\frac{1}{\delta}|\beta|^2+2U\biggr)+\frac
{2nc_2}{r^2}+\frac{16c_2}{r^2}+\frac{4|\beta|\sqrt{c_2}}{\sqrt
{c_1}r}+2|\nabla\beta|.
\end{eqnarray*}
Therefore, we obtain
\begin{eqnarray*}
\frac12 c(\gamma)^2\tau F&\leq&\frac{1}{2c(\gamma)^2}\biggl(|\nabla Q|+\frac{4\sqrt
{c_2}}{\sqrt{c_1}r}|Q|+2\biggr)^2\\
&&{}+c_{\delta}\biggl(\frac{1}{\delta}|\beta
|^2+2U\biggr)+\frac{c}{r^2}+\frac{c|\beta|}{r}+2|\nabla\beta|,
\end{eqnarray*}
with $c_{\delta}=2(c_2+2+\frac{1}{\delta})$ and universal constant $c>0$.
Including the case where $|\beta|^2\geq F$, $|\nabla U|\geq F$, we obtain
\begin{eqnarray*}
F(x_0)&=&\tau(x_0)F(x_0)\leq(\tau F)(x) \\
&\leq&\frac{c}{c(\gamma)^4}\biggl(|\nabla Q|^2+\frac
{1}{r^2}|Q|^2+c\biggr) \\
&&{}+\frac{c'_{\delta}}{c(\gamma)^2}\biggl(|\beta|^2+U+\frac
{|\beta|}{r}+|\nabla\beta|+\frac{1}{r}\biggr)+|\nabla U|,
\end{eqnarray*}
and (\ref{eq3.13}) has been proved.

Now let us prove (\ref{eq3.12}). For each $\rho>0$ take a point $x_{\rho}\in
R^n$ such that $|x_{\rho}|=\rho$. Set
\[
R(x)=c_{\rho}\biggl(1-\frac{4|x-x_{\rho}|^2}{\rho^2}\biggr) \qquad \mbox{in }
D_{\rho}=\biggl\{x; |x-x_{\rho}|\leq\frac{\rho}{2}\biggr\},
\]
where $c_{\rho}$ is a positive constant determined later. Then,
$R(x)\geq0$ in $D_{\rho}$ and $R(x)=0$ on $\partial D_{\rho}$. Set
\[
z(x)={\bar w}(x)-R(x).
\]
Then,
\[
z(x)={\bar w}(x)\geq0, \qquad   x\in\partial D_{\rho}.
\]
Note that
\[
\xi^*\lambda N_{\gamma}^{-1}\lambda^*\xi\leq c_2|\xi|^2,\qquad\xi
\in R^n.
\]
Then we have
\begin{eqnarray*}
&&-\chi-\frac{1}{2}\operatorname{tr}[\lambda\lambda^*D^2z]-\beta
^*Dz \\
&& \qquad =-\frac{1}{2}(D{\bar w})^*\lambda N_{\gamma}^{-1}\lambda^*D{\bar
w}+U_{\gamma}+\frac{1}{2}\operatorname{tr}[\lambda\lambda^*D^2R]+\beta
^*DR \\
&& \qquad =-\frac{1}{2}D({\bar w}+R)^*\lambda N_{\gamma}^{-1}\lambda
^*D({\bar w}-R)-\frac{1}{2}(DR)^*\lambda N_{\gamma}^{-1}\lambda
^*DR+U_{\gamma} \\
&& \qquad  \quad {} +\frac{1}{2}\operatorname{tr}[\lambda\lambda^*D^2R]+\beta^*DR \\
&& \qquad \geq
-\frac{1}{2}D({\bar w}+R)^*\lambda N_{\gamma}^{-1}\lambda^*Dz-\frac
{c_2}{2}|DR|^2+U_{\gamma}+\frac{1}{2}\operatorname{tr}[\lambda\lambda
^*D^2R]+\beta^*DR.
\end{eqnarray*}
Noting that $|\beta_{\gamma}(x)|\leq c\rho,  x\in D_{\rho}$, for a
positive constant independent of $\gamma$,
\begin{eqnarray*}
&&-\chi-\frac{1}{2}\operatorname{tr}[\lambda\lambda^*D^2z]-\beta
^*Dz+\frac{1}{2}D({\bar w}+R)^*\lambda N_{\gamma}^{-1}\lambda^*Dz
 \\
&& \qquad \geq-\frac{c_2}{2}|DR|^2+U(x)-\frac{4c_{\rho}}{\rho^2}\operatorname
{tr}[\lambda\lambda^*]-4c c{\rho} \\
&& \qquad \geq
-\frac{8c_2c_{\rho}^2}{\rho^2}+\frac{-\gamma}{2(1-\gamma
)}c_0\biggl(\frac{|\rho|^2}{4}+1\biggr)-\frac{4c_2nc_{\rho}}{\rho^2}-4cc_{\rho
} \\
&& \qquad \geq-\biggl(8c_2\frac{c_{\rho}^2}{\rho^2}+4c_2n\frac{c_{\rho}}{\rho
^2}+4cc_{\rho}\biggr)+\frac{-\gamma_0c_0\rho^2}{8(1-\gamma_0)}+\frac
{-\gamma_0c_0}{2(1-\gamma_0)}.
\end{eqnarray*}
By setting
$
c_{\rho}=c(\gamma_0)\rho^2
$ with $c(\gamma_0) $ such that\vspace*{-1pt} $8c_2c(\gamma_0)^2+4cc(\gamma
_0)<\frac{-\gamma_0c_0}{8(1-\gamma_0)}$ and $4c_2nc(\gamma_0)<\frac
{-\gamma_0c_0}{2(1-\gamma_0)}$,
we see that
\[
-\tfrac{1}{2}\operatorname{tr}[\lambda\lambda^*D^2z]-\beta^*Dz+\tfrac
{1}{2}D({\bar w}+R)^*\lambda N_{\gamma}^{-1}\lambda^*Dz\geq M>0\qquad
\mbox{in }  D_{\rho}
\]
for some positive constant and sufficiently large $\rho$. Then $z$ is
superharmonic in $D_{\rho}$ and $z(x)\geq0, x \in\partial D_{\rho
}$. Therefore $z(x)\geq0,    x\in D_{\rho}$, from which we have
\[
z(x_{\rho})={\bar w}(x_{\rho})-c_{\rho}\geq0.
\]
Hence, ${\bar w}(x_{\rho})\geq c(\gamma_0)\rho^2$.
\end{pf}

\section{H-J-B equations and related stochastic control problems}\label{sec4}
Let us come back to H-J-B equation (\ref{eq2.16}). According to assumption
(\ref{eq2.11}), we have a~positive constant $c_{\beta}$ such that
\[
|\beta_{\gamma}(x)|^2\leq c_{\beta}(|x|^2+1).
\]
We strengthen condition (\ref{eq2.21}) to (\ref{eq2.22}). Then we have the following lemma.

\begin{lem}\label{lem4.1}Assume (\ref{eq2.11}), (\ref{eq2.12}), (\ref{eq2.22}) and $v_0\geq1$. Then for each
$t<T$ there exist positive constants $k=k(T-t)$ and $k'=k'(T-t)$ such that
%
\begin{equation}\label{eq4.1}
{\bar v}(t,x;T)\geq k|x|^2-k'.
\end{equation}
\end{lem}

\begin{pf} Choose a positive constant $c$ such that
\[
c_{\gamma}-\frac{c}{2}c_{\beta}>0,
\]
and set $b=c_{\gamma}-\frac{c}{2}c_{\beta}$, where $c_{\gamma
}=-\frac{\gamma c_0}{2(1-\gamma)}$, and set
\[
R(t,x):=\tfrac{1}{2}x^*P(t)x+q(t),
\]
where $P(t)$ is a solution to the Riccati equation
%
\begin{equation}\label{eq4.2}
{\dot
P}(t)-\biggl(\frac{c_2}{1-\gamma}+\frac{1}{c}\biggr)P(t)I_nP(t)+bI_n=0,\qquad
P(T)=0,
\end{equation}
and $q(t)$ is a solution to the ordinary equation
%
\begin{equation}\label{eq4.3}
{\dot q}(t)+\frac{c_1}{2}\operatorname{tr}[P(t)]-\frac{c c_{\beta
}}{2}-c_{\gamma}'=0,\qquad q(T)=-\gamma\log v_0,
\end{equation}
where $c_{\gamma}'=-\frac{c_0'\gamma}{2(1-\gamma)}$. Set
\[
z(t,x):={\bar v}(t,x)-R(t,x).
\]
Then
\begin{eqnarray*}
&&
-\frac{\partial z}{\partial t}-\frac{1}{2}\operatorname{tr}[\lambda\lambda
^*D^2z]-\beta_{\gamma}^*Dz \\
&& \qquad = \frac{1}{2}(D{\bar v})^*\lambda N_{\gamma}^{-1}\lambda^*D{\bar
v}+U_{\gamma}+\frac{\partial R}{\partial t}+\frac{1}{2}\operatorname
{tr}[\lambda\lambda^*D^2R]+\beta^*DR \\
&& \qquad =-\frac{1}{2}D({\bar v}+R)^*\lambda N_{\gamma}^{-1}\lambda^*
D({\bar v}-R)-\frac{1}{2}DR^*\lambda N_{\gamma}^{-1}\lambda^*
DR+U_{\gamma} \\
&& \qquad  \quad {}+\frac{\partial R}{\partial t}+\frac{1}{2}\operatorname{tr}[\lambda
\lambda^*D^2R]+\beta^*DR \\
&& \qquad
\geq-\frac{1}{2}D({\bar v}+R)^*\lambda N_{\gamma}^{-1}\lambda^*
D({\bar v}-R)-\frac{c_2}{2(1-\gamma)}(DR)^*I_nDR+c_{\gamma
}|x|^2-c_{\gamma}' \\
&& \qquad  \quad {}+\frac{1}{2}x^*{\dot P}(t)x+{\dot q}(t)+\frac{c_1}{2}\operatorname
{tr}[P(t)]-\frac{c}{2}\beta_{\gamma}^*\beta_{\gamma}-\frac{1}{2c}(DR)^*DR.
\end{eqnarray*}
Therefore,
\begin{eqnarray*}
&&-\frac{\partial z}{\partial t}-\frac{1}{2}\operatorname{tr}[\lambda
\lambda^*D^2 z]-\beta^*Dz+\frac{1}{2}D({\bar v}+R)^*\lambda
N_{\gamma}^{-1}\lambda^*Dz \\
&& \qquad \geq\frac{1}{2}x^*{\dot P}(t)x-\frac{1}{2}\biggl(\frac{c_2}{1-\gamma
}+\frac{1}{c}\biggr)x^*P(t)I_nP(t)x+\biggl(c_{\gamma}-\frac{cc_\beta
}{2}\biggr)|x|^2 \\
&& \qquad  \quad {} +{\dot q}(t)+\frac{c_1}{2}\operatorname{tr}[P(t)]-\frac
{cc_{\beta}}{2}-c_{\gamma}' \\
&& \qquad \geq\frac{1}{2}\biggl( c_{\gamma}-\frac{cc_\beta}{2}\biggr)|x|^2\geq0.
\end{eqnarray*}
Thus we see that $z(t,x)$ is super harmonic in $[0,T)\times R^n$, and
$z(T,x)=0$. Therefore we have $z(t,x)={\bar v}(t,x)-R(t,x)\geq0$, that is,
\[
{\bar v}(t,x)\geq R(t,x)=\tfrac{1}{2}x^*P(t)x+q(t).
\]
Since $P(t)$ is positive definite,
\[
{\bar v}(t,x)\geq k|x|^2-k',\qquad k=k(T-t),  k'=k'(T-t)>0.
\]
\upqed
\end{pf}

Let us rewrite (\ref{eq2.16}) as
%
\begin{equation}\label{eq4.4}
\cases{\displaystyle
0=\frac{\partial{\bar v}}{\partial t}+\frac{1}{2}\operatorname
{tr}[\lambda\lambda^*D^2{\bar v}]+G^*D{\bar v} \cr\displaystyle
 \hphantom{0=}{}-\frac{1}{2}(\lambda
^*D{\bar v}-\Sigma^*{\hat\alpha})^* N_{\gamma}^{-1}(\lambda
^*D{\bar v}-\Sigma^*{\hat\alpha})+\frac{1}{2}{\hat\alpha}^*\Sigma\Sigma^*{\hat
\alpha}, \vspace*{4pt}\cr\displaystyle
 {\bar v}(T,x)=-\gamma\log v_0.
}
\end{equation}
Noting that
\begin{eqnarray*}
&&-\frac{1}{2}(\lambda^*D{\bar v}-\Sigma^*{\hat\alpha})^*N_{\gamma
}^{-1}(\lambda^*D{\bar v}-\Sigma^*{\hat\alpha}) \\
&& \qquad =\inf_{z\in R^{n+m}}\biggl\{\frac{1}{2}z^*N_{\gamma}z-z^*\Sigma^*{\hat
\alpha}+(\lambda z)^*D{\bar v}\biggr\} \\
&& \qquad =\inf_{z\in R^{n+m}}\biggl[\frac{1}{2}\{z+N_{\gamma}^{-1}(\lambda
^*D{\bar v}-\Sigma^*{\hat\alpha})\}^*N_{\gamma}\{z+N_{\gamma
}^{-1}(\lambda^*D{\bar v}-\Sigma^*{\hat\alpha})\} \\
&&\hspace*{109pt} \qquad  \quad {} -\frac{1}{2}(\lambda^*D{\bar v}-\Sigma^*{\hat\alpha
})^*N_{\gamma}^{-1}(\lambda^*D{\bar v}-\Sigma^*{\hat\alpha})\biggr],
\end{eqnarray*}
we can rewrite it again as
%
\begin{equation}\label{eq4.5}\qquad
\cases{\displaystyle
0=\frac{\partial{\bar v}}{\partial t}+\frac{1}{2}\operatorname
{tr}[\lambda\lambda^*D^2{\bar v}]+G^*D{\bar v}+\inf_{z\in R^{n+m}}\{
(\lambda z)^*D{\bar v}+\varphi(x,z)\},\cr\displaystyle
 {\bar v}(T,x)=-\gamma\log v_0,
}
\end{equation}
where
\[
\varphi(x,z)=\tfrac{1}{2}z^* N_{\gamma} z-z^*\Sigma^*{\hat\alpha
}+\tfrac{1}{2}{\hat\alpha}^*\Sigma\Sigma^*{\hat\alpha},\qquad
N_{\gamma}=I-\gamma\Sigma^*(\Sigma\Sigma^*)^{-1}\Sigma.
\]
This H-J-B equation corresponds to the following stochastic control
problem, the value of which is defined as
%
\begin{equation}\label{eq4.6}
\inf_{Z_{\cdot}\in{\tilde{\mathcal A}}(T)}E\biggl[\int_0^T\varphi
(Y_s,Z_s)\,ds-\gamma\log v_0\biggr],
\end{equation}
where $Y_t$ is a controlled process governed by the stochastic
differential equation
%
\begin{equation}\label{eq4.7}
dY_t=\lambda(Y_t)\,dW_t+\{G(Y_t)+\lambda(Y_t)Z_t\}\,dt, \qquad    Y_0=x,
\end{equation}
with control $Z_t\in{\tilde{\mathcal A}}(T)$. Here, ${\tilde
{\mathcal A}}(T)$ is the set of all $R^{n+m}$ valued progressively
measurable processes such that
\[
E\biggl[\int_0^T|Z_s|^2\,ds\biggr]<\infty.
\]
To study this problem, we introduce a value function for $0\leq t\leq T$,
\[
v_*(t,x)=\inf_{Z_{\cdot}\in{\tilde{\mathcal A}}(T-t)}E\biggl[\int
_0^{T-t}\varphi(Y_s,Z_s)\,ds -\gamma\log v_0\biggr].
\]
By the verification theorem, the solution ${\bar v}$ to (\ref{eq4.5}) can be
identified with the value function $v_*$. Indeed, set
\[
{\hat z}(s,x)=-N_{\gamma}^{-1}(\lambda^*D{\bar v}-\Sigma^*{\hat
\alpha})(s,x),
\]
which attains the infimum in (\ref{eq4.5}), and consider the stochastic
differential equation
%
\begin{equation}\label{eq4.8}
d{\hat Y}_t=\lambda({\hat Y}_t)\,dW_t+\{G({\hat Y}_t)+\lambda({\hat
Y}_t){\hat Z}(t,{\hat Y_t})\}\,dt, \qquad    Y_0=x.
\end{equation}
From the estimates obtained in Theorem~\ref{theo2.1}, we see that (\ref{eq4.8}) has a
unique solution. It is also seen by using It\^o's theorem that
\[
{\bar v}(0,x)=E\biggl[\int_0^T\varphi({\hat Y}_s,{\hat Z}_s)\,ds-\gamma\log v_0\biggr]
\]
holds, where ${\hat Z}_s={\hat Z}(s,{\hat Y}_s)$. In a similar way, we
can see that
\[
{\bar v}(0,x)\leq E\biggl[\int_0^T\varphi(Y_s,Z_s)\,ds-\gamma\log v_0\biggr]
\]
for each $Z_{\cdot}\in{\tilde{\mathcal A}}(T)$, hence, ${\bar
v}(0,x)=v_*(0,x)$.\vadjust{\goodbreak}

Let us consider the following stochastic control problem with the
averaging cost criterion:
%
\begin{equation}\label{eq4.9}
\rho(\gamma)=\inf_{Z_{\cdot}\in{\tilde{\mathcal A}}}\limsup_{T\to
\infty}\frac{1}{T}E\biggl[\int_0^T\varphi(Y_s,Z_s)\,ds\biggr],\
\end{equation}
where $Y_t$ is a controlled process governed by controlled stochastic
differential equation (\ref{eq4.7}) with control $Z_t$. The solution $Y_t$ of
(\ref{eq4.7}) is sometimes written as $Y_t^{(Z)}$ to make clear the dependence
on control $Z_t$, and the set ${\tilde{\mathcal A}}$ of all admissible
controls is defined as follows:
\begin{eqnarray*}
 {\tilde{\mathcal A}}&=&\biggl\{Z_{\cdot}; Z_t  \mbox{ is an $R^{n+m}$ valued
progressively measurable process such that}  \\
&&\hspace*{77.5pt} \limsup_{T\to\infty}\frac{1}{T}E\bigl[\bigl|Y_T^{(Z)}\bigr|^2\bigr]=0,
 E\biggl[\int_0^T|Z_s|^2\,ds\biggr]<\infty, \forall T\biggr\}.
\end{eqnarray*}
Corresponding to this stochastic control problem, H-J-B equation of
ergodic type (\ref{eq3.1}) can be written as
%
\begin{equation}\label{eq4.10}
 \qquad -\chi(\gamma)=\frac{1}{2}\operatorname{tr}[\lambda\lambda^*D^2{\bar
w}]+G^*D{\bar w}+\inf_{z\in R^{n+m}}\{(\lambda z)^*D{\bar w}+\varphi
(x,z)\}.
\end{equation}
We then set
%
\begin{equation}\label{eq4.11}
{\hat z}(x)=-N_{\gamma}^{-1}(\lambda^* D{\bar w}-\Sigma^*{\hat
\alpha})(x),
\end{equation}
and consider stochastic differential equation
%
\begin{eqnarray}\label{eq4.12}
d{\bar Y}_t&=&\lambda({\bar Y}_t)\,dW_t+\{G({\bar Y}_t)+\lambda({\bar
Y}_t){\hat z}({\bar Y}_t)\}\,dt\nonumber \\
&=&\lambda({\bar Y}_t)\,dW_t+\{\beta_{\gamma}-\lambda N_{\gamma
}^{-1}\lambda^*D{\bar w}\}({\bar Y}_t)\,dt, \\
{\bar Y}_0&=&x.\nonumber
\end{eqnarray}
We shall prove the following proposition.
%
\begin{prop}\label{prop4.1}
$-\chi(\gamma)=\rho(\gamma)$ and
%
\begin{equation}\label{eq4.13}
\rho(\gamma)=\lim_{T\to\infty}\frac{1}{T}E\biggl[\int_0^T\varphi
({\bar Y}_s,{\bar Z}_s)\,ds\biggr],
\end{equation}
where ${\bar Z}_s={\hat z}({\bar Y}_s)$.
\end{prop}

For the proof of this proposition, we prepare the following lemma.
%
\begin{lem}\label{lem4.2} Under assumptions (\ref{eq2.11}), (\ref{eq2.12}), (\ref{eq2.20}) and (\ref{eq2.22})
the following estimates hold. For each $\gamma_1<\gamma_0<0$ there
exist positive constants $\delta>0$ and $C>0$ independent of $T$ and
$\gamma$ with $\gamma_1\leq\gamma\leq\gamma_0$ such that
%
\begin{equation}\label{eq4.14}
E\bigl[e^{\delta{\bar w}({\bar Y}_T)}\bigr]\leq C\
\end{equation}
and also
%
\begin{equation}\label{eq4.15}
E\bigl[e^{\delta|{\bar Y}_T|^2}\bigr]\leq C.
\end{equation}
\end{lem}

\begin{pf} Let us set
%
\begin{eqnarray}\label{eq4.16}
{\bar L}\psi&=&\tfrac{1}{2}\operatorname{tr}[\lambda\lambda^*D^2\psi
]+(G+\lambda{\hat z})^*D\psi \nonumber
\\[-8pt]
\\[-8pt]
&=&\tfrac{1}{2}\operatorname{tr}[\lambda\lambda^*D^2\psi]+(\beta_{\gamma
}-\lambda N_{\gamma}^{-1}\lambda^*D{\bar w})^*D\psi.\nonumber
\end{eqnarray}
Then we have
\begin{eqnarray*}
-\chi(\gamma)&=&{\bar L}{\bar w}+\varphi(x,{\hat z}(x)) \\
&=&{\bar L}{\bar w}+\frac{1}{2}(\lambda^*D{\bar w}-\Sigma^*{\hat
\alpha})^*N_{\gamma}^{-1}(\lambda^*D{\bar w}-\Sigma^*{\hat\alpha
}) \\
&&{} +(\lambda^*D{\bar w}-\Sigma^*{\hat\alpha
})^*N^{-1}\Sigma^*{\hat\alpha}+\frac{1}{2}{\hat\alpha}\Sigma
\Sigma^*{\hat\alpha} \\
&=&{\bar L}{\bar w}+\frac{1}{2}(\lambda^* D{\bar w})^*N_{\gamma
}^{-1}\lambda^*D{\bar w}-\frac{\gamma}{2(1-\gamma)}{\hat\alpha
}\Sigma\Sigma^*{\hat\alpha}.
\end{eqnarray*}
Therefore, by applying It\^o's formula, we have
\begin{eqnarray*}
 e^{\delta{\bar w}({\bar Y}_t)}-e^{\delta{\bar w}({\bar
Y}_0)}&=&\delta\int_0^t\biggl\{ {\bar L}{\bar w}({\bar Y}_s)+\frac{\delta
}{2} (D{\bar w})^*\lambda\lambda^*D{\bar w}\biggr\}e^{\delta{\bar
w}}({\bar Y}_s)\,ds \\
&& {}+\delta\int_0^te^{\delta{\bar w}}(D {\bar
w})^*\lambda({\bar Y}_s)\,dW_s \\
&=&\delta\int_0^t\biggl\{-\chi-\frac{1}{2}(D{\bar w})^*\lambda N_{\gamma
}^{-1}\lambda^*D{\bar w}\\
&&\hphantom{\delta\int_0^t\biggl\{}{}+\frac{\gamma}{2(1-\gamma)}{\hat\alpha
}\Sigma\Sigma^*{\hat\alpha}+\frac{\delta}{2}(D{\bar w})^*\lambda
\lambda^*D{\bar w}\biggr\}e^{\delta{\bar w}}({\bar Y}_s)\,ds \\
&&{} +\delta\int_0^te^{\delta{\bar w}}(D{\bar w})^*\lambda
({\bar Y}_s)\,dW_s.
\end{eqnarray*}
Thus, for $p>0$,
\begin{eqnarray*}
 d\bigl(e^{\delta{\bar w}({\bar Y}_t)}e^{p\delta t}\bigr)&=&e^{p\delta
t}\,de^{\delta{\bar w}({\bar Y}_t)}+p\delta e^{p\delta t}e^{\delta{\bar
w}({\bar Y}_t)}\,dt \\
&=&e^{p\delta t}\delta\biggl\{-\chi-\frac{1}{2}(D{\bar w})^*\lambda
N_{\gamma}^{-1}\lambda^* D{\bar w}\\
&&\hphantom{e^{p\delta t}\delta\biggl\{}{}+\frac{\gamma}{2(1-\gamma)}{\hat
\alpha}\Sigma\Sigma^*{\hat\alpha}+\frac{\delta}{2} (D{\bar
w})^*\lambda\lambda^*D{\bar w}+ p\biggr\}e^{\delta{\bar w}}({\bar
Y}_t)\,dt \\
&& {}+\delta e^{p\delta t}e^{\delta{\bar w}({\bar
Y}_t)}(D{\bar w})^*\lambda({\bar Y}_t)\,dW_t.
\end{eqnarray*}
Taking into account (\ref{eq2.22}) and (\ref{eq2.17}), for $\delta>0$ such that
$\delta<\frac{c_1}{1-\gamma}$, we have
\begin{eqnarray*}
&&-\chi-\frac{1}{2}(D{\bar w})^*\lambda N_{\gamma}^{-1}\lambda
^*D{\bar w}+\frac{\gamma}{2(1-\gamma)}{\hat\alpha}\Sigma\Sigma
^*{\hat\alpha}+\frac{\delta}{2} (D{\bar w})^*\lambda\lambda
^*D{\bar w}+p\\
&& \qquad \leq-k_1|x|^2+k_2
\end{eqnarray*}
for $k_1, k_2>0$. Thus, we obtain
\begin{eqnarray*}
 e^{\delta{\bar w}({\bar Y}_t)+p\delta t}&\leq& e^{\delta{\bar
w}(x)}+\delta\int_0^te^{p\delta s +\delta{\bar w}({\bar Y}_s)}\{
-k_1|{\bar Y}_s|^2+k_2\}\,ds\\
&&   {}+\delta\int_0^te^{p\delta s+\delta{\bar
w}({\bar Y}_s)}(D {\bar w})^*\lambda({\bar Y}_s)\,dW_s.
\end{eqnarray*}
Therefore, taking
\[
\tau=\tau_R:=\inf\{t;|{\bar Y}_t|\geq R\},
\]
and setting
\[
k_3=\sup_{|y|\leq\sqrt{k_2/k_1}} {\bar w}(y),
\]
we see that
\begin{eqnarray*}
 E\bigl[e^{\delta{\bar w}({\bar Y}_{t\wedge\tau})+p\delta(t\wedge\tau
)}\bigr]&\leq& e^{\delta{\bar w}(x)}+\delta E\biggl[\int_0^{t\wedge\tau
}e^{p\delta s+\delta{\bar w}({\bar Y}_s)}\{-k_1|{\bar Y}_s|^2+k_2\}\,ds\biggr] \\
& \leq& e^{\delta{\bar w}(x)}+\delta k_2E\biggl[\int_0^{t\wedge
\tau}e^{p\delta s+\delta{\bar w}({\bar Y}_s)}{\mathbf1}_{\{|{\bar
Y}_s|^2\leq\fraca{k_2}{k_1}\}}\,ds\biggr] \\
& \leq& e^{\delta{\bar w}(x)}+\delta k_2e^{\delta
k_3}E\biggl[\int_0^{t\wedge\tau}e^{p\delta s}\,ds\biggr] \\
&  =&e^{\delta{\bar w}(x)}+ k_2e^{\delta k_3} E\biggl[\frac
{1}{p}\bigl(e^{p\delta(t\wedge\tau)}-1\bigr)\biggr].
\end{eqnarray*}
By letting $R$ tend to $\infty$, we have
\[
E\bigl[e^{\delta{\bar w}({\bar Y}_t)+p\delta t}\bigr]\leq e^{\delta{\bar
w}(x)}+ k_2e^{\delta k_3}\frac{1}{p}(e^{p\delta t}-1).
\]
Hence,
\begin{eqnarray*}
E\bigl[e^{\delta{\bar w}({\bar Y}_t)}\bigr]&\leq& e^{-p\delta t+\delta{\bar
w}(x)}+ k_2e^{\delta k_3}\frac{1}{p}(1-e^{-p\delta t}) \\
&\leq& e^{\delta{\bar w}(x)}+ k_2e^{\delta k_3}\frac{1}{p}.
\end{eqnarray*}

Finally, we see that (\ref{eq4.15}) follows from (\ref{eq4.14}) because of (\ref{eq3.12}).
\end{pf}

\begin{pf*}{Proof of Proposition~\ref{prop4.1}}
 From (\ref{eq4.10}) it follows that
\[
-\chi(\gamma)\leq\tfrac{1}{2}\operatorname{tr}[\lambda\lambda^*D^2{\bar
w}]+G^*D{\bar w}+(\lambda z)^*D{\bar w}+\varphi(x,z)\vadjust{\goodbreak}
\]
for each $z\in R^{n+m}$. Therefore, for each control $Z_t$ we have
\begin{eqnarray*}
 {\bar w}(Y_t)-{\bar w}(x)&=&\int_0^T( D{\bar w}(Y_s))^*\lambda
(Y_s)\,dW_s+\int_0^T\{G(Y_s)+\lambda(Y_s)Z_s\}^*D{\bar
w}(Y_s)\,ds \\
&&{} +\frac{1}{2}\int_0^T\operatorname{tr}[\lambda\lambda
^*D^2{\bar w}](Y_s)\,ds \\
&\geq&\int_0^T(D{\bar w}(Y_s))^*\lambda(Y_s)\,dW_s-\chi T-\int
_0^T\varphi(Y_s,Z_s)\,ds.
\end{eqnarray*}
Thus,
\[
{\bar w}(x)-\chi T\leq E\biggl[\int_0^T\varphi(Y_s,Z_s)\,ds+{\bar w}(Y_T)\biggr],
\]
from which we obtain
\begin{eqnarray*}
-\chi&\leq&\limsup_{T\to\infty}\frac{1}{T}E\biggl[\int_0^T\varphi
(Y_s,Z_s)\,ds+{\bar w}(Y_T)\biggr] \\
&=&\limsup_{T\to\infty}\frac{1}{T}E\biggl[\int_0^T\varphi(Y_s,Z_s)\,ds\biggr]
\end{eqnarray*}
since $|{\bar w}(y)|^2\leq c|y|^2+c'$. Namely, we have $-\chi(\gamma
)\leq\rho(\gamma)$. On the other hand, by taking $Z_t={\bar Z}_t$,
we have
\[
{\bar w}({\bar Y}_t)-{\bar w}(x)=\int_0^T( D{\bar w}({\bar
Y}_s))^*\lambda({\bar Y}_s)\,dW_s-\chi T-\int_0^T\varphi({\bar
Y}_s,{\bar Z}_s)\,ds,
\]
and thus
\[
{\bar w}(x)-\chi T=E\biggl[\int_0^T\varphi({\bar Y}_s,{\bar Z}_s)\,ds+{\bar
w}({\bar Y}_T)\biggr].
\]
Lemma~\ref{lem4.2} implies
\[
-\chi=\limsup_{T\to\infty}\frac{1}{T}E\biggl[\int_0^T\varphi({\bar
Y}_s,{\bar Z}_s)\,ds\biggr]\geq\rho,
\]
and we see that $-\chi(\gamma)=\rho(\gamma)$.
\end{pf*}

Let us define
%
\begin{equation}\label{eq4.17}
 \qquad {\bar\chi}(\gamma)=\limsup_{T\to\infty}\frac{1}{T}\inf_{Z\in
{\tilde{\mathcal A}}}E\biggl[\int_0^T\varphi(Y_s,Z_s)\,ds\biggr]=\limsup_{T\to
\infty}\frac{1}{T}{\bar v}(0,x;T).
\end{equation}
Then we can see that
\[
{\bar\chi}\leq\rho(\gamma)=-\chi(\gamma).
\]
%
\begin{prop}\label{prop4.2}
Assume (\ref{eq2.11}), (\ref{eq2.12}), (\ref{eq2.20}) and (\ref{eq2.22}). Then
\[
{\bar\chi}(\gamma)=\rho(\gamma)=-\chi(\gamma).
\]
\end{prop}

The proof of Theorem~\ref{theo2.2} follows directly from this proposition since
${\bar\chi}(\gamma)=-{\hat\chi}(\gamma)$ because of Proposition~\ref{prop2.1}.

For the proof of the present proposition, we prepare some lemmas.

For each $T>0$, we take the controlled process ${\hat Y}_t={\hat
Y}_t^{(T)}$ defined by (\ref{eq4.8}) and control ${\hat z}(t,{\hat
Y}_t^{(T)})$. Taking a sequence $\{T_n\}$ such that
\[
{\bar\chi}=\lim_{T_n\to\infty}\frac{1}{T_n}{\bar v}(0,x;T_n)=\lim
_{T_n\to\infty}\frac{1}{T_n}E\biggl[\int_0^{T_n}\varphi\bigl({\hat
Y}^{(T_n)}_t,{\hat z}\bigl(t,{\hat Y}_t^{(T_n)}\bigr)\bigr)\,dt\biggr],
\]
we have the following lemma.
%
\begin{lem} Under the assumptions of Proposition~\ref{prop4.2},
for each $t>0$ we have
%
\begin{equation}\label{eq4.18}
\liminf_{T_n\to\infty}\frac{1}{T_n}E\bigl[\bigl|{\hat
Y}_{T_n-t}^{(T_n)}\bigr|^2\bigr]=0.
\end{equation}
\end{lem}

\begin{pf} Set
\[
{\hat L}\psi:=\tfrac{1}{2}\operatorname{tr}[\lambda\lambda^*D^2\psi
]+\bigl(G+\lambda{\hat z}(t,x)\bigr)^*D\psi.
\]
Then
\begin{eqnarray*}
&&{\bar v}(T-t,{\hat Y}_{T-t};T)-{\bar v}(0,x;T) \\
&& \qquad =\int_0^{T-t}\biggl(\frac{\partial{\bar v}}{\partial t}+{\hat L}(s,{\hat
Y}_s)\biggr)\,ds+\int_0^{T-t}(D {\bar v})(s,{\hat Y}_s)\,dW_s \\
&& \qquad =-\int_0^{T-t}\varphi({\hat Y}_s,{\hat z}(s,{\hat Y}_s))\,ds+\int
_0^{T-t}(D{\bar v}(s,{\hat Y}_s))^*\lambda({\hat Y}_s)\,dW_s.
\end{eqnarray*}
Therefore,
\[
{\bar v}(0,x;T_n)=E\biggl[\int_0^{T_n-t}\varphi\bigl({\hat Y}_s^{(T_n)},{\hat
z}\bigl(s,{\hat Y}_s^{(T_n)}\bigr)\bigr)\,ds+{\bar v}\bigl(T_n-t,{\hat Y}_{T_n-t}^{(T_n)};T_n\bigr)\biggr].
\]
Since
\[
\limsup_{T_n\to\infty}\frac{1}{T_n}E\biggl[\int_0^{T_n-t}\varphi\bigl({\hat
Y}_s^{(T_n)},{\hat z}\bigl(s,{\hat Y}_s^{(T_n)}\bigr)\,ds\biggr]\geq{\bar\chi},
\]
we have
%
\begin{equation}\label{eq4.19}
\liminf_{T_n\to\infty}\frac{1}{T_n}E\bigl[{\bar v}\bigl(T_n-t,{\hat
Y}_{T_n-t}^{(T_n)};T_n\bigr)\bigr]=0.
\end{equation}
Noting that ${\bar v}(T_n-t,x;T_n)\geq k|x|^2-k',   k=k(t)$ and
$k'=k'(t)$, because of Lemma~\ref{lem4.1}, we obtain
\[
0\leq\liminf_{T_n\to\infty}\frac{1}{T_n}kE\bigl[\bigl|{\hat
Y}_{T_n-t}^{(T_n)}\bigr|^2\bigr]\leq0,
\]
and our lemma has been proved.
\end{pf}

%
\begin{lem}
Under the assumptions of Proposition~\ref{prop4.2}, there exists a~subsequence $\{
T_n'\}\subset\{T_n\}$ such that
\[
\lim_{T_n'\to\infty}\frac{1}{T_n'}E\bigl[\bigl|{\hat Y}_{T_n'}^{(T_n')}\bigr|^2\bigr]=0.
\]
\end{lem}

\begin{pf}
\begin{eqnarray*}
 \bigl|{\hat Y}_T^{(T)}\bigr|^2-\bigl|{\hat Y}_{T-t}^{(T)}\bigr|^2&=&2\int_{T-t}^T\bigl({\hat
Y}_s^{(T)}\bigr)^*\lambda\bigl({\hat Y}^{(T)}_s\bigr)\,dW_s \\[-2pt]
&&{}+2\int_{T-t}^T{\hat Y}^{(T)}_s\bigl\{G\bigl({\hat Y}^{(T)}\bigr)+\lambda\bigl({\hat
Y}_s^{(T)}\bigr){\hat z}\bigl(s,{\hat Y}_s^{(T)}\bigr)\bigr\}\,ds\\[-2pt]
&&{} +\int_{T-t}^T\operatorname
{tr}\bigl[\lambda\lambda\bigl({\hat Y}_s^{(T)}\bigr)\bigr]\,ds.
\end{eqnarray*}
Therefore,
\begin{eqnarray*}
 E\bigl[\bigl|{\hat Y}_T^{(T)}\bigr|^2\bigr]&=&E\bigl[\bigl|{\hat Y}_{T-t}^{(T)}\bigr|^2\bigr]+2E\biggl[\int
_{T-t}^T{\hat Y}^{(T)}_s\bigl\{G\bigl({\hat Y}^{(T)}\bigr)+\lambda\bigl({\hat
Y}_s^{(T)}\bigr){\hat z}\bigl(s,{\hat Y}_s^{(T)}\bigr)\bigr\}\,ds  \\[-2pt]
&&\hphantom{E\bigl[\bigl|{\hat Y}_{T-t}^{(T)}\bigr|^2\bigr]+2E\biggl[}\hspace*{91pt}{}
+\int_{T-t}^T\operatorname{tr}\bigl[\lambda\lambda\bigl({\hat Y}^{(T)}_s\bigr)\bigr]\,ds\biggr].
\end{eqnarray*}
By using the gradient estimates in Theorem~\ref{theo2.1} and (\ref{eq2.20}) we obtain
\[
y^*G(y)+y^*\lambda(y){\hat z}(s,y)\leq c(|y|^2+1)
\]
for some positive constant $c$ and
\begin{eqnarray*}
E\bigl[\bigl|{\hat Y}_T^{(T)}\bigr|^2\bigr] &\leq&E\bigl[\bigl|{\hat Y}_{T-t}^{(T)}\bigr|^2\bigr]+cE\biggl[\int
_{T-t}^T\bigl|{\hat Y}_s^{(T)}\bigr|^2\,ds\biggr] \\[-2pt]
&\leq&E\bigl[\bigl|{\hat Y}_{T-t}^{(T)}\bigr|^2\bigr]+c'E\biggl[\int_{T-t}^T\varphi\bigl({\hat
Y}_s^{(T)},{\hat z}\bigl(s,{\hat Y}_s^{(T)}\bigr)\bigr)\,ds\biggr]
\end{eqnarray*}
since (\ref{eq2.22}) is assumed and
\begin{eqnarray*}
\varphi(x,{\hat z}(t,x))&=&\frac{1}{2}{\hat z}(t,x)^* N_{\gamma}
{\hat z}(t,x)-{\hat z}(t,x)^*\Sigma^*{\hat\alpha}+\frac{1}{2}{\hat
\alpha}\Sigma\Sigma^*{\hat\alpha} \\[-2pt]
&=&\frac{1}{2}(\lambda^*D{\bar v}-\Sigma^*{\hat\alpha})^*
N_{\gamma}( \lambda^*D{\bar v}-\Sigma^*{\hat\alpha}) \\[-2pt]
&&{} +(\lambda^*D{\bar v}-\Sigma^*{\hat\alpha
})^*N_{\gamma}^{-1}\Sigma^*{\hat\alpha}+\frac{1}{2}{\hat\alpha
}\Sigma\Sigma^*{\hat\alpha} \\[-2pt]
&=&\frac{1}{2}(D{\bar v})^*\lambda N_{\gamma}^{-1}\lambda^*D{\bar
v}-\frac{\gamma}{2(1-\gamma)}{\hat\alpha}\Sigma\Sigma^*{\hat
\alpha}.
\end{eqnarray*}
By It\^o's formula,
\begin{eqnarray*}
&&{\bar v}(T,{\hat Y}_T;T)-{\bar v}(T-t,{\hat Y}_{T-t};T) \\
&& \qquad =-\int_{T-t}^T\varphi({\hat Y}_s,{\hat z}(s,{\hat Y}_s))\,ds+\int
_{T-t}^T(D{\bar v}(s,{\hat Y}_s))^*\lambda({\hat Y}_s)\,dW_s,
\end{eqnarray*}
and we have
\[
E\biggl[\int_{T-t}^T\varphi({\hat Y}_s,{\hat z}(s,{\hat Y}_s))\,ds\biggr]=
E[{\bar v}(T-t,{\hat Y}_{T-t};T)].
\]
Take a subsequence $\{T_n'\}\subset\{T_n\}$ such that
\[
\lim_{T_n'\to\infty}\frac{1}{T_n'}E[{\bar v}(T_n'-t,{\hat
Y}_{T_n'-t};T_n')]=0
\]
and
\[
\lim_{T_n'\to\infty}\frac{1}{T_n'}E\bigl[\bigl|{\hat Y}_{T_n'-t}^{(T_n')}\bigr|^2\bigr]=0.
\]
Then we have
\begin{eqnarray*}
 0\leq\limsup_{T_n'\to\infty}\frac{1}{T_n'}E\bigl[\bigl|{\hat
Y}_{T_n'}^{(T_n')}\bigr|^2\bigr]&\leq&\lim_{T_n'\to\infty}\frac
{1}{T_n'}E\bigl[\bigl|{\hat Y}_{T_n'-t}^{(T_n')}\bigr|^2\bigr] \\
&&{}
+c'\lim_{T_n'\to\infty}\frac{1}{T_n'}E\biggl[\int_{T_n'-t}^{T_n'}\varphi
\bigl({\hat Y}_s^{(T_n')},{\hat z}\bigl(s,{\hat
Y}_s^{(T_n')}\bigr)\bigr)\,ds\biggr]=0,
\end{eqnarray*}
and thus the lemma has been proved.
\end{pf}

\begin{pf*}{Proof of Proposition~\ref{prop4.2}}
 For each $\epsilon$ there exists
$T_{\epsilon}$ such that
\begin{eqnarray*}
&\displaystyle E\bigl[\bigl|{\hat
Y}_{T_{\epsilon}}^{(T_{\epsilon})}\bigr|^2\bigr]< \epsilon T_{\epsilon
},\qquad{\bar w}(x)< \epsilon T_{\epsilon},&
\\
&\displaystyle
\biggl|{\bar\chi}-\frac{1}{T_{\epsilon}}E\biggl[\int_0^{T_{\epsilon}}\varphi
\bigl({\hat Y}_s^{(T_{\epsilon})},{\hat
Z}_s^{(T_{\epsilon})}\bigr)\,ds\biggr]\biggr|<\epsilon.&
\end{eqnarray*}
Set
\[
Z_s^{(\epsilon)}=\cases{\displaystyle
{\hat Z}_s^{(T_{\epsilon})}, & \quad  $s<T_{\epsilon}$,\cr\displaystyle
0,& \quad  $T_{\epsilon}\leq s$,
}
\]
and consider $Y_t^{(\epsilon)}$ defined by
\[
dY_t^{(\epsilon)}=\lambda\bigl(Y_t^{(\epsilon)}\bigr)\,dW_t+\bigl\{G\bigl(Y_t^{(\epsilon
)}\bigr)+\lambda\bigl(Y_t^{(\epsilon)}\bigr)Z_t^{(\epsilon)}\bigr\}\,dt, \qquad   Y_0^{(\epsilon)}=x.
\]
Then, for $t\geq T_{\epsilon}$,
\[
Y_t^{(\epsilon)}=Y_{T_{\epsilon}}^{(\epsilon)}+\int_{T_{\epsilon
}}^t\lambda\bigl(Y_s^{(\epsilon)}\bigr)\,dW_s+\int_{T_{\epsilon
}}^tG\bigl(Y_s^{(\epsilon)}\bigr)\,ds
\]
and
\begin{eqnarray*}
d\bigl|Y_t^{(\epsilon)}\bigr|^2&=&2\bigl(Y_t^{(\epsilon)}\bigr)^*\lambda\bigl(Y_t^{(\epsilon
)}\bigr)\,dW_t+2\bigl(Y_t^{(\epsilon)}\bigr)^*G\bigl(Y_t^{(\epsilon)}\bigr)+\int_{T_{\epsilon
}}^tG\bigl(Y_t^{(\epsilon)}\bigr)\,dt\\
&&{}+\operatorname{tr}\bigl[\lambda\lambda^*\bigl(Y_t^{(\epsilon)}\bigr)\bigr]\,dt.
\end{eqnarray*}
Therefore,
\begin{eqnarray*}
 &&e^{pt}\bigl|Y_t^{(\epsilon)}\bigr|^2\\
 && \qquad =e^{p T_{\epsilon}}\bigl|Y_{T_{\epsilon
}}^{(\epsilon)}\bigr|^2+2\int_{T_{\epsilon}}^te^{ps}\bigl(Y_s^{(\epsilon
)}\bigr)^*\lambda\bigl(Y_s^{(\epsilon)}\bigr)\,dW_s \\
&&  \quad \qquad {}+\int_{T_{\epsilon}}^te^{ps}\bigl\{2\bigl(Y_s^{(\epsilon
)}\bigr)^*G\bigl(Y_s^{(\epsilon)}\bigr)+p\bigl|Y_s^{(\epsilon)}\bigr|^2+\operatorname{tr}\bigl[\lambda
\lambda^*\bigl(Y_s^{(\epsilon)}\bigr)\bigr]\bigr\}\,ds \\
&& \qquad \leq e^{p T_{\epsilon}}\bigl|Y_{T_{\epsilon}}^{(\epsilon)}\bigr|^2+2\int
_{T_{\epsilon}}^te^{ps}\bigl(Y_s^{(\epsilon)}\bigr)^*\lambda\bigl(Y_s^{(\epsilon
)}\bigr)\,dW_s+\int_{T_{\epsilon}}^te^{ps}\bigl\{-k_1\bigl|Y_s^{(\epsilon)}\bigr|^2+k_2\bigr\}\,ds
\end{eqnarray*}
for some positive constants $k_1, k_2>0$. By using stopping time
arguments, for $t\geq T_{\epsilon}$,
we have
\begin{eqnarray*}
E\bigl[e^{pt}\bigl|Y_t^{(\epsilon)}\bigr|^2\bigr]&\leq& E[e^{pT_{\epsilon
}}|Y_{T_{\epsilon}}|^2]+E\biggl[\int_{T_{\epsilon}}^te^{ps}k_2 \,ds\biggr] \\
&=& E\bigl[e^{pT_{\epsilon}}\bigl|Y_t^{(\epsilon)}\bigr|^2\bigr]+\frac
{k_2}{p}(e^{pt}-e^{pT_{\epsilon}}).
\end{eqnarray*}
Thus, we see that
\[
E\bigl[\bigl|Y_t^{(\epsilon)}\bigr|^2\bigr]\leq E\bigl[\bigl|Y_{T_{\epsilon}}^{(\epsilon
)}\bigr|^2\bigr]+\frac{k_2}{p},
\]
from which we obtain $\limsup_{t\to\infty}\frac
{1}{t}E[|Y_t^{(\epsilon)}|^2]=0$.
Hence, $Z^{(\epsilon)}\in{\tilde{\mathcal A}}_{\infty}$. Now, by
applying It\^o's formula, we have
\begin{eqnarray*}
 {\bar w}\bigl(Y_{T_{\epsilon}}^{(\epsilon)}\bigr)-{\bar w}(x)&=&\int
_0^{T_{\epsilon}}\bigl(D{\bar w}\bigl(Y_s^{(\epsilon)}\bigr)\bigr)^*\lambda
\bigl(Y_s^{(\epsilon)}\bigr)\,dW_s \\
&&{}+\int_0^{T_{\epsilon}}\bigl\{G\bigl(Y_s^{(\epsilon)}\bigr)+\lambda
\bigl(Y_s^{(\epsilon)}\bigr)Z_s^{(\epsilon)}\bigr\}^*D{\bar w}\bigl(Y_s^{(\epsilon)}\bigr)\,ds\\
&&{}+
\frac{1}{2}\int_0^{T_{\epsilon}}\operatorname{tr}[\lambda\lambda
^*D^2{\bar w}]\bigl(Y_s^{(\epsilon)}\bigr)\,ds \\
&
\geq&\int_0^{T_{\epsilon}}\bigl(D{\bar w}\bigl(Y_s^{(\epsilon)}\bigr)\bigr)^* \lambda
\bigl(Y_s^{(\epsilon)}\bigr)\,dW_s-\chi T_{\epsilon}-\int_0^{T_{\epsilon
}}\varphi\bigl(Y_s^{(\epsilon)},Z_s^{(\epsilon)}\bigr)\,ds.
\end{eqnarray*}
Therefore,
\[
{\bar w}(x)-\chi T_{\epsilon} \leq E\biggl[\int_0^{T_{\epsilon}}\varphi
\bigl(Y_s^{(\epsilon)},Z_s^{(\epsilon)}\bigr)\,ds+{\bar w}\bigl(Y_{T_{\epsilon
}}^{(\epsilon)}\bigr)\biggr],
\]
from which we have
\begin{eqnarray*}
-\chi&\leq&\frac{1}{T_{\epsilon}}E\biggl[\int_0^{T_{\epsilon}}\varphi
\bigl(Y_s^{(\epsilon)},Z_s^{(\epsilon)}\bigr)\,ds\biggr]+\frac{1}{T_{\epsilon
}}E\bigl[{\bar w}\bigl(Y_{T_{\epsilon}}^{(\epsilon)}\bigr)\bigr] \\
&\leq&\epsilon+{\bar\chi}+c\epsilon
\end{eqnarray*}
for some positive constant $c>0$. Therefore, $-\chi\leq{\bar\chi
}+c\epsilon$ for any $\epsilon$, and we have $-\chi\leq{\bar\chi
}$. This completes the proof of the proposition.
\end{pf*}

The following is a direct consequence of Proposition~\ref{prop4.1}.
%
\begin{cor}\label{cor4.1}
Under the assumptions of Proposition~\ref{prop4.2}, $\rho(\gamma)$ is a~concave
function on $(-\infty,0)$, and ${\hat\chi}(\gamma)$ is a convex function.
\end{cor}

Indeed,
\[
\varphi=\frac{1}{2}z^*z-\frac{\gamma}{2}z^*\sigma^*(\sigma\sigma
^*)^{-1}\sigma z-z^*\Sigma^*{\hat\alpha}+\frac{1}{2}{\hat\alpha
}\Sigma\Sigma^*{\hat\alpha}
\]
is a concave function of $\gamma$, and the infimum of a family of
concave functions~$\rho(\gamma)$ is concave.
%
\begin{prop}\label{prop4.3}
Under the assumptions of Proposition~\ref{prop3.1}, ${\bar L}$ is ergodic.
\end{prop}

\begin{pf}
\[
{\bar L}{\bar w}=-\frac{1}{2}(D{\bar w})^*\lambda N_{\gamma
}^{-1}\lambda^*D{\bar w}+\frac{\gamma}{2(1-\gamma)}{\hat\alpha
}^*\Sigma\Sigma^*{\hat\alpha}-\chi\to-\infty
\]
as $|x|\to\infty$, and ${\bar L}{\bar w}\leq-c,   |x|\gg1$ and
$c>0$. Moreover, ${\bar w}(x)\to\infty,   |x|\to\infty$, and the
Hasminskii conditions (cf. \cite{Has}) hold.
\end{pf}

%
\section{Derived Poisson equation}\label{sec5}
We are going to consider a Poisson equation formally obtained by
differentiating H-J-B equation (\ref{eq3.1}) of ergodic type with respect to
$\gamma$. Namely, we consider
\begin{eqnarray*}
 -\theta(\gamma)&=&\frac{1}{2}\operatorname{tr}[\lambda\lambda^* D^2
u]+G^*Du-(\lambda^*D{\bar w}-\Sigma^*{\hat\alpha})^*N_{\gamma
}^{-1}\lambda^*Du \\
&&{}  -\frac{1}{2(1-\gamma)^2}(\lambda^*D{\bar w}-\Sigma
^*{\hat\alpha})^*\Sigma^*(\Sigma\Sigma^*)^{-1}\Sigma(\lambda
^*D{\bar w}-\Sigma^*{\hat\alpha}).
\end{eqnarray*}
Since
\begin{eqnarray*}
&&-\frac{1}{2(1-\gamma)^2}(\lambda^*D{\bar w}-\Sigma^*{\hat\alpha
})^*\Sigma^*(\Sigma\Sigma^*)^{-1}\Sigma(\lambda^*D{\bar w}-\Sigma
^*{\hat\alpha}) \\
&& \qquad  =-\frac{1}{2(1-\gamma)^2}(\sigma\lambda^*D{\bar
w}-{\hat\alpha})^*(\sigma\sigma^*)^{-1} (\sigma\lambda^*D{\bar
w}-{\hat\alpha}),
\end{eqnarray*}
we can write
%
\begin{equation}\label{eq5.1}\quad
-\theta(\gamma)={\bar L}u-\frac{1}{2(1-\gamma)^2}(\sigma\lambda
^*D{\bar w}-{\hat\alpha})^*(\sigma\sigma^*)^{-1} (\sigma\lambda
^*D{\bar w}-{\hat\alpha}).
\end{equation}
Note that ${\bar L}$ is ergodic in light of Proposition~\ref{prop4.3}, and the
pair $(u,\theta(\gamma))$ of a function $u$ and a constant $\theta
(\gamma)$ is the solution to (\ref{eq5.1}). Let us set
\[
{\mathcal D}=B_{R_0}=\{x\in R^n; |x|< R_0\},
\]
and let $R_0$ be sufficiently large so that
%
\begin{eqnarray}\label{eq5.2}
K(x;{\bar w}):=\frac{1}{2}(D{\bar w})^*\lambda N_{\gamma}^{-1}\lambda
^* D{\bar w}-\frac{\gamma}{2(1-\gamma)}{\hat\alpha}^*\Sigma\Sigma
^*{\hat\alpha}+\chi>0,\nonumber
\\[-8pt]
\\[-8pt]
\eqntext{\displaystyle x\in{\mathcal D}^c,}
\end{eqnarray}
for $\gamma\leq\gamma_0<0$, which is possible because of assumption (\ref{eq2.22}).
Therefore, we see that ${\bar L}$ and ${\bar w}$ satisfy assumption
(\ref{eq9.3}) in the \hyperref[appm]{Appendix} and also
\[
\sup_{x\in{\mathcal D}^c}\frac{|f^{(\gamma)}(x)|}{K(x\dvtx {\bar
w})}<\infty,
\]
for
\[
f^{(\gamma)}=-\frac{1}{2(1-\gamma)^2}(\sigma\lambda^*D{\bar
w}-{\hat\alpha})^*(\sigma\sigma^*)^{-1} (\sigma\lambda^*D{\bar
w}-{\hat\alpha}).
\]
In the following we always take a solution ${\bar w}$ to (\ref{eq3.1}) such
that ${\bar w}(x)>0$. Thus, according to Proposition~\ref{prop9.4} we can show
the existence of the solution $(u,\theta(\gamma))$ to (\ref{eq5.1}).
%
\begin{cor}\label{cor5.1}
Equation (\ref{eq5.1}) has a solution $(u, \theta(\gamma))$ such that
\[
\sup_{x\in{\mathcal D}^c}\frac{|u|}{{\bar w}}<\infty,\qquad u\in
W^{2,p}_{\mathrm{loc}},
\]
and
\[
\theta(\gamma)=\int\frac{1}{2(1-\gamma)^2}(\sigma\lambda^*D{\bar
w}-{\hat\alpha})^*(\sigma\sigma^*)^{-1}(\sigma\lambda^*D{\bar
w}-{\hat\alpha})m_{\gamma}(y)\,dy.
\]
Moreover, this solution $u$ is unique up to additive constants.
\end{cor}

\begin{pf} It can be clearly seen that
\[
\frac{1}{2(1-\gamma)^2}(\sigma\lambda^*D{\bar w}-{\hat\alpha
})^*(\sigma\sigma^*)^{-1}(\sigma\lambda^*D{\bar w}-{\hat\alpha
})\in F_K
\]
and Proposition~\ref{prop9.4} applies.
\end{pf}

%
\section{Differentiability of H-J-B equation}\label{sec6}
%
\begin{lem}\label{lem6.1}\label{lem2.6}
Under the assumptions of Proposition~\ref{prop4.2},
%
\begin{equation}\label{eq6.1}
\int e^{\delta|x|^2}m_{\gamma}(dx)\leq c,
\end{equation}
where $c$ and $\delta$ are positive constants independent of $\gamma
_1\leq\gamma\leq\gamma_0<0$.
\end{lem}

\begin{pf} Inequality (\ref{eq6.1}) is a direct consequence of (\ref{eq4.15}) in Lemma~\ref{lem4.2} since
${\bar Y}_t$ is an ergodic diffusion process with the invariant measure
$m_{\gamma}(dx)$.
\end{pf}

In the following, we always work under the assumptions of
Theorem~\ref{theo2.2}
(Proposition~\ref{prop4.2}).
%
\begin{lem}\label{lem6.2} Let $({\bar w}^{(\gamma)},\chi(\gamma))$ and $  ({\bar
w}^{(\gamma+\Delta)}, \chi(\gamma+\Delta))$ be solutions to~(\ref{eq3.1})
with $\gamma$ and $\gamma+\Delta$, respectively, such that ${\bar
w}^{(\gamma)}(0)={\bar w}^{(\gamma+\Delta)}(0)=c_w>0$. Then
${\bar w}^{(\gamma+\Delta)}$ converges to ${\bar w}^{(\gamma)}$ in
$H^1_{\mathrm{loc}}$ strongly and uniformly on each compact~set.
\end{lem}

\begin{pf} We have
\[
\bigl\|{\bar w}^{(\gamma+\Delta)}\bigr\|_{L^{\infty
}(B_{2r})}\leq2r\bigl\|\nabla{\bar w}^{(\gamma+\Delta)}
\bigr\|_{L^{\infty}(B_{2r})}\leq c_1(\gamma,r)
\]
and
%
\begin{eqnarray}
\bigl|{\bar w}^{(\gamma+\Delta)}(x)-{\bar w}^{(\gamma+\Delta)}(y)\bigr|\leq
|x-y|\bigl\|\nabla{\bar w}^{(\gamma+\Delta)} \bigr\|
_{L^{\infty}(B_{2r})}\leq c_2(\gamma,r), \nonumber\\
\eqntext{\displaystyle  x, y\in B_{2r},}
\end{eqnarray}
for each $r$ in light of (\ref{eq3.11}) and (\ref{eq3.15}), where $c_i(\gamma,r)$ is a
positive constant independent of $\Delta$, $i=1,2$. Therefore it
follows that $\{{\bar w}^{(\gamma+\Delta)}\}_{\Delta}$ is bounded in
$H^1(B_{2r})$ and converges to some ${\tilde w}$ in $H^1(B_{2r})$
weakly for each $r$ and also uniformly on each compact set by taking a
subsequence if necessary. Note that $\chi(\gamma+\Delta)$ converges
to $\chi(\gamma)$ because $\chi(\gamma)$ is convex on $(-\infty
,0)$ and thus continuous. Take a function $\tau\in C_0^{\infty
}(B_{2r})$ such that $\tau(x)\equiv1,  x\in B_r$, and $0\leq\tau
\leq1$.
With $({\bar w}^{(\gamma+\Delta)}-{\tilde w})\tau$, we test
\begin{eqnarray*}
 -\chi(\gamma+\Delta)&=&\tfrac{1}{2}\operatorname{tr}\bigl[\lambda\lambda
^*D^2{\bar w}^{(\gamma+\Delta)}\bigr]+\beta_{\gamma+\Delta}^*D{\bar
w}^{(\gamma+\Delta)} \\
&&{} -\tfrac{1}{2}\bigl(D{\bar w}^{(\gamma+\Delta)}\bigr)^*\lambda
N_{\gamma+\Delta}^{-1}\lambda^*D{\bar w}^{(\gamma+\Delta
)}+U_{\gamma+\Delta}
\end{eqnarray*}
and obtain
\begin{eqnarray*}
&&-\int_{B_{2r}}\chi(\gamma+\Delta)\bigl({\bar w}^{(\gamma+\Delta
)}-{\tilde w}\bigr)\tau \,dx \\
&& \qquad
=-\int_{B_{2r}} \frac{1}{2}(\lambda\lambda^*)^{ij}D_i{\bar
w}^{(\gamma+\Delta)}D_j\bigl(\bigl({\bar w}^{(\gamma+\Delta)}-{\tilde w}\bigr)\tau
\bigr)\,dx \\
&& \qquad  \quad {}  +\int_{B_{2r}}{\tilde\beta}_{\gamma+\Delta
}^*D{\bar w}^{(\gamma+\Delta)}\bigl({\bar w}^{(\gamma+\Delta)}-{\tilde
w}\bigr)\tau \,dx \\
&& \qquad  \quad {}  -\frac{1}{2}\int_{B_{2r}}\bigl(D{\bar w}^{(\gamma+\Delta
)}\bigr)^*\lambda N_{\gamma+\Delta}^{-1}\lambda^*D{\bar w}^{(\gamma
+\Delta)}\bigl({\bar w}^{(\gamma+\Delta)}-{\tilde w}\bigr)\tau \,dx \\
&& \qquad  \quad {} +\int_{B_{2r}}U_{\gamma+\Delta}\bigl({\bar w}^{(\gamma
+\Delta)}-{\tilde w}\bigr)\tau \,dx,
\end{eqnarray*}
where ${\tilde\beta}_{\gamma}^i=\beta_{\gamma}^i-\frac{1}{2}\sum
_j\frac{\partial(\lambda\lambda^*)^{ij}}{\partial x^j}$.
Therefore,
\begin{eqnarray*}
&&\int_{B_{2r}} \frac{1}{2}(\lambda\lambda^*)^{ij}D_i\bigl({\bar
w}^{(\gamma+\Delta)}-{\tilde w}\bigr)D_j\bigl({\bar w}^{(\gamma+\Delta
)}-{\tilde w}\bigr) \tau \,dx \\
&& \qquad =-\int_{B_{2r}} \frac{1}{2}(\lambda\lambda^*)^{ij}D_i{\tilde
w}D_j\bigl({\bar w}^{(\gamma+\Delta)}-{\tilde w}\bigr)\tau \,dx \\
&& \qquad  \quad {} -\int_{B_{2r}} \frac{1}{2}(\lambda\lambda
^*)^{ij}D_i{\bar w}^{(\gamma+\Delta)}D_j\tau\bigl({\bar w}^{(\gamma
+\Delta)}-{\tilde w}\bigr)\,dx \\
&& \qquad  \quad {}  +\int_{B_{2r}}{\tilde\beta}_{\gamma+\Delta
}^*D{\bar w}^{(\gamma+\Delta)}\bigl({\bar w}^{(\gamma+\Delta)}-{\tilde
w}\bigr)\tau \,dx \\
&& \qquad  \quad {}  -\frac{1}{2}\int_{B_{2r}}\bigl(D{\bar w}^{(\gamma+\Delta
)}\bigr)^*\lambda N_{\gamma+\Delta}^{-1}\lambda^*D{\bar w}^{(\gamma
+\Delta)}\bigl({\bar w}^{(\gamma+\Delta)}-{\tilde w}\bigr)\tau \,dx \\
&& \qquad  \quad {} +\int_{B_{2r}}U_{\gamma+\Delta}\bigl({\bar w}^{(\gamma
+\Delta)}-{\tilde w}\bigr)\tau \,dx+\int_{B_{2r}}\chi(\gamma+\Delta
)\bigl({\bar w}^{(\gamma+\Delta)}-{\tilde w}\bigr)\tau \,dx.
\end{eqnarray*}
Since all terms of the right-hand side converge to $0$, we see that
\mbox{$D({\bar w}^{(\gamma+\Delta)}\!-\!{\tilde w})$} converges strongly to $0$
in $L^2(B_r)$ and
${\bar w}^{(\gamma+\Delta)}$ to ${\tilde w}$ strongly in $H^1(B_r)$.
Thus, we obtain our present lemma because $({\bar w},\chi(\gamma))$
satisfies (\ref{eq3.1}), and the solution is unique up to additive constants
with respect to ${\bar w}$.
\end{pf}

%
\begin{lem}\label{lem6.8}\label{lem6.3} Let $(u^{(\gamma+\Delta)}, \theta(\gamma+\Delta))$ be
a solution to
%
\begin{equation}\label{eq6.2}
-\theta(\gamma+\Delta)={\bar L}(\gamma+\Delta)u^{(\gamma+\Delta
)}+f^{(\gamma+\Delta)},
\end{equation}
where
\begin{eqnarray*}
f^{(\gamma+\Delta)}
&=&-\frac{1}{2(1-\gamma-\Delta)^2}\bigl(\lambda
^*D{\bar w}^{(\gamma+\Delta)}-\Sigma^*{\hat\alpha}\bigr)^*\\
&&{}\times \Sigma
^*(\Sigma\Sigma^*)^{-1} \Sigma\bigl(\lambda^*D{\bar w}^{(\gamma+\Delta
)}-\Sigma^*{\hat\alpha}\bigr).
\end{eqnarray*}
Then, as $|\Delta|\to0$, $\theta(\gamma+\Delta)$ converges to
$\theta(\gamma)$ and $u^{(\gamma+\Delta)}$ converges to $u^{(\gamma
)}$ in $H^1_{\mathrm{loc}}$ strongly and uniformly on each compact set, where
$(u^{(\gamma)}, \theta(\gamma))$ is a~solution to (\ref{eq5.1}).
\end{lem}

\begin{pf} Note that $|f^{(\gamma+\Delta)}(x)|\leq c(1+|x|^2),
\exists c>0$, and that $f^{(\gamma+\Delta)}(x)\to f^{(\gamma)}$
almost everywhere by taking a subsequence, if necessary, since\vadjust{\goodbreak}
$w^{(\gamma+\Delta)}$ converges strongly in $H^1_{\mathrm{loc}}$ to
$w^{(\gamma)}$ by Lemma~\ref{lem6.2}. Moreover,
we note that $\{m_{\gamma+\Delta}(dx)\}=\{m_{\gamma+\Delta}(x)\,dx\}$
is tight because of Lemma~\ref{lem6.1}. Therefore, it converges weakly to some
probability measure ${\tilde m}(dx)$ by taking a subsequence if
necessary. The limit can be identified with $m_{\gamma}(dx)=m_{\gamma
}(x)\,dx$, and~$m_{\gamma}(x)$ is the only function satisfying (\ref{eq9.26})
for ${\bar L}(\gamma)$ and $\int m_{\gamma}(x)\,dx =1$. Thus $m_{\gamma
+\Delta}(x)\,dx$ converges to $m_{\gamma}(x)\,dx$ weakly. Therefore,
\[
\theta(\gamma+\Delta)=-\int f^{(\gamma+\Delta)}(x)m_{\gamma
+\Delta}(x)\,dx
\]
converges to $\theta(\gamma)$.

On the other hand, since $u^{(\gamma+\Delta)}$ is a solution to (\ref{eq6.2})
it satisfies
\[
\sup_{x\in{\mathcal D}^c}\frac{|u^{(\gamma+\Delta)}|}{{\bar
w}^{(\gamma+\Delta)}}<\infty,
\]
and we have $|{{\bar w}^{(\gamma+\Delta)}}|\leq c(1+|x|^2)$.
Therefore, we see that $u^{(\gamma+\Delta)}$ is locally bounded by a
constant independent of $\Delta$. Then, by testing (\ref{eq6.2}) with
$u^{(\gamma+\Delta)}\tau$, we can see that $\| u^{(\gamma
+\Delta)}\|_{H^1(B_r)}$ is bounded for each $r$. Therefore, by
taking a subsequence if necessary, $u^{(\gamma+\Delta)}$ converges in
$H^1_{\mathrm{loc}}$ weakly to a function ${\tilde u}$, which turns out to be
the solution $u^{(\gamma)}$ to (\ref{eq5.1}). By similar arguments as in the
proof of Lemma~\ref{lem6.2}, we can see it converges in $H^1_{\mathrm{loc}}$ strongly.
Furthermore, by Theorem~9.11 in \cite{GT}, we see that
\begin{eqnarray*}
\bigl\| u^{(\gamma+\Delta)}\bigr\|_{W^{2,p}(B_r)}&\leq
&C\bigl(\bigl\| u^{(\gamma+\Delta)}\bigr\|_{L^p(B_{2r})}+\bigl\|
f^{(\gamma+\Delta)}+\theta(\gamma+\Delta)\bigr\|
_{L^p(B_{2r})}\bigr) \\
&\leq&C\bigl(\bigl\| u^{(\gamma+\Delta)}\bigr\|_{L^{\infty
}(B_{2r})}+\bigl\| f^{(\gamma+\Delta)}+\theta(\gamma+\Delta
)\bigr\|_{L^p(B_{2r})}\bigr)
\end{eqnarray*}
for each $r>0$, where $C$ is a constant depending on $c_1, c_2$ and
the $L^{\infty}(B_{2r})$ norms of the coefficients of ${\bar L}(\gamma
+\Delta)$. Thus, by
the Sobolev imbedding theorem, $\{u^{(\gamma+\Delta)}\}$ is
equicontinuous, and thus $u^{(\gamma+\Delta)}$ converges uniformly
to~$u^{(\gamma)}$ on each compact set.
\end{pf}

%
\begin{lem}\label{lem6.4} Let $({\bar w}^{(\gamma)},\chi(\gamma))$ and $({\bar
w}^{(\gamma+\Delta)}, \chi(\gamma+\Delta))$ be solutions to~(\ref{eq3.1})
with $\gamma$ and $\gamma+\Delta$, respectively. Set $\chi^{(\Delta
)}=\frac{\chi(\gamma+\Delta)-\chi(\gamma)}{\Delta}$ and
$\zeta^{(\Delta)}=\frac{{\bar w}^{(\gamma+\Delta)}- {\bar
w}^{(\gamma)}}{\Delta}$. Then
\[
\lim_{|\Delta|\to0} \chi^{(\Delta)}= \theta(\gamma)
\]
and
\[
\lim_{|\Delta|\to0}\zeta^{(\Delta)}(x)=u^{(\gamma)}(x), \qquad
x\in R^n.
\]
Here, $(u^{(\gamma)}, \theta(\gamma))$ is the solution to (\ref{eq5.1}).
\end{lem}

\begin{pf} Here we abbreviate $u^{(\gamma)}$ to $u$. From (\ref{eq3.1}) it
follows that
\begin{eqnarray*}
&&-\chi(\gamma+\Delta)+\chi(\gamma)\\[-2pt]
&& \qquad =\tfrac{1}{2}\operatorname{tr}\bigl[\lambda
\lambda^*D^2\bigl({\bar w}^{(\gamma+\Delta)}-{\bar w}^{(\gamma)}\bigr)\bigr]+\beta
_{\gamma+\Delta}^*D{\bar w}^{(\gamma+\Delta)} \\[-2pt]
&& \qquad  \quad {} -\beta_{\gamma}^*D{\bar w}^{(\gamma)}-\tfrac
{1}{2}\bigl(D{\bar w}^{(\gamma+\Delta)}\bigr)^*\lambda N_{\gamma+\Delta
}^{-1}\lambda^*D{\bar w}^{(\gamma+\Delta)} \\[-2pt]
&& \qquad  \quad {} +\tfrac{1}{2}\bigl(D{\bar w}^{(\gamma)}\bigr)^*\lambda N_{\gamma
}^{-1}\lambda^*D{\bar w}^{(\gamma)}+U_{\gamma+\Delta}-U_{\gamma
} \\[-2pt]
&& \qquad =\tfrac{1}{2}\operatorname{tr}\bigl[\lambda\lambda^*D^2\bigl({\bar w}^{(\gamma
+\Delta)}-{\bar w}^{(\gamma)}\bigr)\bigr]+\beta_{\gamma}^*D\bigl({\bar w}^{(\gamma
+\Delta)}-{\bar w}^{(\gamma)}\bigr) \\[-2pt]
&& \qquad  \quad {} -\bigl(D{\bar w}^{(\gamma)}\bigr)^*\lambda N_{\gamma
}^{-1}\lambda^*D\bigl({\bar w}^{(\gamma+\Delta)}-{\bar w}^{(\gamma)}\bigr)
-\tfrac{1}{2}\bigl(D{\bar w}^{(\gamma)}\bigr)^*\lambda N_{\gamma}^{-1}\lambda
^*D{\bar w}^{(\gamma)} \\[-2pt]
&& \qquad  \quad {} +\bigl(D{\bar w}^{(\gamma)}\bigr)^*\lambda N_{\gamma
}^{-1}\lambda^*D{\bar w}^{(\gamma+\Delta)}-\tfrac{1}{2}\bigl(D{\bar
w}^{(\gamma+\Delta)}\bigr)^*\lambda N_{\gamma+\Delta}^{-1}\lambda
^*D{\bar w}^{(\gamma+\Delta)} \\[-2pt]
&& \qquad  \quad {}
+(\beta_{\gamma+\Delta}-\beta_{\gamma})^*D{\bar w}^{(\gamma
+\Delta)}+U_{\gamma+\Delta}-U_{\gamma} \\[-2pt]
&& \qquad ={\bar L}(\gamma)\bigl({\bar w}^{(\gamma+\Delta)}-{\bar w}^{(\gamma
)}\bigr)\\[-2pt]
&& \qquad  \quad {}-\tfrac{1}{2}D\bigl({\bar w}^{(\gamma+\Delta)}-{\bar w}^{(\gamma
)}\bigr)^*\lambda N_{\gamma+\Delta}^{-1}\lambda^* D\bigl({\bar w}^{(\gamma
+\Delta)}-{\bar w}^{(\gamma)}\bigr) \\[-2pt]
&& \qquad  \quad {} +(\beta_{\gamma+\Delta}-\beta_{\gamma})^*D{\bar
w}^{(\gamma+\Delta)}+\tfrac{1}{2}\bigl(D{\bar w}^{(\gamma)}\bigr)^*\lambda(
N_{\gamma+\Delta}^{-1}-N_{\gamma}^{-1})\lambda^*D{\bar w}^{(\gamma
)} \\[-2pt]
&& \qquad  \quad {} -\bigl(D{\bar w}^{(\gamma)}\bigr)^*\lambda( N_{\gamma+\Delta
}^{-1}-N_{\gamma}^{-1})\lambda^*D{\bar w}^{(\gamma+\Delta)}
+U_{\gamma+\Delta}-U_{\gamma}.
\end{eqnarray*}
Therefore we have
%
\begin{equation}\label{eq6.3}
-\chi^{(\Delta)}={\bar L}(\gamma)\zeta^{(\Delta)}+f_1^{(\Delta
)}(x)-g^{(\Delta)}(x),
\end{equation}
where
\begin{eqnarray*}
 f_1^{(\Delta)}(x)&=&\frac{(\beta_{\gamma+\Delta}-\beta_{\gamma
})^*}{\Delta}D{\bar w}^{(\gamma+\Delta)}+\frac{1}{2}\bigl(D{\bar
w}^{(\gamma)}\bigr)^*\lambda\frac{( N_{\gamma+\Delta}^{-1}-N_{\gamma
}^{-1})}{\Delta}\lambda^*D{\bar w}^{(\gamma)} \\
&&{} -\bigl(D{\bar w}^{(\gamma)}\bigr)^*\lambda\frac{( N_{\gamma
+\Delta}^{-1}-N_{\gamma}^{-1})}{\Delta}\lambda^*D{\bar w}^{(\gamma
+\Delta)}
+\frac{U_{\gamma+\Delta}-U_{\gamma}}{\Delta}
\end{eqnarray*}
and
\[
g^{(\Delta)}(x)=\frac{1}{2\Delta}D\bigl({\bar w}^{(\gamma+\Delta
)}-{\bar w}^{(\gamma)}\bigr)^*\lambda N_{\gamma+\Delta}^{-1}\lambda^*
D\bigl({\bar w}^{(\gamma+\Delta)}-{\bar w}^{(\gamma)}\bigr).
\]
Note that $f_1^{(\Delta)}$ is dominated by $c(1+|x|^2)$ with a certain
positive constant~$c$ and that it converges almost everywhere to
$f^{(\gamma)}$
by taking a subsequence, if necessary, because
\begin{eqnarray*}
\frac{\partial\beta_{\gamma}}{\partial\gamma}&=&\frac{1}{(1-\gamma
)^2}\lambda\Sigma^*{\hat\alpha},
\\
\frac{\partial N_{\gamma}^{-1}}{\partial\gamma}&=&\frac{1}{(1-\gamma
)^2}\Sigma^*(\Sigma\Sigma^*)^{-1}\Sigma,
\\
\frac{\partial U_{\gamma}}{\partial\gamma}&=&-\frac{1}{2(1-\gamma
)^2}\Sigma^*(\Sigma\Sigma^*)^{-1}\Sigma.
\end{eqnarray*}

Therefore,
%
\begin{equation}\label{eq6.4}
 \qquad -\chi_1^{(\Delta)}:=\int f_1^{(\Delta)}(x)m_{\gamma}(dx)\to\int
f(x)m_{\gamma}(dx)=-\theta(\gamma), \qquad   |\Delta|\to0.
\end{equation}
Moreover, for $\Delta>0$,
%
\begin{equation}\label{eq6.5}
-\chi_1^{(\Delta)}\geq-\chi^{(\Delta)}.
\end{equation}
Let us consider the equation
%
\begin{equation}\label{eq6.6}
-\chi_1^{(\Delta)}={\bar L}(\gamma)u_1^{(\Delta)}+f_1^{(\Delta)}.
\end{equation}

We shall see that $u_1^{(\Delta)}-u \to0$ as $\Delta\to0$ by
specifying suitable ambiguity constants. For that, set
\[
z^{(\Delta)}:= u_1^{(\Delta)}-u,  \qquad   F^{(\Delta)}:=f_1^{(\Delta
)}-f+\chi_1^{(\Delta)}-\theta(\gamma).
\]
Then we have
\[
{\bar L}(\gamma)z^{(\Delta)}+F^{(\Delta)}=0,\qquad z^{(\Delta)}\in
W^{2,p}_{\mathrm{loc}},  \qquad  \sup_{{\mathcal D}^c}\frac{|z^{(\Delta)}|}{{\bar
w}^{(\gamma)}}<\infty.
\]
By considering constructing the solution to this equation according to
the proof of Proposition~\ref{prop9.4} in the \hyperref[appm]{Appendix}, we see that $z^{(\Delta
)}\to0$ as $\Delta\to0$. To begin,
let $\Psi^{(\Delta)}$ be the solution to (\ref{eq9.21}) for $L_0={\bar
L}(\gamma)$, $f=F^{(\Delta)}$ and $\xi^{(\Delta)}$ the solution to
(\ref{eq9.22}) for $L_0={\bar L}(\gamma)$, $f=F^{(\Delta)}$ and $\Psi=\Psi
^{(\Delta)}$. The operator~$T$ is defined as
\[
TF^{(\Delta)}(x)=\xi^{(\Delta)}(x),  \qquad  x\in\Gamma_1,
\]
and the operator $P$ is defined in the same way as in (\ref{eq9.11}) by
replacing $L_0$ with ${\bar L}(\gamma)$ in (\ref{eq9.4}) and (\ref{eq9.10}). Then
starting with $\zeta^{(\Delta)}_0=\Psi^{(\Delta)}$, $\eta^{(\Delta
)}_0=\xi^{(\Delta)}$, we define the sequence $\zeta^{(\Delta)}_k,
\eta^{(\Delta)}_k$, $k=1,2,\ldots,$ successively as the solution to
(\ref{eq9.10}) with $\phi=\eta^{(\Delta)}_{k-1}$ and $L_0={\bar L}(\gamma
)$, and as the solution to (\ref{eq9.4}) with $h=\zeta^{(\Delta)}_k$ and
$L_0={\bar L}(\gamma)$, respectively. Then we obtain
\[
{\bar\eta}^{(\Delta)}(x)=\sum_{k=0}^{\infty}  \eta
_k^{(\Delta)} \big|_{\Gamma_1}=\sum_{k=0}^{\infty}P^k\bigl(TF^{(\Delta)}\bigr)(x)
\]
and the estimate for
${\bar\eta}^{(\Delta)}$,
\[
\bigl\|{\bar\eta}^{(\Delta)}\bigr\|_{L^{\infty}(\Gamma
_1)}\leq K\bigl\| TF^{(\Delta)}\bigr\|_{L^{\infty}(\Gamma
_1)}\frac{1}{1-e^{-\rho}}.
\]
To estimate $\| TF^{(\Delta)}\|_{L^{\infty}(\Gamma
_1)}$, we set
$\xi_1^{(\Delta)}$ to be the solution to
\[
\cases{\displaystyle
{\bar L}(\gamma)\xi_1^{(\Delta)}+F^{(\Delta)}=0, & \quad  $\overline
{{\mathcal D}}^c$,\vspace*{2pt} \cr\displaystyle
  \xi_1^{(\Delta)} \big|_{\Gamma}=0,
}
\]
and $\xi_2^{(\Delta)}=\xi^{(\Delta)}- \xi_1^{(\Delta)}$. Then
$\xi_2^{(\Delta)}$ satisfies
\[
\cases{\displaystyle
{\bar L}(\gamma)\xi_2^{(\Delta)}=0,  \qquad  \overline{{\mathcal
D}}^c,\vspace*{2pt} \cr\displaystyle
  \xi_2^{(\Delta)}\big |_{\Gamma}=  \Psi^{(\Delta
)} |_{\Gamma},
}
\]
and we have
\[
\bigl\|\xi_2^{(\Delta)}\bigr\|_{L^{\infty}(\Gamma_1)}\leq
\bigl\|\xi_2^{(\Delta)}\bigr\|_{L^{\infty}({\mathcal
D}^c)}\leq\bigl\|\Psi^{(\Delta)} \bigr\|_{L^{\infty}(\Gamma
)}\leq K_1\bigl\| F^{(\Delta)}\bigr\|_{L^{\infty}({\mathcal D}_1)}
\]
for some constant $K_1>0$. On the other hand, to estimate $ \|
\xi_1^{(\Delta)}\|_{L^{\infty}(\Gamma_1)}$, we set
\[
\xi_{1}^{(\Delta)}:=v^{(\Delta)}\bigl({\bar w}^{(\gamma)}\bigr)^{\alpha
},\qquad\alpha>1.
\]
We can assume that ${\mathcal D}$ is sufficiently large so that
%
\begin{eqnarray}
-U_{\gamma}-\frac{1}{2}\bigl(D{\bar w}^{(\gamma)}\bigr)^*\lambda N_{\gamma
}^{-1}\lambda^*D{\bar w}^{(\gamma)}-\chi(\gamma)+\frac{\alpha
-1}{{\bar w}^{(\gamma)}}\bigl(D{\bar w}^{(\gamma)}\bigr)^*\lambda\lambda
^*D{\bar w}^{(\gamma)}<-M,\nonumber\\
\eqntext{\displaystyle x\in{\mathcal D}^c,}
\end{eqnarray}
for some $M>0$. Since
\begin{eqnarray*}
  {\bar L}(\gamma)\xi_{1}^{(\Delta)}&=&\bigl({\bar w}^{(\gamma)}\bigr)^{\alpha
}{\bar L}(\gamma)v^{(\Delta)}+\alpha v^{(\Delta)}\bigl({\bar w}^{(\gamma
)}\bigr)^{\alpha-1}{\bar L}(\gamma){\bar w}^{(\gamma)} \\
&&{}+\alpha\bigl(Dv^{(\Delta)}\bigr)^*\lambda\lambda^*\frac{D{\bar w}^{(\gamma
)}}{{\bar w}^{(\gamma)}}\bigl({\bar w}^{(\gamma)}\bigr)^{\alpha}\\
&&{}+\alpha
(\alpha-1)v^{(\Delta)}\biggl(\frac{D{\bar w}^{(\gamma)}}{{\bar
w}^{(\gamma)}}\biggr)^*\lambda\lambda^*\frac{D{\bar w}^{(\gamma)}}{{\bar
w}^{(\gamma)}}\bigl({\bar w}^{(\gamma)}\bigr)^{\alpha},
\end{eqnarray*}
$v^{(\Delta)}$ satisfies
\[
\cases{\displaystyle
{\bar L}(\gamma)v^{(\Delta)}+\alpha\biggl(\frac{D{\bar w}^{(\gamma
)}}{{\bar w}^{(\gamma)}}\biggr)^*\lambda\lambda^*Dv^{(\Delta)} \cr\displaystyle
   \quad {}+\frac{\alpha}{{\bar w}^{(\gamma)}}\biggl\{{\bar L}(\gamma
){\bar w}^{(\gamma)}+\frac{\alpha-1}{{\bar w}^{(\gamma)}}\bigl(D{\bar
w}^{(\gamma)}\bigr)^*\lambda\lambda^*D{\bar w}^{(\gamma)}\biggr\}v^{(\Delta
)}  =-\frac{F^{(\Delta)}}{({\bar w}^{(\gamma)})^{\alpha}},
\vspace*{5pt}\cr\displaystyle
  v^{(\Delta)} \big|_{\partial{\mathcal D}}=0.
}
\]
Noting that
\begin{eqnarray*}
&&{\bar L}(\gamma){\bar w}^{(\gamma)}+\frac{\alpha-1}{{\bar
w}^{(\gamma)}}\bigl(D{\bar w}^{(\gamma)}\bigr)^*\lambda\lambda^*D{\bar
w}^{(\gamma)} \\
&& \qquad
=-U_{\gamma}-\frac{1}{2}\bigl(D{\bar w}^{(\gamma)}\bigr)^*\lambda N_{\gamma
}^{-1}\lambda^*D{\bar w}^{(\gamma)}-\chi(\gamma)+\frac{\alpha
-1}{{\bar w}^{(\gamma)}}\bigl(D{\bar w}^{(\gamma)}\bigr)^*\lambda\lambda
^*D{\bar w}^{(\gamma)} \\
&& \qquad <-M,  \qquad    x\in{\mathcal D}^c,
\end{eqnarray*}
we have
\[
\bigl|v^{(\Delta)}(x)\bigr|\leq K_2\sup_{{\mathcal D}^c}\frac{|F^{(\Delta
)}|}{({\bar w}^{(\gamma)})^{\alpha}},
\]
and thus
\[
\bigl\|\xi_1^{(\Delta)}\bigr\|_{L^{\infty}(\Gamma_1)}\leq
K_2\sup_{{\mathcal D}^c}\frac{|F^{(\Delta)}|}{({\bar w}^{(\gamma
)})^{\alpha}}\bigl\|\bigl({\bar w}^{(\gamma)}\bigr)^{\alpha}\bigr\|
_{L^{\infty}(\Gamma_1)}.
\]
Moreover, $\xi_1^{(\Delta)}=v^{(\Delta)}({\bar w}^{(\gamma
)})^{\alpha}\to0$ uniformly on each compact set as $\Delta\to0$.
Therefore, $\xi^{(\Delta)}\to0$ uniformly on each compact set and we
also obtain the estimates
\[
\bigl\|\xi^{(\Delta)}\bigr\|_{L^{\infty}(\Gamma_1)}\leq
K_1\bigl\| F^{(\Delta)}\bigr\|_{L^{\infty}({\mathcal D}_1)}+
K_2\sup_{{\mathcal D}^c}\frac{|F^{(\Delta)}|}{({\bar w}^{(\gamma
)})^{\alpha}}
\bigl\|\bigl({\bar w}^{(\gamma)}\bigr)^{\alpha}\bigr\|_{L^{\infty
}(\Gamma_1)}
\]
and
\[
\bigl\|{\bar\eta}^{(\Delta)}\bigr\|_{L^{\infty}(\Gamma
_1)}\leq K_1'\bigl\| F^{(\Delta)}\bigr\|_{L^{\infty}({\mathcal
D}_1)}+ K_2'\sup_{{\mathcal D}^c}\frac{|F^{(\Delta)}|}{({\bar
w}^{(\gamma)})^{\alpha}}
\bigl\|\bigl({\bar w}^{(\gamma)}\bigr)^{\alpha}\bigr\|_{L^{\infty
}(\Gamma_1)}.
\]

Let ${\tilde\zeta}^{(\Delta)}$ be the solution to (\ref{eq9.28}) for
$L_0={\bar L}(\gamma)$, $f=F^{(\Delta)}$ and ${\bar\eta}={\bar\eta
}^{(\Delta)}$.
Then ${\tilde\zeta}^{(\Delta)} \to0$ uniformly as $\Delta\to0$
since $\|{\bar\eta}^{(\Delta)}\|$ is estimated as
shown above.
Let ${\tilde\eta}^{(\Delta)}$ be the solution to (\ref{eq9.22}) for
$L_0={\bar L}(\gamma)$, $f=F^{(\Delta)}$ and $\Psi={\tilde\zeta
}^{(\Delta)}$. Then, as above,
$
{\tilde\eta}^{(\Delta)} \to0
$
uniformly on each compact set as $\Delta\to0$. Since $z^{(\Delta
)}={\tilde\zeta}^{(\Delta)}$ in ${\mathcal D}_1$ and
$z^{(\Delta)}={\tilde\eta}^{(\Delta)}$ in ${\mathcal D}^c$, we
conclude that
$
z^{(\Delta)} \to0
$
uniformly on each compact set.

In a similar manner, we have
\begin{eqnarray*}
-\chi^{(\Delta)}&=&{\bar L}(\gamma+\Delta)\zeta^{(\Delta)}+\frac
{(\beta_{\gamma+\Delta}-\beta_{\gamma})^*}{\Delta}Dw^{(\gamma
)} \\
&&{} +\frac{1}{2}\bigl(D{\bar w}^{(\gamma)}\bigr)^*\lambda\frac
{N_{\gamma}^{-1}-N_{\gamma+\Delta}^{-1}}{\Delta}\lambda^*D{\bar
w}^{(\gamma)}+\frac{U_{\gamma+\Delta}-U_{\gamma}}{\Delta} \\
&&{} +\frac{1}{2\Delta}D\bigl({\bar w}^{(\gamma+\Delta)}-{\bar
w}^{(\gamma)}\bigr)^*\lambda N_{\gamma+\Delta}^{-1}\lambda^*D\bigl({\bar
w}^{(\gamma+\Delta)}-{\bar w}^{(\gamma)}\bigr).
\end{eqnarray*}
By setting
\begin{eqnarray*}
f_2^{(\Delta)}(x)&:=&\frac{(\beta_{\gamma+\Delta}-\beta_{\gamma
})^*}{\Delta}Dw^{(\gamma)}+\frac{1}{2}\bigl(D{\bar w}^{(\gamma
)}\bigr)^*\lambda\frac{N_{\gamma}^{-1}-N_{\gamma+\Delta}^{-1}}{\Delta
}\lambda^*D{\bar w}^{(\gamma)}\\
&&{}+\frac{U_{\gamma+\Delta}-U_{\gamma
}}{\Delta},
\end{eqnarray*}
we have
%
\begin{equation}\label{eq6.7}
-\chi^{(\Delta)}={\bar L}(\gamma+\Delta)\zeta^{(\Delta
)}+f_2^{(\Delta)}(x)+g^{(\Delta)}(x).
\end{equation}
Since $m_{\gamma+\Delta}(x)\,dx$ converges to $m_{\gamma}(x)\,dx$
weakly, and $f_2^{(\Delta)}$ converges almost everywhere to
$f(x)$ by taking a subsequence if necessary, as above, we have
%
\begin{eqnarray}\label{eq6.8}
-\chi_2^{(\Delta)}:=\int f_2^{(\Delta)}(x)m_{\gamma+\Delta}(x)\,dx
\to\int f(x)m_{\gamma}(x)\,dx=-\theta(\gamma)   \nonumber
\\[-8pt]
\\[-8pt] \eqntext{\mbox{as }  |\Delta
|\to0.}
\end{eqnarray}
Moreover, for $\Delta>0$,
%
\begin{equation}\label{eq6.9}
-\chi_2^{(\Delta)}\leq-\chi^{(\Delta)}.
\end{equation}
We consider
%
\begin{equation}\label{eq6.10}
-\chi_2^{(\Delta)}={\bar L}(\gamma+\Delta)u_2^{(\Delta
)}+f_2^{(\Delta)}.
\end{equation}
Then, in the same manner as above, we see that $u_2^{(\Delta
)}-u^{(\gamma+\Delta)}\to0$, as $|\Delta|\to0$, by specifying
ambiguity constants. Since $u^{(\gamma+\Delta)}$ converges to $u$,
$u_2^{(\Delta)}$ does the same.

From (\ref{eq6.4}), (\ref{eq6.5}), (\ref{eq6.8}) and (\ref{eq6.9}), it follows that
%
\begin{equation}\label{eq6.11}
\lim_{\Delta\downarrow0}-\chi^{(\Delta)}= -\theta(\gamma).
\end{equation}
The converse inequalities of (\ref{eq6.5}) and (\ref{eq6.9}) hold for $\Delta<0$, and
we have
%
\renewcommand{\theequation}{6.11$'$}
\begin{equation}\label{eq6.11'}
\lim_{\Delta\uparrow0}-\chi^{(\Delta)}= -\theta(\gamma).
\renewcommand{\theequation}{\arabic{section}.\arabic{equation}}
\setcounter{equation}{11}
\end{equation}
\noindent
Hence, we see that $-\chi^{(\Delta)}\to-\theta(\gamma)$ as
$|\Delta|\to0$.
From (\ref{eq6.3}) and (\ref{eq6.6}), we have
\[
-\chi_1^{(\Delta)}+\chi^{(\Delta)}={\bar L}(\gamma)\bigl(u_1^{(\Delta
)}-\zeta^{(\Delta)}\bigr)+g^{(\Delta)},
\]
and through arguments similar to those above, we see that
%
\begin{equation}\label{eq6.12}
\liminf_{\Delta\downarrow0} \bigl(u_1^{(\Delta)}(x)-\zeta^{(\Delta
)}(x)\bigr)\geq0,
\end{equation}
since $g^{(\Delta)}(x)\geq0$ for $\Delta>0$ and $\chi^{(\Delta
)}-\chi_1^{(\Delta)}\to0$ as $|\Delta|\to0$. Similarly, from~(\ref{eq6.7})
and (\ref{eq6.10}), we have
\[
-\chi_2^{(\Delta)}+\chi^{(\Delta)}={\bar L}(\gamma+\Delta
)\bigl(u_2^{(\Delta)}-\zeta^{(\Delta)}\bigr)-g^{(\Delta)}
\]
and see that
%
\begin{equation}\label{eq6.13}
\limsup_{\Delta\downarrow0}\bigl(u_2^{(\Delta)}(x)-\zeta^{(\Delta
)}(x)\bigr)\leq0.
\end{equation}
Therefore,
\[
\lim_{\Delta\downarrow0}\zeta^{(\Delta)}(x)=u(x).
\]
We likewise have
\[
\limsup_{\Delta\uparrow0} \bigl(u_1^{(\Delta)}(x)-\zeta^{(\Delta
)}(x)\bigr)\leq0
\]
and
\[
\liminf_{\Delta\uparrow0}\bigl(u_2^{(\Delta)}(x)-\zeta^{(\Delta
)}(x)\bigr)\geq0
\]
since $g^{(\Delta)}(x)\leq0$ for $\Delta<0$.
Therefore, we obtain
\[
\lim_{\Delta\uparrow0}\zeta^{(\Delta)}(x)=u(x)
\]
and Lemma~\ref{lem6.4} follows.
\end{pf}

%
\begin{rem}\label{rem6.1}
Since $u=u^{(\gamma)}=\frac{\partial{\bar w}}{\partial\gamma}$ is
a solution to (\ref{eq5.1}), it has a polynomial growth order. More precisely,
we have
\[
|u(x)|\leq C(1+|x|^2),\qquad\exists C>0;
\]
cf. Corollary~\ref{cor5.1} and (\ref{eq3.11}). Furthermore, we can see that $\frac
{\partial u}{\partial x_l}$\vspace*{-1pt} also has a~polynomial growth order for each
$l$. Indeed, $u_l:= \frac{\partial u}{\partial x_l}$ satisfies
\[
0=\tfrac{1}{2}D_i(a^{ij}D_ju_l)+B^iD_iu_l-D_lf+\tfrac
{1}{2}D_i(a^{ij}_lD_ju)+B_l^iD_iu,
\]
where
\begin{eqnarray*}
 a^{ij}(x)&=&(\lambda\lambda^*)^{ij}(x), \\
     B^j(x)&=&\beta_{\gamma
}^j(x)-\tfrac{1}{2}D_i((\lambda\lambda)^{ij})-(\lambda N_{\gamma
}^{-1}\lambda^*D{\bar w})^j, \\
f&=&\frac{1}{2(1-\gamma)^2}(\lambda^*D{\bar w}-\Sigma^*{\hat\alpha
})^*\Sigma^*(\Sigma\Sigma^*)^{-1}\Sigma(\lambda^*D{\bar w}-\Sigma
^*{\hat\alpha}), \\
a_l^{ij}&=&\frac{\partial a^{ij}}{\partial x_l}, \qquad    B_l^{i}=\frac
{\partial B^{i}}{\partial x_l}.
\end{eqnarray*}
Therefore, if we set $f_i=-\frac{1}{2}a^{ij}_lD_ju, i\neq l$, and
$f_l=-\frac{1}{2}a^{lj}_lD_ju-f$, then $u_l$ satisfies
%
\begin{equation}\label{eq6.14}
\int\biggl(\frac{1}{2}a^{ij}D_ju_l-f_i\biggr)\xi_{x_i}\,dx+\int
(B^iD_iu_l+B_l^iD_iu)\xi \,dx =0
\end{equation}
for each $\xi\in W^{1,2}_0(B_{\rho}(x_0))$, $\rho>0,   x_0\in
R^n$. We note that
%
\begin{eqnarray}
 \| f_i \|_{L^{p/2}(B_{\rho}(x_0))},    \|
B^i\|_{L^{p/2}(B_{\rho}(x_0))},   \| B^i_lD_iu\|
_{L^{p/2}(B_{\rho}(x_0))}\!\leq\!\mu(x_0)   \!\leq\! C(1+|x_0|^{m_0}),  \nonumber\\
\eqntext{\displaystyle \exists C\,{>}\,0,  m_0\,{>}\,0,\hspace*{3pt}}
\end{eqnarray}
which can be seen in light of our assumptions since $u$ is a solution
to (\ref{eq5.1}), and ${\bar w}$ is a solution to (\ref{eq3.1}). Equation (\ref{eq6.14})
corresponds to (13.4) in \cite{Lad},  Chapter 3, Section 13.
Therefore, by the same arguments as in that work,
\[
\int_{A_k,\rho}|\nabla u_l|^2\zeta^2\,dx \leq\gamma(x_0)\biggl[\int
_{A_{k,\rho}}(u_l-k)^2|\nabla\zeta|^2\,dx+(k^2+1)|A_{k,\rho}|^{1-2/p}\biggr]
\]
is seen to hold, where $\zeta$ is a cut-off function supported by
$B_{\rho}(x_0)$, $A_{k,\rho}=\{x\in B_{\rho}(x_0) ; u_l(x) > k\}$,
and $\gamma(x_0)$ is a constant dominated by
$C(1+|x_0|^{m_1})$,  $C>0$, $m_1>0$. From this inequality we obtain
inequality ({5.12}) for $u=u_l$ in \cite{Lad}, Chapter 2, Section 5.
Hence, similarly to the proof of Lemma 5.4 in~\cite{Lad}, Chapter~2,
we see that $u_l$ has a polynomial growth order.
\end{rem}

\section{\texorpdfstring{Proof of Theorem~\protect\ref{theo2.4}}{Proof of Theorem 2.4}}\label{sec7}
We first state the following lemma.
%
\begin{lem}\label{lem7.1} Under the assumptions of Theorem~\ref{theo2.3},
%
\begin{equation}\label{eq7.1}
\chi'(-\infty):=\lim_{\gamma\to-\infty}\chi'(\gamma)=0.
\end{equation}
\end{lem}

\begin{pf} Note that
\[
0\leq\lim_{\gamma\to-\infty}\chi'(\gamma)\leq\chi'(\gamma_0)
\]
for $\gamma_0<0$ and that $\chi'(\gamma)$ is nondecreasing.
Therefore, $\chi'(-\infty)$ exists. Furthermore,
\[
\frac{\chi(\gamma)}{\gamma}=-\frac{1}{\gamma}\int_{\gamma
}^{\gamma_0}\chi'(t)\,dt+\frac{\chi(\gamma_0)}{\gamma}
\]
and $\lim_{\gamma\to-\infty}\frac{\chi(\gamma)}{\gamma}=0$
since $-\chi_0\leq\chi(\gamma)\leq0$. Here, $\chi_0$ is a
constant defined by (\ref{eq3.10}). Hence, we obtain (\ref{eq7.1}).
\end{pf}

We next give the proof of Theorem~\ref{theo2.4}. For $\gamma<0$, we have
\begin{eqnarray*}
\log P\biggl(\frac{\log V_T(h)-\log S_T^0}{T}\leq\kappa\biggr)&=&\log P\biggl(\biggl(\frac
{V_T(h)}{S_T^0}\biggr)^{\gamma}\geq e^{\gamma\kappa T}\biggr) \\
&\leq&\log\biggl\{ E\biggl[\biggl(\frac{V_T(h)}{S_T^0}\biggr)^{\gamma}\biggr]e^{-\gamma\kappa
T}\biggr\} \\
&=&\log E\biggl[\biggl(\frac{V_T(h)}{S_T^0}\biggr)^{\gamma}\biggr] -\gamma\kappa T.
\end{eqnarray*}
Therefore,
\begin{eqnarray*}
\inf_{h}\log P\biggl(\frac{\log V_T(h)-\log S_T^0}{T}\leq\kappa\biggr)&\leq&
\inf_{h}\log E\biggl[\biggl(\frac{V_T(h)}{S_T^0}\biggr)^{\gamma}\biggr] -\gamma\kappa
T\\
&\leq& v(0,x;T)-\gamma\kappa T,
\end{eqnarray*}
from which we obtain
\[
\liminf_{T\to\infty}\frac{1}{T}\inf_{h}\log P\biggl(\frac{\log
V_T(h)-\log S_t^0}{T}\leq\kappa\biggr)\leq\chi(\gamma)-\gamma\kappa
\]
for all $\gamma<0$. Hence, we have
\[
\liminf_{T\to\infty}\frac{1}{T}\inf_{h}\log P\biggl(\frac{\log
V_T(h)-\log S_T^0}{T}\leq\kappa\biggr)\leq\inf_{\gamma<0}\{ \chi(\gamma
)-\gamma\kappa\}.
\]

 The converse inequality is more difficult to prove. Take a
constant $\kappa$ and $\epsilon>0$ such that $\kappa-\epsilon>0$.
Then there exists $ \gamma_{\epsilon}$ such that
%
\begin{equation}\label{eq7.2}
\inf_{\gamma<0}\{\chi(\gamma)-\gamma(\kappa-\epsilon)\}=\chi
(\gamma_{\epsilon})-\gamma_{\epsilon}\chi'(\gamma_{\epsilon
}).
\end{equation}
We write $ \gamma_{\epsilon}$ as $\gamma$ for simplicity in the
following. Let us introduce\vspace*{1pt} a probability measure ${\tilde P}$ defined by
\[
  \frac{d {\tilde P}}{d P} \bigg|_{{\mathcal G}_T}=e^{
M^{\gamma}_T-\fracd{1}{2}\langle M^{\gamma}\rangle_T},
\]
where
\[
M^{\gamma}_t=\int_0^t\biggl\{ \frac{\gamma}{1-\gamma}{\hat\alpha
}^*\Sigma+(Dw)^*\lambda N_{\gamma}^{-1}\biggr\}(X_s)\,dW_s.
\]
Then ${\tilde W}_t=W_t-\int_0^t\{\frac{\gamma}{1-\gamma}\Sigma
^*{\hat\alpha}+N_{\gamma}^{-1}\lambda^*Dw\}(X_s)\,ds$ is a martingale
under the probability measure ${\tilde P}$ and
%
\begin{eqnarray}\label{eq7.3}
dX_t&=&\beta(X_t)\,dt+\lambda(X_t)\,dW_t \nonumber
\\[-8pt]
\\[-8pt]
&=&\biggl\{\beta+\frac{\gamma}{1-\gamma}\lambda\Sigma^*{\hat\alpha}+
\lambda N^{-1}_{\gamma}\lambda^*Dw\biggr\}(X_t)\,dt+\lambda(X_t)\,d{\tilde W}_t.
\nonumber
\end{eqnarray}
Note that $(X_t,{\tilde P})$ has the same law as the diffusion process
$({\bar Y}_t,P)$ governed by stochastic differential equation (\ref{eq4.12}).
We further note H-J-B equation of ergodic type (\ref{eq2.19}) can be written as
%
\begin{equation}\label{eq7.4}
\chi(\gamma)={\bar L}w-\tfrac{1}{2}(Dw)^*\lambda N_{\gamma
}^{-1}\lambda^*Dw-U_{\gamma}
\end{equation}
by using ${\bar L}$ defined by (\ref{eq4.16}).
On the other hand, (\ref{eq5.1}) is written as
%
\begin{eqnarray}\label{eq7.5}
\chi'(\gamma)&=&{\bar L}w_{\gamma}+\frac{1}{2(1-\gamma)^2}(\sigma
\lambda^*Dw+{\hat\alpha})^*(\sigma\sigma^*)^{-1}(\sigma\lambda
^*Dw+{\hat\alpha}) \nonumber\\
&=&{\bar L}w_{\gamma}+\frac{1}{2(1-\gamma)^2}(\lambda^*D w+\Sigma
^*{\hat\alpha})^*\Sigma^*(\Sigma\Sigma^*)^{-1}\Sigma(\lambda^*D
w+\Sigma^*{\hat\alpha})
 \\
&\hspace*{2.3pt}=:&{\bar L}w_{\gamma}+V_1(x),\nonumber
\end{eqnarray}
owing to Lemma~\ref{lem6.3}, where $w_{\gamma}=\frac{\partial w}{\partial
\gamma}$. Now we have
\begin{eqnarray*}
 &&\log V_T(h)-\log S_T^0\\
 && \qquad =\log v+\int_0^T\biggl\{h_s^*{\hat\alpha
}(X_s)-\frac{1}{2}h_s^*\sigma\sigma^*(X_s)h_s\biggr\}\,ds+\int
_0^Th_s^*\sigma(X_s)\,dW_s
 \\
&& \qquad =\log v+\int_0^Th_s^*\sigma(X_s)\,d{\tilde W}_s \\
&& \qquad  \quad {}+\int_0^T\biggl\{h_s^*{\hat\alpha}(X_s)-\frac{1}{2}h_s^*\sigma
\sigma^*(X_s)h_s+\frac{\gamma}{1-\gamma}h_s^*\sigma\Sigma^*{\hat
\alpha}(X_s)\\
&& \qquad \quad \hspace*{120pt}\hphantom{{}+\int_0^T\biggl\{} {}+h_s^*\sigma N_{\gamma}^{-1}\lambda^*Dw(X_s)\biggr\}\,ds \\
&& \qquad =\log v+\int_0^Th_s^*\sigma(X_s)\,d{\tilde W}_s \\
&& \qquad  \quad {}+\int_0^T\biggl\{
\frac{1}{1-\gamma}\bigl(h_s^*{\hat\alpha}(X_s)+h_s^*\sigma\lambda
^*Dw(X_s)\bigr)-\frac{1}{2}h_s^*\sigma\sigma^*(X_s)h_s\biggr\}\,ds \\
&& \qquad =\log v+\int_0^T\biggl\{h_s-\frac{1}{1-\gamma}(\sigma\sigma
^*)^{-1}({\hat\alpha}+\sigma\lambda^*Dw)\biggr\}^*\sigma(X_s)\,d{\tilde
W}_s \\
&& \qquad  \quad {}+\int_0^T\frac{1}{1-\gamma}\biggl\{(\sigma\sigma
^*)^{-1}({\hat\alpha}+\sigma\lambda^*Dw)\biggr\}^*\sigma(X_s)\,d{\tilde
W}_s \\
&& \qquad  \quad {}-\frac{1}{2}\int_0^T\biggl\{h_s-\frac{1}{1-\gamma}(\sigma
\sigma^*)^{-1}({\hat\alpha}+\sigma\lambda^*Dw)\biggr\}^*\\
&& \hphantom{-\frac{1}{2}\int_0^T\,}\qquad  \quad {}\times\sigma\sigma
^*\biggl\{h_s-\frac{1}{1-\gamma}(\sigma\sigma^*)^{-1}({\hat\alpha
}+\sigma\lambda^*Dw)\biggr\}\,ds \\
&& \qquad  \quad {}+\frac{1}{2(1-\gamma)^2}\int_0^T({\hat\alpha}+\sigma
\lambda^*Dw)^*(\sigma\sigma^*)^{-1}({\hat\alpha}+\sigma\lambda
^*Dw)(X_s)\,ds \\
&& \qquad =\log v+M_T^h-\frac{1}{2}\langle M^h\rangle_T \\
&& \qquad  \quad {}+\int_0^T\frac
{1}{1-\gamma}\{(\sigma\sigma^*)^{-1}({\hat\alpha}+\sigma\lambda
^*Dw)\}^*\sigma(X_s)\,d{\tilde W}_s \\
 &&\qquad  \quad {} +\int_0^TV_1(X_s)\,ds,
\end{eqnarray*}
and we set
\[
M_t^h:=\int_0^t\biggl\{h_s-\frac{1}{1-\gamma}(\sigma\sigma^*)^{-1}({\hat
\alpha}+\sigma\lambda^*Dw)\biggr\}^*\sigma(X_s)\,d{\tilde W}_s.
\]
Note that $\kappa-\epsilon=\chi'(\gamma)$. Then it follows that
\begin{eqnarray*}
&&{\tilde P}\biggl(\frac{1}{T}\bigl( \log V_T(h)-\log S_T^0\bigr) >\kappa\biggr)  \\
&& \qquad  \leq{\tilde P}\biggl(\frac{1}{T}\log v+\frac{1}{T}\int
_0^TV_1(X_s)\,ds>\chi'(\gamma)+\frac{\epsilon}{3}\biggr) \\
&& \qquad  \quad {} +{\tilde P}\biggl(\frac{1}{T}\biggl(M^h_T-\frac{1}{2}\langle
M^h\rangle_T\biggr)>\frac{\epsilon}{3}\biggr) \\
&&  \qquad  \quad {}+{\tilde P}\biggl(\frac{1}{T}\int_0^T\frac{1}{1-\gamma}\{
(\sigma\sigma^*)^{-1}({\hat\alpha}+\sigma\lambda^*Dw)\}^*\sigma
(X_s)\,d{\tilde W}_s>\frac{\epsilon}{3}\biggr).
\end{eqnarray*}

Taking Lemma~\ref{lem4.2} into account, we have
\begin{eqnarray*}
&&{\tilde P}\biggl(\frac{1}{T}\int_0^T\frac{1}{1-\gamma}\{(\sigma\sigma
^*)^{-1}({\hat\alpha}+\sigma\lambda^*Dw)\}^*\sigma(X_s)\,d{\tilde
W}_s>\frac{\epsilon}{3}\biggr)
 \\[-2pt]
&& \qquad  \leq\frac{9}{\epsilon^2 T^2}{\tilde E}\biggl[\int
_0^TV_1(X_s)\,ds\biggr] \\[-2pt]
&& \qquad  \leq\frac{C}{\epsilon^2 T}
\end{eqnarray*}
for some positive constant $C$ and
\[
{\tilde P}\biggl(\frac{1}{T}\biggl(M^h_T-\frac{1}{2}\langle M^h\rangle_T\biggr)>\frac
{\epsilon}{3}\biggr)\leq e^{-\fraca{\epsilon T}{3}}{\tilde E}\bigl[e^{M^h_T-\fracd
{1}{2}\langle M^h\rangle_T}\bigr]\leq e^{-\fraca{\epsilon T}{3}}.
\]
Thus, by using the following lemma, we can see that
%
\begin{equation}\label{eq7.6}
{\tilde P}\biggl(\frac{1}{T}\bigl( \log V_T(h)-\log S_T^0\bigr) >\kappa\biggr)<\epsilon
\end{equation}
for sufficiently large $T$.
%
\begin{lem}\label{lem7.2} For sufficiently large $T$ we have
\[
{\tilde P}\biggl(\frac{\log v}{T}+\frac{1}{T}\int_0^TV_1(X_s)\,ds>\chi
'(\gamma)+\frac{\epsilon}{3}\biggr)\leq\frac{\epsilon}{2}.
\]
\end{lem}

\begin{pf} By It\^o's formula,
\begin{eqnarray*}
&&
w_{\gamma}(X_T)-w_{\gamma}(X_0) \\[-2pt]
&& \qquad  =\int_0^T\{{\bar L}(\gamma)w_{\gamma}(X_s)\}\,ds+\int
_0^T(\nabla w_{\gamma}(X_s))^*\lambda(X_s)\,d{\tilde W}_s \\[-2pt]
&& \qquad  =-\int_0^TV_1(X_s)\,ds+\chi'(\gamma)T+\int_0^T(\nabla
w_{\gamma}(X_s))^*\lambda(X_s)\,d{\tilde W}_s.
\end{eqnarray*}
Therefore,
\begin{eqnarray*}
\frac{1}{T}\int_0^TV_1(X_s)\,ds&=&\chi'(\gamma)+\frac{1}{T}\{w_{\gamma
}(x)-w_{\gamma}(X_T)\}\\[-2pt]
&&{}+\frac{1}{T}\int_0^T(\nabla w_{\gamma
}(X_s))^*\lambda(X_s)\,d{\tilde W}_s.
\end{eqnarray*}
Thus,
\begin{eqnarray*}
&&{\tilde P}\biggl(\frac{1}{T}\log v+\frac{1}{T}\int_0^TV_1(X_s)\,ds>\chi
'(\gamma)+\frac{\epsilon}{3}\biggr) \\[-2pt]
&& \qquad ={\tilde P}\biggl(\frac{1}{T}\log v+\frac{1}{T}\{w_{\gamma
}(x)-w_{\gamma}(X_T)\}+
\frac{1}{T}\int_0^T(\nabla w_{\gamma}(X_s))^*\,d{\tilde W}_s>\frac
{\epsilon}{3}\biggr) \\[-2pt]
&& \qquad \leq{\tilde P}\biggl(\frac{1}{T}\log v>\frac{\epsilon
}{9}\biggr)+{\tilde P}\biggl(\frac{1}{T}\{w_{\gamma}(x)-w_{\gamma}(X_T)\}>\frac
{\epsilon}{9}\biggr)
\\[-2pt]
&& \qquad  \quad {}+{\tilde P}\biggl(\frac{1}{T}\int_0^T(\nabla w_{\gamma}(X_s))^*\,d{\tilde
W}_s>\frac{\epsilon}{9}\biggr) \\[-2pt]
&& \qquad \leq\frac{81}{\epsilon^2T^2}E[|w_{\gamma}(x)-w_{\gamma}(X_T)|^2]
+\frac{81}{\epsilon^2T^2}E\biggl[\int_0^T(D w_{\gamma})^*\lambda\lambda
^*Dw_{\gamma}(X_s)\,ds\biggr].
\end{eqnarray*}
Hence, by taking $T$ and $R$ to be sufficiently large, we obtain our
present lemma because of Lemma~\ref{lem4.2}; cf. Remark~\ref{rem6.1}.
\end{pf}

Let us complete the proof of Theorem~\ref{theo2.4} for $0<\kappa<\chi'(0-)$. Set
\[
{\tilde M}^{\gamma}_t=\int_0^t\biggl\{ \frac{\gamma}{1-\gamma}{\hat
\alpha}^*\Sigma+(Dw)^*\lambda N_{\gamma}^{-1}\biggr\}(X_s)\,d{\tilde W}_s
\]
and
\begin{eqnarray*}
 A_1&=&\{  -{\tilde M}^{\gamma}_T\geq-\epsilon T\}, \\[-2pt]
 A_2&=&\bigl\{ -\tfrac{1}{2}\langle M^{\gamma}\rangle_T\geq\bigl(\chi(\gamma
)-\gamma\chi'(\gamma)-\epsilon\bigr)T\bigr\}, \\[-2pt]
 A_3&=&\biggl\{ \frac{1}{T}\bigl(\log V_T(h)-\log S_T^0\bigr)\leq\kappa\biggr\}.
\end{eqnarray*}
Then
\begin{eqnarray*}
&&P\biggl(\frac{1}{T}\bigl(\log V_T(h)-S_T^0\bigr)\leq\kappa\biggr)\\[-2pt]
&& \qquad ={\tilde E}\biggl[e^{-{\tilde
M}^{\gamma}_T-\fracd{1}{2}\langle M^{\gamma}\rangle_T};\frac
{1}{T}\bigl(\log V_T(h)-S_T^0\bigr)\leq\kappa\biggr] \\[-2pt]
&& \qquad \geq{\tilde E}\bigl[e^{-{\tilde M}^{\gamma}_T-\fracd{1}{2}\langle
M^{\gamma}\rangle_T}; A_1\cap A_2\cap A_3\bigr] \\[-2pt]
&& \qquad \geq e^{(\chi(\gamma)-\gamma\chi'(\gamma)-2\epsilon)T}{\tilde
P}(A_1\cap A_2\cap A_3) \\[-2pt]
&& \qquad \geq e^{(\chi(\gamma)-\gamma\chi'(\gamma)-2\epsilon)T}\{
1-{\tilde P}(A_1^c)-{\tilde P}(A_2^c)-{\tilde P}(A_3^c)\}.
\end{eqnarray*}
We have seen that ${\tilde P}(A_3^c)<\epsilon$ holds for sufficiently
large $T$ in (\ref{eq7.6}),\vspace*{-1pt} and it is straightforward to see that likewise
${\tilde P}(A_1^c)<\epsilon$ for sufficiently large $T$. Therefore,
taking the following lemma into account as well, we have
\[
P\biggl(\frac{1}{T}\bigl(\log V_T(h)-\log S_T^0\bigr)\leq\kappa\biggr)\geq e^{(\chi
(\gamma)-\gamma\chi'(\gamma)-2\epsilon)T}(1-3\epsilon)\qquad
\forall h\in{\mathcal H}(T),
\]
for sufficiently large $T$, and
\begin{eqnarray*}
\varliminf_{T\to\infty}\frac{1}{T}\inf_{h\in{\mathcal H}(T)}\log
P\biggl(\frac{\log V_T(h)-\log S_T^0}{T}\leq\kappa\biggr)&\geq& \chi(\gamma
)-\gamma\chi'(\gamma)-2\epsilon \\[-2pt]
&=&\chi(\gamma)-\gamma(\kappa-\epsilon)-2\epsilon \\[-2pt]
&\geq& \inf_{\gamma<0}\{\chi(\gamma)-\gamma(\kappa-\epsilon)\}
-2\epsilon
\end{eqnarray*}
for each $\epsilon$. Since $\chi(\gamma)$ is smooth and convex,
$J(\kappa)=\inf_{\gamma<0}\{\chi(\gamma)-\gamma\kappa\}$,
$\kappa>0$, is strictly concave, and thus continuous. Hence,
\[
\hspace*{2pt}\qquad\varliminf_{T\to\infty}\frac{1}{T}\inf_{h\in{\mathcal H}(T)}\log
P\biggl(\frac{\log V_T(h)-\log S_T^0}{T}\leq\kappa\biggr)\geq\inf_{\gamma<0}\{
\chi(\gamma)-\gamma\kappa\}.\qquad\hspace*{2pt}\qed
\]

\begin{lem} Under the assumptions of Theorem~\ref{theo2.4},
\[
{\tilde P}\bigl(\tfrac{1}{2}\langle M^{\gamma}\rangle_T\geq-\bigl(\chi(\gamma
)-\gamma\chi'(\gamma)-\epsilon\bigr)T\bigr)<\epsilon
\]
holds for sufficiently large $T$.
\end{lem}

\begin{pf} First note that
\[
\frac{1}{2}N_{\gamma}^{-1}-\frac{1}{2}(N_{\gamma}^{-1})^2=-\frac
{\gamma}{2(1-\gamma)^2}\Sigma^*(\Sigma\Sigma^*)^{-1}\Sigma
\]
and
\[
\frac{\gamma}{2(1-\gamma)}-\frac{\gamma}{2(1-\gamma)^2}=-\frac
{\gamma^2}{2(1-\gamma)^2}.
\]
Then, from (\ref{eq7.4}) and (\ref{eq7.5}), we have
\begin{eqnarray*}
 \chi(\gamma)-\gamma\chi'(\gamma)&=&{\bar L}(w-\gamma w_{\gamma
})-\frac{1}{2}(Dw)^*\lambda N_{\gamma}^{-1}\lambda^*Dw-U_{\gamma
}-\gamma V_1(x) \\
&=&{\bar L}(w-\gamma w_{\gamma})-\frac{1}{2}(Dw)^*\lambda N_{\gamma
}^{-1}\lambda^*Dw-\frac{\gamma}{(1-\gamma)^2}(\lambda\Sigma
^*{\hat\alpha})^*Dw \\
&&{}+\frac{1}{2}(Dw)^*\lambda\bigl(N_{\gamma}^{-1}-(N_{\gamma
}^{-1})^2\bigr)\lambda^*Dw-{\gamma^2}{2(1-\gamma)^2}{\hat\alpha}\Sigma
\Sigma^*{\hat\alpha} \\
&
=&{\bar L}(w-\gamma w_{\gamma})\\
&&{}-\frac{1}{2}\biggl(\frac{\gamma}{1-\gamma
}\Sigma^*{\hat\alpha}+N_{\gamma}^{-1}\lambda^*Dw\biggr)^*\biggl(\frac{\gamma
}{1-\gamma}\Sigma^*{\hat\alpha}+N_{\gamma}^{-1}\lambda^*Dw\biggr).
\end{eqnarray*}
Set
\[
V_2(x)=\frac{1}{2}\biggl(\frac{\gamma}{1-\gamma}\Sigma^*{\hat\alpha
}+N_{\gamma}^{-1}\lambda^*Dw\biggr)^*\biggl(\frac{\gamma}{1-\gamma}\Sigma
^*{\hat\alpha}+N_{\gamma}^{-1}\lambda^*Dw\biggr).
\]
Then
\[
\frac{1}{2}\langle M^{\gamma}\rangle_t=\int_0^tV_2(X_s)\,ds.
\]
By It\^o's formula, we have
\begin{eqnarray*}
&&(w-\gamma w_{\gamma})(X_t)-(w-\gamma w_{\gamma})(X_0) \\
&& \qquad  =\int_0^t{\bar L}(\gamma)(w-\gamma w_{\gamma})(X_s)\,ds+\int
_0^tD(w-\gamma w_{\gamma})(X_s)^*\lambda(X_s)\,d{\tilde W}_s \\
&& \qquad  =\{\chi(\gamma)-\gamma\chi'(\gamma)\}t+\int_0^tV_2(X_s)\,ds+\int
_0^tD(w-\gamma w_{\gamma})(X_s)^*\lambda(X_s)\,d{\tilde W}_s.
\end{eqnarray*}
Thus,
\begin{eqnarray*}
&&{\tilde P}\biggl(\frac{1}{2}\langle M^{\gamma}\rangle_T+\{\chi(\gamma
)-\gamma\chi'(\gamma)\}T>\epsilon T\biggr) \\
&&\qquad=
{\tilde P}\biggl((w-\gamma w_{\gamma})(X_T)-(w-\gamma w_{\gamma})(x)\\
&&\qquad\hphantom{={\tilde P}\biggl(}    {}-\int
_0^tD(w-\gamma w_{\gamma})(X_s)^*\lambda(X_s)\,d{\tilde W}_s>\epsilon
T\biggr) \\
&& \qquad
\leq
{\tilde P}\biggl((w-\gamma w_{\gamma})(X_T)>\frac{\epsilon T}{3}\biggr)
+{\tilde P}\biggl(-(w-\gamma w_{\gamma})(x)>\frac{\epsilon T}{3}\biggr) \\
&& \qquad  \quad {}+{\tilde P}\biggl(-\int_0^tD(w-\gamma w_{\gamma})(X_s)^*\lambda
(X_s)\,d{\tilde W}_s>\frac{\epsilon T}{3}\biggr).
\end{eqnarray*}
Hence, we obtain the present lemma in the same way as Lemma~\ref{lem7.2}.
\end{pf}

For $\kappa<0$, we shall prove Theorem~\ref{theo2.4}. By convexity, we have
\[
\chi(-1)\geq\chi(\gamma)+\chi'(\gamma)(-1-\gamma),\qquad\gamma<-1.
\]
That is,
\[
\chi(\gamma)-\gamma\kappa\leq\chi(-1)+\chi'(\gamma)+\gamma
\bigl(\chi'(\gamma)-\kappa\bigr).
\]
$\chi'(\gamma)$ is monotonically nondecreasing, and $\chi'(\gamma
)\to0$ as $\gamma\to-\infty$. Therefore, we see that
\[
\chi(\gamma)-\gamma\kappa\to-\infty \qquad \mbox{as }   \gamma
\to-\infty.
\]
Hence,
\[
\inf_{\gamma<0}\{\chi(\gamma)-\gamma\kappa\}=-\infty.
\]
On the other hand, by taking $h=0$, we have $V_T(h)=v\exp(rT)$ and
\[
P\biggl(\frac{\log V_T(h)-\log S_T^0}{T}\leq\kappa\biggr)=0
\]
for sufficiently large $T$. Thus, $J(\kappa)=-\infty$.

%
\section{\texorpdfstring{Proof of Theorem~\protect\ref{theo2.5}}{Proof of Theorem 2.5}}\label{sec8}
For a given constant $0<\kappa<\chi'(0-)$, take $\gamma(\kappa)$
such that $\chi'(\gamma(\kappa))=\kappa$, namely,
\[
\inf_{\gamma<0}\{\chi(\gamma)-\gamma\kappa\}=\chi(\gamma(\kappa
))-\gamma(\kappa)\kappa.
\]
Then, since
\[
\inf_{h_{\cdot}}\log P\biggl(\frac{\log V_T(h)-\log S_T^0}{T}\leq\kappa\biggr)\leq
\inf_{h_{\cdot}}\log E\biggl[\biggl(\frac{V_T(h)}{S_T^0}\biggr)^{\gamma(\kappa)}\biggr]-\gamma
(\kappa)\kappa T,
\]
we have
\[
J(\kappa)\leq\varliminf_{T\to\infty}\frac{1}{T}\inf_{h_{\cdot}}\log
E\biggl[\biggl(\frac{V_T(h)}{S_T^0}\biggr)^{\gamma(\kappa)}\biggr]-\gamma(\kappa)\kappa.
\]
Therefore, if we prove that
%
\begin{eqnarray}\label{eq8.1}
 \qquad \varliminf_{T\to\infty}\frac{1}{T}\inf_{h_{\cdot}}\log E\biggl[\biggl(\frac
{V_T(h)}{S_T^0}\biggr)^{\gamma(\kappa)}\biggr]&=&\lim_{T\to\infty}\frac
{1}{T}\log E\biggl[\biggl(\frac{V_T(h^{(\gamma(\kappa))})}{S_T^0}\biggr)^{\gamma
(\kappa)}\biggr]\nonumber
\\[-8pt]
\\[-8pt] \quad &=&\chi(\gamma(\kappa)),
\nonumber
\end{eqnarray}
then we complete the proof of the present theorem because
\[
J(\kappa)\leq J_{\infty}(\kappa)\leq\lim_{T\to\infty}\frac
{1}{T}\log E\biggl[\biggl(\frac{V_T(h^{(\gamma(\kappa))})}{S_T^0}\biggr)^{\gamma
(\kappa)}\biggr]-\gamma(\kappa)\kappa
\]
and $J(\kappa)=\inf_{\gamma<0}\{\chi(\gamma)-\gamma\kappa\}$ by
Theorem~\ref{theo2.4}.
(\ref{eq8.1}) is proved in the following proposition.
%
\begin{prop}
Under the assumptions of Theorem~\ref{theo2.5}, (\ref{eq8.1}) holds.
\end{prop}

\begin{pf} Let $w=w^{(\gamma(\kappa))}$ be a solution to (\ref{eq2.19}) for
$\gamma=\gamma(\kappa)$ and ${\bar h}_t^{(\gamma)}={\bar h}(X_t)$,
where $X_t$ is the solution to (\ref{eq2.25}). Noting that
%
\begin{eqnarray}\label{eq8.2}
\eta(x,{\bar h})&=&{\bar h}^*{\hat\alpha}-\frac{1-\gamma}{2}{\bar
h}^*\sigma\sigma^*{\bar h} \nonumber
\\[-8pt]
\\[-8pt]
&=&\frac{1}{2(1-\gamma)}{\hat\alpha}^*\sigma\sigma^*{\hat\alpha
}-\frac{1}{2(1-\gamma)}(Dw)^*\lambda\sigma^*(\sigma\sigma
^*)^{-1}\sigma\lambda^*Dw,
\nonumber
\end{eqnarray}
we have
\begin{eqnarray*}
 &&w(X_t)-w(X_0)\\
 && \qquad =\int_0^t\biggl\{\frac{1}{2}\operatorname{tr}[\lambda\lambda
^*D^2w]+(\beta+\gamma\lambda\sigma^*{\bar h})^*Dw\biggr\}(X_s)\,ds \\
&& \qquad  \quad {} +\int_0^t(Dw)^*\lambda(X_s)\,dW_s^{{\bar h}} \\
&& \qquad =\int_0^t\biggl\{\frac{1}{2}\operatorname{tr}[\lambda\lambda^*D^2w]+\beta
_{\gamma}^*Dw\\
&& \qquad\hphantom{=\int_0^t\biggl\{} {}+\frac{\gamma}{1-\gamma}(Dw)^*\lambda\sigma^*(\sigma
\sigma^*)^{-1}\sigma\lambda^*Dw\biggr\}(X_s)\,ds \\
&& \qquad  \quad {} +\int_0^t(Dw)^*\lambda(X_s)\,dW_s^{{\bar h}} \\
&& \qquad =\int_0^t\biggl\{\chi+U_{\gamma}+\frac{\gamma}{2(1-\gamma
)}(Dw)^*\lambda\sigma^*(\sigma\sigma^*)^{-1}\sigma\lambda
^*Dw\\
&&\hspace*{179.5pt}{}-\frac{1}{2}(Dw)^*\lambda\lambda^*Dw\biggr\}(X_s)\,ds \\
&& \qquad  \quad {} +\int_0^t(Dw)^*\lambda(X_s)\,dW_s^{{\bar h}} \\
&& \qquad =\int_0^t\biggl\{\chi-\gamma\eta\bigl(X_s,{\bar h}_s^{(\gamma)}\bigr)-\frac
{1}{2}(Dw)^*\lambda\lambda^*Dw(X_s)\biggr\}\,ds\\
&& \qquad  \quad {}+\int_0^t(Dw)^*\lambda
(X_s)\,dW_s^{{\bar h}}.
\end{eqnarray*}
Thus,
\begin{eqnarray*}
&&E^{{\bar h}}\bigl[e^{\int_0^T\gamma\eta(X_s,{\bar h}_s^{(\gamma
)})\,ds}\bigr] \\
&& \qquad =E^{{\bar h}}\bigl[e^{\chi T+w(x)-w(X_T)+\int_0^T(Dw)^*\lambda
(X_s)\,dW_s^{\bar h}-\fracd{1}{2}\int_0^T(Dw)^*\lambda\lambda^*Dw(X_s)\,ds}\bigr].
\end{eqnarray*}
Let us introduce a new measure ${\check P}$ defined by
\[
\frac{d{\check P}}{d P^{\bar h}}=e^{\int_0^T(Dw)^*\lambda
(X_s)\,dW_s^{\bar h}-\fracd{1}{2}\int_0^T(Dw)^*\lambda\lambda^*Dw(X_s)\,ds}.
\]
Then
\[
{\check W}_t=W_t^{\bar h}-\int_0^t\lambda^*Dw(X_s)\,ds
\]
is a Brownian motion process under ${\check P}$ and
\[
dX_t=\{\beta(X_t)+\gamma\lambda\sigma^*{\bar h}(X_t)+\lambda
\lambda^*Dw(X_t)\}\,dt+\lambda(X_t)\,d{\check W}_t.
\]
Therefore,
\begin{eqnarray*}
 &&e^{-w(X_T)}-e^{-w(x)}\\[-0.5pt]
 && \qquad =-\int_0^Te^{-w(X_s)}(Dw)^*\lambda
(X_s)\,d{\check W}_s \\[-0.5pt]
&& \qquad  \quad {}+\int_0^Te^{-w(X_s)}\biggl\{-\frac{1}{2}\operatorname{tr}[\lambda\lambda
^*D^2w]-(\beta+\gamma\lambda\sigma^*{\bar h})^*Dw\\[-0.5pt]
&&\hspace*{183pt}{}-\frac
{1}{2}(Dw)^*\lambda\lambda^*Dw\biggr\}(X_s)\,ds \\[-0.5pt]
&& \qquad =-\int_0^Te^{-w(X_s)}(Dw)^*\lambda(X_s)\,d{\check W}_s \\[-0.5pt]
&& \qquad  \quad {}-\int_0^Te^{-w(X_s)}\biggl\{\frac{1}{2}\operatorname{tr}[\lambda\lambda
^*D^2w]+\beta_{\gamma}^*Dw\\[-0.5pt]
&&\hphantom{\qquad  \quad {}-\int_0^Te^{-w(X_s)}\biggl\{}{}+\frac{\gamma}{1-\gamma}(Dw)^*\lambda
\sigma^*(\sigma\sigma^*)^{-1}\sigma\lambda^*Dw \\[-0.5pt]
&&\hspace*{149pt} \qquad  \quad {} +\frac{1}{2}(Dw)^*\lambda\lambda^*Dw\biggr\}(x_s)\,ds \\[-0.5pt]
&& \qquad =-\int_0^Te^{-w(X_s)}(Dw)^*\lambda(X_s)\,d{\check W}_s-\int
_0^Te^{-w(X_s)}\bigl\{\chi-\gamma\eta\bigl(X_s,{\bar h}_s^{(\gamma)}\bigr)\bigr\}\,ds.
\end{eqnarray*}
Then, by the arguments using the stopping time, we have
\begin{eqnarray*}
&&
E^{{\bar h}}\bigl[e^{\int_0^T\gamma\eta(X_s,{\bar h}_s^{(\gamma
)})\,ds}\bigr]\\[-0.5pt]
&&\qquad=e^{\chi T+w(x)}{\check E}\bigl[e^{-w(X_T)}\bigr] \\[-0.5pt]
&&\qquad=e^{\chi T+w(x)}{\check
E}\biggl[e^{-w(x)}+\int_0^Te^{-w(X_s)}\bigl\{\gamma \eta\bigl(X_s,{\bar
h}_s^{(\gamma)}\bigr)-\chi\bigr\}\,ds\biggr].
\end{eqnarray*}
Hence, we obtain
\[
\lim_{T\to\infty}\frac{1}{T}\log E^{{\bar h}}\bigl[e^{\int_0^T\gamma
\eta(X_s,{\bar h}_s^{(\gamma)})\,ds}\bigr]\leq\chi(\gamma)
\]
by taking into account (\ref{eq2.26}) and (\ref{eq8.2}). The converse inequality holds
since ${\bar h}_s^{(\gamma)}\in{\mathcal A}(T)$.
\end{pf}

\setcounter{equation}{0}
\setcounter{prop}{0}
\begin{appendix}
\section*{Appendix}\label{appm}
Let $L_0$ be an elliptic operator defined by
%
\begin{eqnarray}\label{eq9.1}
L_0u&:=&\frac{1}{2}\sum_{i,j}a^{ij}(x)D_{ij}u+\sum_{i}b^i(x)D_iu\nonumber
\\[-8pt]
\\[-8pt]&\hspace*{2.3pt}=&\frac
{1}{2}\sum_{i,j}D_i(a^{ij}(x)D_ju)+\sum_i{\tilde b}^i(x)D_iu,
\nonumber
\end{eqnarray}
where $a^{ij}(x)$ and $b^i(x)$ are Lipschitz continuous functions such that
%
\begin{equation}\label{eq9.2}
k_0|y|^2\leq y^*a(x)y\leq k_1|y|^2 \qquad \forall y\in R^N,  k_0,
k_1>0
\end{equation}
and ${\tilde b}^i=b^i-\frac{1}{2}\sum_jD_ja^{ji}$. We assume that
there exists a positive function $\psi\in C^2(R^N)$ such that
%
\begin{equation}\label{eq9.3}
\cases{\displaystyle
\psi(x)\to\infty,&\quad as    $|x|\to\infty$, \vspace*{2pt}\cr\displaystyle
  -L_0\psi-\frac{c_a}{\psi}(D\psi)^*aD\psi> 0,&\quad $x\in
B_{R_0}^c,  \exists R_0>0, c_a>0$, \vspace*{2pt}\cr\displaystyle
L_0\psi<-1,&\quad $x\in B_{R_0}^c$.
}
\end{equation}
Set $K(x;\psi)=-L_0 \psi$,
\begin{eqnarray*}
F_{\psi}&=&\biggl\{u\in W^{2,p}_{\mathrm{loc}};\esssup\limits_{x\in B_{R_0}^c}\frac
{|u(x)|}{\psi(x)}<\infty\biggr\},\\
  F_K&=&\biggl\{f\in L^{\infty}_{\mathrm{loc}}; \esssup\limits_{x\in B_{R_0}^c}\frac
{|f(x)|}{K(x;\psi)}<\infty\biggr\}
\end{eqnarray*}
and
\[
{\mathcal D}=B_{R_0}=\{x\in R^N; |x|< R_0\}.
\]
Then we consider the following exterior Dirichlet problem for a given
bounded continuous function $h$ on $\Gamma=\partial{\mathcal D}$:
%
\begin{equation}\label{eq9.4}
\cases{\displaystyle
-L_0\xi=0,& \quad $x\in{\overline{\mathcal D}}^c$,\cr\displaystyle
  \xi |_{\Gamma}=h.
}
\end{equation}

\begin{prop}\label{prop9.1}
Exterior Dirichlet problem (\ref{eq9.4}) has a unique bounded solution $\xi\in
W^{2,p}_{\mathrm{loc}}\cap L^{\infty}$, $1<p<\infty$.
\end{prop}

\begin{pf} We first show uniqueness. Note that
\[
-L_0\psi=K(x;\psi)>0,\qquad x\in{\mathcal D}^c,
\]
and set $\xi=\mu\psi$. Then
\[
0=L_0\xi=(-L_0\mu)\psi-(L_0\psi)\mu-(D\mu)^*aD\psi.
\]
Therefore, $\mu$ satisfies
%
\begin{equation}\label{eq9.5}
\cases{\displaystyle
 -L_0\mu-\biggl(\frac{D\psi}\psi\biggr)^*aD\mu-\frac{L_0\psi}\psi\mu=0,\cr\displaystyle
  \mu |_{\Gamma}=\frac{h}\psi \quad    \mbox{and} \quad   \mu
(x)\to0   \qquad \mbox{as }  |x|\to\infty.
}
\end{equation}
Let $\mu_1$ and $\mu_2$ be solutions to (\ref{eq9.5}). Then $g:=\mu_1-\mu
_2$ satisfies
%
\begin{equation}\label{eq9.6}
\cases{\displaystyle
-L_0g-\biggl(\frac{D\psi}{\psi}\biggr)^*a Dg-\frac{L_0\psi}{\psi} g=0,\vspace*{2pt}\cr\displaystyle
  g |_{\Gamma}=0 \quad    \mbox{and} \quad   g(x)\to0 \qquad   \mbox
{as }  |x|\to\infty.
}
\end{equation}
To prove uniqueness, it is sufficient to show that the solution $g$ to
(\ref{eq9.6}) is trivial. For each $\epsilon>0$, there exists $R_{\epsilon}$
such that $|g|\leq\epsilon,   B_{R_{\epsilon}}^c$. Take $R\geq
R_{\epsilon}\vee R_0$. Then we see that
\[
|g|\leq\epsilon,\qquad B_R\cap{\mathcal D}^c,
\]
since $\psi> 0$ and $K(x;\psi)>0$ in ${\mathcal D}^c$. Thus, we see
that $g=0$ because $\epsilon$ is arbitrary.

Let us show the existence of the solution to (\ref{eq9.4}). We can assume
$h\geq0$. Consider the following Dirichlet problem for $R>R_0$:
%
\begin{equation}\label{eq9.7}
\cases{\displaystyle
-L_0\xi_R=0,& \quad  $B_R\cap{\overline{\mathcal D}}^c$, \cr\displaystyle
  \xi_R |_{\Gamma}=h,& \quad    $\xi_R |_{\partial B_R}=0$.
}
\end{equation}
Then we have
%
\begin{equation}\label{eq9.8}
\|\xi_R\|_{L^{\infty}(B_R\cap{\mathcal D}^c)}\leq
\| h\|_{L^{\infty}(\Gamma)}.
\end{equation}
It is clear that $\xi_R\leq\xi_{R'}$ and $R<R'$ by the maximum
principle. Therefore, there exists $\xi\in L^{\infty}(R^n\cap
{\mathcal D}^c)$ and
\[
\xi_R\to\xi, \qquad\|\xi\|_{L^{\infty}({\mathcal
D}^c)}\leq\| h\|_{L^{\infty}(\Gamma)}.
\]
When taking
\[
{\mathcal D}^*\subset\subset{\tilde{\mathcal D}}\subset B_R\cap
{\mathcal D}^c,
\]
we see that
\[
\|\xi_R\|_{W^{2,p}({\mathcal D}^*)}\leq c\|\xi
_R\|_{L^p({\tilde{\mathcal D}})}\leq c'\|\xi
_R\|_{L^{\infty}({\tilde{\mathcal D}})}\leq c'\|
h\|_{L^{\infty}(\Gamma)}.
\]
Thus, $\xi_R$ converges to $\xi$ weakly in $W^{1,q}_{\mathrm{loc}}$.
Regularity theorems show that $\xi\in W^{2,q}_{\mathrm{loc}}$.
\end{pf}

Let us take a bounded domain ${\mathcal D}_1$ such that ${\mathcal
D}\subset{\mathcal D}_1$ and a bounded Borel function $\phi$ on
$\Gamma_1=\partial{\mathcal D}_1$. We consider a Dirichlet problem
%
\begin{equation}\label{eq9.10}
\cases{\displaystyle
-L_0\zeta=0,& \quad ${\mathcal D}_1$, \cr\displaystyle
  \zeta |_{\Gamma_1}=\phi,
}
\end{equation}
which admits a solution $\zeta\in W^{2,p}({\mathcal D}_1)\cap
L^{\infty}$. For this solution, we consider exterior Dirichlet problem
(\ref{eq9.4}) with $h=\zeta$. Then we introduce an operator
$
P\dvtx  {\mathbf B}(\Gamma_1)\mapsto{\mathbf B}(\Gamma_1)
$
defined by
%
\begin{equation}\label{eq9.11}
P\phi(x)=\xi(x), \qquad   x\in\Gamma_1,
\end{equation}
where $\xi(x)$ is the solution to (\ref{eq9.4}) with $h=\zeta$. In a similar
manner to Lem\-ma~5.1 in   \cite{Bens0}, Chapter II, we have
\[
\sup_{B\in{\mathcal B}(\Gamma_1),  x,y\in\Gamma_1}\lambda_{x,y}(B)<1,
\]
where
\[
\lambda_{x,y}(B)=P\chi_B(x)-P\chi_B(y),\qquad B\in{\mathcal
B}(\Gamma_1).
\]
Moreover, we have the following proposition; cf. Theorem 4.1, Chapter
II in~\cite{Bens0}.
%
\begin{prop} Operator $P$ defined above satisfies the following properties.
%
\begin{equation}\label{eq9.12}
\| P\phi\|_{L^{\infty}(\Gamma_1)}\leq\|\phi
\|_{L^{\infty}(\Gamma_1)},\qquad P1(x)=1,
\end{equation}
and for some $\delta>0$,
%
\begin{equation}\label{eq9.13}
P\chi_B(x)-P\chi_B(y)\leq1-\delta,  \qquad   x,y\in\Gamma_1, B\in
{\mathcal B}(\Gamma_1).
\end{equation}
Furthermore, there exists a probability measure $\pi(dx)$ on $(\Gamma
_1,{\mathcal B}(\Gamma_1))$ such that
%
\begin{eqnarray}\label{eq9.14}
   \biggl|P^n\phi(x)-\int\phi(x)\pi(dx)\biggr|\leq K\|\phi\|
_{L^{\infty}}e^{-\rho n}, \nonumber
\\[-8pt]
\\[-8pt]
 \eqntext{\displaystyle  \rho=\log\frac{1}{1-\delta},  K=\frac{2}{1-\delta},}
\end{eqnarray}
and
%
\begin{equation}\label{eq9.15}
\int\phi(x)\pi(dx)=\int P\phi(x)\pi(dx)
\end{equation}
for any bounded Borel function $\phi$.
\end{prop}

Consider an exterior Dirichlet problem for a given function $f\in F_K$:
%
\begin{equation}\label{eq9.16}
\cases{\displaystyle
-L_0u=f,& \quad  $x\in{\overline{\mathcal D}}^c$, \cr\displaystyle
  u |_{\Gamma}=0.
}
\end{equation}
Then we have the following proposition.
%
\begin{prop}
For a given function $f\in F_K$, there exists a unique solution $u\in
W^{2,p}_{\mathrm{loc}}$, $1<p<\infty$, to (\ref{eq9.16}) such that
\[
\sup_{x\in{\mathcal D}^c}\frac{|u(x)|}{\psi(x)}<\infty.
\]
\end{prop}

\begin{pf} Assume that $f \geq0$, $f\in F_K$. For $R>R_0$ we
consider a Dirichlet problem on $B_R\cap{\mathcal D}^c$:
%
\begin{equation}\label{eq9.17}
\cases{\displaystyle
-L_0u_R=f, \qquad  x\in B_R\cap{\overline{\mathcal D}}^c, \cr\displaystyle
  u_R |_{\Gamma}=0, \qquad   u_R |_{\partial B_R}=0.
}
\end{equation}
There exists a unique solution $u_R\in W_0^{2,p}(B_R\cap{\overline
{\mathcal D}}^c)$. Set
\[
c_f=\esssup\limits_{x\in{\mathcal D}^c}\frac{|f(x)|}{K(x;\psi)}.
\]
Then we have
%
\begin{equation}\label{eq9.18}
0\leq u_R\leq c_f\psi.
\end{equation}
To see that, set ${\tilde u}_R:=u_R-c_f\psi$. Then
\begin{eqnarray*}
-L_0{\tilde u}_R&=&-L_0u_R+c_fL_0\psi \\
&=&f-c_fK(x;\psi)\leq0,\qquad{\overline{\mathcal D}}^c\cap B_R.
\end{eqnarray*}
Therefore,
\[
{\tilde u}_R\leq0, \qquad{\overline{\mathcal D}}^c\cap B_R,
\]
since ${\tilde u}_R$ is subharmonic in ${\overline{\mathcal D}}^c\cap
B_R$ and ${\tilde u}_R\leq0$ on $\Gamma\cup\partial B_R$. Hence, we
have $u_R\leq c_f\psi$.

On the other hand,
\[
-L_0u_R=f\geq0, \qquad B_R\cap{\overline{\mathcal D}}^c.
\]
Hence, $u_R$ is superharmonic and $u_R=0$ on $\Gamma\cap\partial
B_R$. Thus,
\[
u_R\geq0, \qquad B_R\cap{\overline{\mathcal D}}^c,
\]
and (\ref{eq9.18}) holds.

If $f\leq0,  f\in F_K$ and $c_f=\esssup_{x\in{\mathcal D}^c}\frac
{|f(x)|}{K(x;\psi)}$, then, through the same arguments for $-f$, we obtain
\[
-c_f\psi\leq u_R^-\leq0,
\]
where $-u_R^-$ is the corresponding solution to (\ref{eq9.17}). Therefore, for
general $f=f^+-f^-$, we have
\[
-c_f\psi\leq u_R\leq c_f \psi.
\]

 Let $u_R^+$ be a solution to (\ref{eq9.17}) for $f^+$. Then, $u_R^+$ is
nondecreasing with respect to $R$ because of the maximum principle.
Indeed, for $R<R'$, we have
\begin{eqnarray*}
 -L_0(u_{R'}^+-u_R^+)&=&0, \qquad B_R\cap{\overline{\mathcal
D}}^c, \\
 u_{R'}^+-u_R^+&\geq&0, \qquad\Gamma\cup\partial B_R.
\end{eqnarray*}
Since $u_R^+$ is dominated by $c_f\psi$, there exists $u^+$ such that
\[
u^+(x)=\lim_{R\to\infty}u_R^+(x),\qquad  u^+(x)
|_{\Gamma}=0.
\]
Let us show that $u^+(x)$ satisfies
\[
-L_0u^+=f^+, \qquad{\overline{\mathcal D}}^c.
\]
Set
\[
{\mathcal D}^*:=B_R\cap{\overline{\mathcal D}}^c,\qquad{\mathcal
D}'\subset\subset{\mathcal D}^*\cup\partial{\mathcal D}^*.
\]
Then we have
\begin{eqnarray*}
\| u^+\|_{2,p;{\mathcal D}'}&\leq&c(\|
u^+\|_{p;{\mathcal D}^*}+\| f^+\|_{p;{\mathcal
D}^*}) \\
&\leq&c'(\| u^+\|_{\infty;{\mathcal D}^*}+\|
f^+\|_{p;{\mathcal D}^*}.
\end{eqnarray*}
For ${\mathcal D}''\subset{\mathcal D}'$ injection
$W^{2,p}({\mathcal D}')\hookrightarrow W^{1,q}({\mathcal D}''),
1\leq q\leq\frac{np}{n-p}$, is\vspace*{-2pt} compact. Therefore,
$u_R^+\to u^+$ weakly in $W^{1,q}_{\mathrm{loc}}$ for each $1\leq q<\infty$ and
$u^+$ is a weak solution to
\[
\cases{\displaystyle
 -L_0u^+=f^+, &\quad $R^n\cap{\overline{\mathcal D}}^c$,  \cr\displaystyle
   u^+ |_{\Gamma}=0.
}
\]
By the regularity theorem $u^+\in W^{2,p}_{\mathrm{loc}}$, $\forall p>1$.

Similarly, we have $u^-\in W^{2,p}_{\mathrm{loc}}$, which is a solution to
\[
\cases{\displaystyle
 -L_0u^-=f^-, &\quad $R^n\cap{\overline{\mathcal D}}^c$,  \cr\displaystyle
    u^- |_{\Gamma}=0.
}
\]

Now let us prove uniqueness. For $i=1,2$, we assume that $u_i$ is a
solution to (\ref{eq9.16}) such that
\[
-c_f\psi\leq u_i\leq c_f \psi, \qquad u_i\in W^{2,p}_{\mathrm{loc}}.
\]
Then $u=u_1-u_2$ satisfies
%
\begin{equation}\label{eq9.19}
\cases{\displaystyle
 -L_0u=0,&\quad${\overline{\mathcal D}}^c$,  \cr\displaystyle
    u |_{\Gamma}=0,& \quad   $-2c_f\psi\leq u\leq2c_f\psi,
  u\in W^{2,p}_{\mathrm{loc}}$.
}
\end{equation}
We shall prove that $u$ satisfying (\ref{eq9.19}) is trivial, $u\equiv0$. For
this purpose, we set
\[
u=v\psi^{\alpha} ,\qquad\alpha=1+c_a>1,
\]
where $c_a>0$ is the constant that appears in (\ref{eq9.3}).
Since $-L_0u=0$ we have
\[
L_0v+2a\biggl(\frac{D\psi}{\psi}\biggr)^*a Dv+\frac{\alpha v}\psi\biggl\{L_0\psi+
\frac{\alpha-1}{\psi}(D\psi)^*a D\psi\biggr\}=0.
\]
Note that
\[
-L_0\psi-\frac{\alpha-1}{\psi}(D\psi)^*a D\psi=
K(x;\psi)-\frac{\alpha-1}{\psi}(D\psi)^*a D\psi\geq0
\]
for $|x|\gg1$ under assumption (\ref{eq9.3}). Moreover,
\[
v=\frac{u}{\psi^{\alpha}}\to0\qquad\mbox{as }  |x|\to\infty.
\]
Hence, from the maximum principle, we see that $v\equiv0$ as in the
proof of Proposition~\ref{prop9.1}.
\end{pf}

Let $f$ be a function on $R^n$ such that $f$ is bounded in ${\mathcal
D}$ and $f\in F_K({\mathcal D}^c)$, and ${\mathcal D}_1$ a bounded
domain such that ${\mathcal D}\subset{\mathcal D}_1$. We consider
%
\begin{equation}\label{eq9.21}
\cases{\displaystyle
-L_0\Psi=f, & \quad  ${\mathcal D}_1$,  \cr\displaystyle
  \Psi |_{\Gamma_1}=0,
}
\end{equation}
and
%
\begin{equation}\label{eq9.22}
\cases{\displaystyle
-L_0\xi=f, & \quad  $R^n\cap{\overline{\mathcal D}}^c$,  \cr\displaystyle
  \xi |_{\Gamma}= \Psi |_{\Gamma}.
}
\end{equation}
Then we set
\[
Tf(x)=\xi(x),\qquad x\in \Gamma_1,
\]
and
%
\begin{equation}\label{eq9.23}
\nu(f)=\frac{\int_{\Gamma_1}Tf(\sigma)\pi(d\sigma)}{\int
_{\Gamma_1}T1(\sigma)\pi(d\sigma)}.
\end{equation}
We further consider
%
\begin{equation}\label{eq9.24}
\cases{\displaystyle
 -L_0z=f, \cr\displaystyle
  z\in W^{2,p}_{\mathrm{loc}},& \quad $\displaystyle     \sup_{x\in{\mathcal D}^c}\frac{|z|}{\psi
}<\infty$.
}
\end{equation}
Then, as in the proof of Theorem 5.3, in \cite{Bens0}, Chapter II, we
obtain the following proposition. Here, we only give the proof of the
existence of the solution for use in Section~\ref{sec6}.
%
\begin{prop}\label{prop9.4}
Equation (\ref{eq9.24}) has a solution unique up to additive constants if and only if
$\nu(f)=0$. Moreover,
%
\begin{equation}\label{eq9.25}
\nu(f)=\int m(y)f(y)\,dy
\end{equation}
for $m\in L^1(R^n)$, $m\geq0$ and $-L_0^*m=0$ in distribution sense
%
\begin{equation}\label{eq9.26}
\int m(y)(-L_0z)\,dy=0,\qquad z\in W^{2,p}_{\mathrm{loc}},
\end{equation}
such that $z\in F_\psi$ and $-L_0z\in F_K$. Furthermore, $m(x)$ is the
only function in $L^1$ satisfying (\ref{eq9.26}) and
\[
\int m(x)\,dx =1.
\]
\end{prop}

\begin{pf*}{Proof of existence} Let $\zeta_0=\Psi$ and $ \eta_0=\xi$,
where $\Psi$ (resp., $\xi$) is the solution to (\ref{eq9.21}) [resp., (\ref{eq9.22})].
For each $k=1,2,\ldots$ define
$\zeta_k$ and $\eta_k$ as follows. Let $\zeta_k$ be the solution to
(\ref{eq9.10}) for $\phi=\eta_{k-1}$, and $\eta_k$ the solution to~(\ref{eq9.4}) for
$h=\zeta_k$. Then
\[
  \eta_0(x) |_{\Gamma_1}=Tf(x),\qquad  \eta
_n(x) |_{\Gamma_1}=P^n(Tf)(x),  \qquad   n=1,2,\ldots.
\]
Since $\int_{\Gamma_1}Tf(x)\pi(dx)=0$, we have
\[
|P^n(Tf)(x)|\leq K\| Tf\|_{L^{\infty}(\Gamma
_1)}e^{-\rho n}
\]
by (\ref{eq9.14}). Set ${\tilde\eta}_n(x)=\sum_{k=0}^n\eta_k(x), {\tilde
\zeta}_n(x)=\sum_{k=0}^n\zeta_k(x)$. Then
\[
  {\tilde\eta}_n |_{\Gamma_1}=Tf+P(Tf)+\cdots+P^n(Tf).
\]
Therefore, we see that there exists ${\bar\eta}\in C(\Gamma_1)$ such
that $\|{\tilde\eta}_n-{\bar\eta}\|_{L^{\infty
}(\Gamma_1)}\to0$, \mbox{$n\to0$}.
Moreover, we have
\[
\|{\tilde\eta}_n\|_{L^{\infty}(\Gamma_1)}\leq
K\| Tf\|_{L^{\infty}(\Gamma_1)}\frac{1}{1-e^{-\rho}}.
\]
Note that ${\tilde\zeta}_n$ is the solution to
%
\begin{equation}\label{eq9.27}
\cases{\displaystyle
-L_0{\tilde\zeta}_n=f,  \qquad   {\mathcal D}_1, \vspace*{2pt}\cr\displaystyle
  {\tilde\zeta_n} |_{\Gamma_1}=  {\tilde\eta
}_{n-1} |_{\Gamma_1},
}
\end{equation}
and ${\tilde\eta}_n$ the solution to (\ref{eq9.22}) with $\Psi(x)={\tilde
\zeta}_n$. Noting that
\[
\|{\tilde\zeta}_n-{\tilde\zeta}_m\|_{L^{\infty
}({\mathcal D}_1)}\leq\|{\tilde\zeta}_n-{\tilde\zeta
}_m\|_{L^{\infty}(\Gamma_1)}\leq\|{\tilde\eta
}_{n-1}-{\tilde\eta}_{m-1}\|_{L^{\infty}(\Gamma_1)},
\]
we see that ${\tilde\zeta}_n$ converges in $C({\overline{\mathcal
D}_1})$ and weakly in $W^{1,q}_{\mathrm{loc}}$ since its $W^{2,p}_{\mathrm{loc}}$ norm is
bounded; cf. Theorem 9.11 in \cite{GT}. By the regularity theorems, the
limit ${\bar\zeta}\in W^{2,p}_{\mathrm{loc}}\cap C({\overline{\mathcal
D}_1})$ and satisfies
%
\begin{equation}\label{eq9.28}
\cases{\displaystyle
-L_0{\tilde\zeta}=f, & \quad  ${\mathcal D}_1$,\vspace*{2pt} \cr\displaystyle
  {\tilde\zeta} |_{\Gamma_1}={\bar\eta}.
}
\end{equation}
On the other hand, ${\tilde\eta}_n-\xi$ is the solution to (\ref{eq9.4})
with $h= {\tilde\zeta}_n-\Psi |_{\Gamma}=  \sum
_{i=1}^n\zeta_i |_{\Gamma}$ and
\[
\|{\tilde\eta}_n-\xi\|_{L^{\infty}(\Gamma)}\leq
\Biggl\|\sum_{j=1}^n\zeta_j\Biggr\|_{L^{\infty}(\Gamma)}\leq
\Biggl\|\sum_{j=0}^{n-1}\eta_j\Biggr\|_{L^{\infty}(\Gamma
_1)}\leq\|{\tilde\eta}_{n-1}\|_{L^{\infty}(\Gamma_1)}.
\]
Thus, we see that
\[
\|{\tilde\eta}_n-\xi\|_{L^{\infty}({\mathcal
D}^c)}\leq\|{\tilde\eta}_{n-1}\|_{L^{\infty}(\Gamma_1)},
\]
and ${\tilde\eta}_n$ converges in $C({\mathcal D}^c)$ and weakly in
$W^{1,q}_{\mathrm{loc}}({\mathcal D}^c)$. The limit ${\tilde\eta}\in
W^{2,p}_{\mathrm{loc}}({\mathcal D}^c)\cap C( {\mathcal D}^c)$ and satisfies
(\ref{eq9.22}) with $\Psi={\tilde\zeta}$.
Setting $z={\tilde\zeta}$ in ${\mathcal D}_1$ and $z={\tilde\eta}$
in ${\mathcal D}^c$, we have a solution to (\ref{eq9.24}).
\end{pf*}
\end{appendix}

\section*{Acknowledgments}
 The author thanks the anonymous referee for
careful reading the manuscript and giving valuable suggestions. The
author also thanks Dr. Y.~Watanabe for fruitful discussions.


%

\printaddresses

\end{document}